\documentclass[12pt]{article}

\usepackage[colorlinks=true, pdfstartview=FitV, linkcolor=blue, citecolor=blue, urlcolor=blue]{hyperref}

\usepackage{amssymb,amsmath, amscd}
\usepackage{times, verbatim}
\usepackage{graphicx}
\input xy
\xyoption{all}

\newtheorem{thm}{Theorem}[section]

\newtheorem{prop}[thm]{Proposition}

\newtheorem{lemma}[thm]{Lemma}
\newtheorem{definition}[thm]{Definition}
\newtheorem{cor}[thm]{Corollary}

\newcommand\proof
{\par\noindent{\bf Proof:\ } }
\newcommand\qed{\hfill$\blacksquare$}
\newcommand{\half}{{\textstyle{ \frac{1}{2} } }}

\newcommand\aff{\mathsf{\!aff}}
\newcommand\cox{\mathsf{cox}}

\newcommand{\wtW}{\widetilde{W}}
        
\DeclareMathOperator{\Ad}{Ad}
\DeclareMathOperator{\Aut}{Aut}
\DeclareMathOperator{\ad}{ad}

\DeclareMathOperator{\End}{End}

\DeclareMathOperator{\Gal}{Gal}
\DeclareMathOperator{\GL}{GL}
\DeclareMathOperator{\Hom}{Hom}

\DeclareMathOperator{\Kac}{Kac}

\DeclareMathOperator{\lcm}{lcm}

\DeclareMathOperator{\tr}{tr}

\DeclareMathOperator{\PGL}{PGL}
\DeclareMathOperator{\rank}{rank}

\DeclareMathOperator{\SL}{SL}

\DeclareMathOperator{\Spec}{Spec}
\DeclareMathOperator{\Spin}{Spin}
\DeclareMathOperator{\Sp}{Sp}

\DeclareMathOperator{\Sym}{Sym}
\DeclareMathOperator{\sym}{sym}

\DeclareMathOperator{\un}{un}

\newcommand{\br}{\mathbb{R}}

\newcommand{\bc}{\mathbb{C}}

\newcommand{\bq}{\mathbb{Q}}
\newcommand{\bz}{\mathbb{Z}}

\newcommand{\al}{\alpha}
\newcommand{\be}{\beta}
\newcommand{\ga}{\gamma}

\newcommand{\De}{\Delta}
\newcommand{\lam}{\lambda}
\newcommand{\Lam}{\Lambda}
\newcommand{\om}{\omega}
\newcommand{\Om}{\Omega}

\newcommand{\vep}{\varepsilon}
\newcommand{\vp}{\varphi}

\newcommand{\Th}{\Theta}

\newcommand{\vt}{\vartheta}

\newcommand{\sca}{\mathcal{A}}

\newcommand{\scc}{\mathcal{C}}
\newcommand{\scd}{\mathcal{D}}

\newcommand{\scr}{\mathcal{R}}

\newcommand{\la}{\langle}
\newcommand{\ra}{\rangle}
\newcommand{\lra}{\longrightarrow}
\newcommand{\hra}{\hookrightarrow}

\newcommand{\bG}{\mathbf{G}}

\newcommand{\fc}{\mathfrak{c}}

\newcommand{\fe}{\mathfrak{e}}

\newcommand{\fg}{\mathfrak{g}}
\newcommand{\fh}{\mathfrak{h}}
\newcommand{\fk}{\mathfrak{k}}

\newcommand{\fl}{\mathfrak{l}}
\newcommand{\fm}{\mathfrak{m}}

\newcommand{\fs}{\mathfrak{s}}
\newcommand{\fsl}{\mathfrak{sl}}
\newcommand{\fso}{\mathfrak{so}}
\newcommand{\fsp}{\mathfrak{sp}}

\newcommand{\ft}{\mathfrak{t}}
\newcommand{\fu}{\mathfrak{u}}
\newcommand{\fv}{\mathfrak{v}}
\newcommand{\fz}{\mathfrak{z}}

\newcommand{\bfmu}{\boldsymbol{\mu}}

\newcommand{\twoAtwo}{
\SelectTips{eu}{}
\xymatrix@M=0pt{
\underset{\hphantom{1}}{\overset{1}{\bullet}}\ar@{=>}[r]\ar@<0.65ex>@{-}[r]+<-1.2ex,0ex>\ar@<-0.65ex>@{-}[r]+<-1.2ex,0ex>&
\underset{\hphantom{2}}{\overset{2}\circ}}
}

\newcommand{\twoAtwoo}{
\SelectTips{eu}{}
\xymatrix@M=0pt{
{1\ }\ar@{=>}[r]\ar@<0.65ex>@{-}[r]+<-1.2ex,0ex>\ar@<-0.65ex>@{-}[r]+<-1.2ex,0ex>&{\ 0}}
}
\newcommand{\twoAtwooo}{
\SelectTips{eu}{}
\xymatrix@M=0pt{
{1\ }\ar@{=>}[r]\ar@<0.65ex>@{-}[r]+<-1.2ex,0ex>\ar@<-0.65ex>@{-}[r]+<-1.2ex,0ex>&{\ 1}}
}

\newcommand{\outEVI}[5]{#1 #2 #3\!\Leftarrow\!#4 #5}
\newcommand{\FIV}[5]{#1 #2 #3\!\Rightarrow\!#4 #5}

\newcommand{\EVI}[7]{
{\scriptsize\begin{matrix}\!\!\!\!#2\!\!\!\!&\!\!\!\!#4\!\!\!\!&\!\!\!#5\!\!\!&\!\!\!\!#6\!\!\!\!
&\!\!\!\!#7
\!\!\!\!\\
&&\!\!\!#3\!\!\!&&\\&&\!\!\!#1\! \!\!&&\end{matrix}} }

\newcommand{\EVII}[8]{
{\scriptsize\begin{matrix}
\!\!\!#1\!\!\! \!&\!\!\!\!\!#2\!\!\!\!\!\!&\!\!\!\!\!#4\!\!\!\!\!&\!\!\!\!\!#5\!\!\!\!\!&\!\!\!\!\!#6\!\!\!\!\!&\!\!\!\!\!#7
\!\!\!\!\!&\!\!\!\!\!\!#8\!\!\!\\
&&&\!\!\!#3\!\!\!&&&\end{matrix}} }

\newcommand{\eVII}[8]{
{\begin{matrix}
 #1 & \!\!#2 \!& \!#3 \!\!& \!\!#4 \!\!& \!\!#5 \!& \!#6
 \!& \!#7\\
&&& #8 &&&\end{matrix}} }

\newcommand{\E}[9]{
{\scriptsize\begin{matrix}
\!\!\!#2\!\!\! \!\!&\!\!\!\!\!#4\!\!\!\!\!&\!\!\!\!\!#5\!\!\!\!\!&\!\!\!\!\!#6\!\!\!\!\!&\!\!\!\!\!#7\!\!\!\!\!&\!\!\!\!\!#8\!\!\!\!\!&\!\!\!\!\!#9
\!\!\!\!\!&\!\!\!\!\!#1\!\!\!\\
&&\!\!\!#3\!\!\!&&&&&\end{matrix}} }

\begin{document}

\title{Gradings of positive rank on simple Lie algebras}

\author{
Mark Reeder\thanks{Supported by NSF grants DMS-0801177 and DMS-0854909}\\
\texttt{reederma@bc.edu} 
\and
Paul Levy\\
\texttt{ p.d.levy@lancaster.ac.uk} 
\and
Jiu-Kang Yu\thanks{Supported by NSF grant DMS-0854909}\\
\texttt{jyu@math.purdue.edu} 
\and
Benedict H. Gross\thanks{Supported by NSF grant DMS-0901102}\\
\texttt{gross@math.harvard.edu} }

\date{\today}
\maketitle
\tableofcontents
\section{Introduction}
Let $\fg$ be the  Lie algebra of a connected simple algebraic group $G$ of adjoint type over an algebraically closed field $k$. 
A {\it grading} on $\fg$ is a decomposition 
$$\fg=\bigoplus_{i\in \bz/m}\fg_i$$
where $m$ is an integer $\geq 0$  and $[\fg_i,\fg_j]\subset \fg_{i+j}$ for all $i,j$. 
The summand $\fg_0$ is a Lie subalgebra of $\fg$ and we let $G_0$ denote the corresponding connected subgroup of $G$. 
The adjoint action of $G$ on $\fg$ restricts to an action of $G_0$  on each summand $\fg_i$. 
 We are interested in the invariant theory of this action, for which there is no loss of generality if we assume that $i=1$. 

If $m=1$ this is the invariant theory of the adjoint representation, first developed by Chevalley, 
who showed that the restriction $k[\fg]^G\to k[\ft]^W$ of $G$-invariant polynomials on $\fg$ 
to polynomials on a Cartan subalgebra $\ft$ invariant under the Weyl group $W$ is an isomorphism. 
This and other aspects of Chevalley's theory were generalized to the case $m=2$ by Kostant and Rallis \cite{kostant-rallis}. 
Soon after, Vinberg \cite{vinberg:graded} showed that for any $m\geq 0$ the invariant theory of the $G_0$-action on $\fg_1$ has similar parallels with the adjoint representation of $G$ on $\fg$. Vinberg worked over $\bc$, but in 
\cite{levy:thetap}, Vinberg's theory was extended to fields of good odd positive characteristic 
not dividing $m$. 

Some highlights of Vinberg theory are as follows. 
A {\it Cartan subspace} is a linear subspace $\fc\subset\fg_1$ 
which is abelian as a Lie algebra, consists of semisimple elements, and is maximal with these two properties.  All Cartan subspaces are conjugate under $G_0$. 
Hence the dimension of $\fc$ is an invariant of the grading, called the {\it rank}, which we denote in this introduction  by $r$. 
The {\it little Weyl group} is the subgroup $W_\fc$ of $\GL(\fc)$ arising from the action of the normalizer of $\fc$ in $G_0$. The group $W_\fc$ is finite and is generated by semisimple transformations of $\fc$ fixing a hyperplane and we have an isomorphism of invariant polynomial rings
$$k[\fg_1]^{G_0}\overset\sim\lra k[\fc]^{W_\fc},$$
given by restriction. Finally $k[\fg_1]^{G_0}\simeq k[f_1,\dots, f_r]$ is a polynomial algebra generated by $r$ algebraically independent polynomials  $f_1,\dots,f_r$ whose degrees $d_1,\dots, d_r$ are determined by the grading. In particular the product of these degrees is the order of $W_\fc$. 

We have a dichotomy: either the rank $r=0$, in which case $\fg_1$ consists entirely of nilpotent elements of $\fg$, or $r>0$, in which case $m>0$ and $\fg_1$ contains semisimple elements of $\fg$.  A basic problem is to classify all gradings of rank $r>0$ and to compute the little Weyl groups $W_\fc$ in each case.
Another open question is {\it Popov's conjecture:\ } $\fg_1$ should contain a
{\it Kostant section}: an affine subspace $\fv$ of $\fg_1$ with $\dim\fv=r$, such that the restriction map 
$k[\fg_1]^{G_0}\lra k[\fv]$
is an isomorphism. 
 
The classification of positive-rank gradings and their little Weyl groups, along with verification of Popov's conjecture was given in \cite{levy:thetap} and \cite{levy:exceptional} for gradings of  Lie algebras of classical  type and those of types $G_2$ and $F_4$. In this paper we complete this work by proving analogous  results for  types $E_6,\ E_7$ and $E_8$, using new methods which apply
to the Lie algebras of general simple algebraic groups $G$. 

The main idea is to compute Kac coordinates of lifts of automorphisms of the root system $R$ of $\fg$, as we shall now explain. 
Choosing a base in $R$ and a pinning in $\fg$ (defined in section \ref{mum}), we may write
the automorphism groups $\Aut(R)$ and $\Aut(\fg)$ as semidirect products:
$$\Aut(R)=W\rtimes \Theta,\qquad\Aut(\fg)=G\rtimes\Theta,$$
where $W$ is the Weyl group of $R$ and $\Theta$, the symmetry group of the Dynkin graph $D(R)$ of $R$, is identified with the group of automorphisms of $\fg$ fixing the chosen pinning. 
To each 
$\vt\in \Th$ one can associate an affine root system $\Psi=\Psi(R,\vt)$ consisting of affine functions on an affine space $\sca$ of dimension 
equal to the number of $\vt$-orbits on the nodes of the diagram $D(R)$. 
Kac' original construction of $\Psi$ uses infinite dimensional Lie algebras and works over $\bc$; our approach constructs $\Psi$ directly from the pair $(R,\vt)$ and works over any algebraically closed field in which the order $e$ of $\vt$ is nonzero. 
The choice of  pinning on $\fg$ determines a rational structure on $\sca$ and a basepoint $x_0\in\sca$. 
Following an idea of Serre \cite{serre:kac}, we associate to each rational point $x\in\sca_\bq$ an embedding $\varrho_x:\bfmu_m\hra G$ of group schemes over $k$, where $m$ is the denominator of $x$.
If $m$ is nonzero in $k$ and we choose a root of unity $\zeta\in k^\times$ of order $m$, 
then $x$ determines an actual automorphism $\theta_x\in G\vt$ of order $m$. If $x$ lies in the closure $\overline C$ of the fundamental alcove  of $\sca$ then the affine coordinates of $x$ are those defined by Kac (when $k=\bc$ and $\zeta=e^{2\pi i/m}$); we call these {\it normalized Kac coordinates}, since we also consider points $x$ outside $\overline C$ having some affine coordinates negative. Any $x\in\sca_\bq$ can be moved into $\overline C$ via operations of  the affine Weyl group $W(\Psi)$, and this can be done effectively, using a simple algorithm. 
See also \cite{levy:exceptional}, which gives a different way of extending Kac coordinates to positive characteristic. 

The  half-sum of the positive co-roots is a vector $\check\rho$ belonging to the translation subgroup of $\sca$. In the  {\it principal segment} $[x_0,x_0+\check\rho]\subset \sca$ we are especially interested in the points 
$$x_m:=x_0+\tfrac{1}{m}\check\rho\in\sca_\bq,$$ 
where $m$ is the order of an elliptic $\bz$-regular automorphism $\sigma\in\Aut(R)$. Here $\sigma$ is {\it elliptic} if $\sigma$ has no nonzero fixed-points in the reflection representation, and we say $\sigma$ is {\it $\bz$-regular} if the group generated by $\sigma$ acts freely on $R$. 
(This is almost equivalent to Springer's notion of regularity, and for our purposes it is the correct one. See section \ref{Zregular}.) 

Now assume that the characteristic of $k$ is not a torsion prime for $\fg$. 

Choose a Cartan subalgebra $\ft$ of $\fg$, let $T$ be the maximal torus of $G$ centralizing $\ft$ with normalizer $N$ in $G$ and let $\Aut(\fg,\ft)$ be the subgroup of $\Aut(\fg)$ preserving $\ft$.
The groups $\Aut(R)$ and $\Aut(\fg,\ft)/T$ are isomorphic and we may canonically identify $W$-conjugacy classes in $\Aut(R)$ with $N/T$-conjugacy classes in $\Aut(\fg,\ft)/T$. Let 
$\sigma\in\Aut(R)$ be an elliptic $\bz$-regular automorphism whose order $m$ is nonzero in $k$. Write $\sigma=w\cdot\vt$ with $w\in W$ and $\vt\in\Th$. 
Then there is a unique $G$-conjugacy class $C_\sigma\subset G\vt$ such that 
$C_\sigma\cap \Aut(\fg,\ft)$ projects to the class of $\sigma$ in $\Aut(R)$. 
Using results of Panyushev in \cite{panyushev:theta}, we show that $C_\sigma$ contains the automorphism $\theta_{x_m}$, where $x_m$ is the point
on the principal segment defined above. 
The  (un-normalized) Kac coordinates of $x_m$ are all $=1$ except one coordinate is $1+(m-h_\vt)/e$, where $h_\vt$ is the twisted Coxeter number of $(R,\vt)$. Translating by the affine Weyl group we obtain the normalized Kac coordinates of the class $C_\sigma\subset G\vt$. The automorphisms in $C_\sigma$ have positive rank equal to the multiplicity of the cyclotomic polynomial $\Phi_m$ in the characteristic polynomial of $\sigma$. They are exactly the semisimple automorphisms of $\fg$ for which $G_0$ has stable orbits in $\fg_1$, in the sense of Geometric Invariant Theory.

Every $G$-conjugacy class of positive-rank automorphisms $\theta\in\Aut(\fg)$ whose order is nonzero in $k$ contains a lift of a $W$-conjugacy class in $\Aut(R)$. For any particular group $G$ we can tabulate the Kac coordinates of such lifts; these are exactly the Kac coordinates of positive rank gradings.
For this purpose it is enough to consider only the lifts of certain classes in $\Aut(R)$, almost all of which are elliptic and $\bz$-regular in $\Aut(R')$ for some root subsystem of $R$, whose Kac coordinates are easily found, as above.  

These tables are only preliminary because they contain some Kac diagrams more than once, reflecting the fact that a given class in $\Aut(\fg)$ may contain lifts of several classes of 
$\sigma\in \Aut(R)$. However, each class in $\Aut(\fg)$ has a ``best" $\sigma$ whose properties tell us about other aspects of the grading, for example the little Weyl group $W(\fc)$. Our final tables for $E_6, E_7$ and $E_8$ list each positive rank Kac diagram once and contain this additional data. 

Besides its contributions to Vinberg theory {\it per se}, this paper was motivated by connections between Vinberg theory and the structure and representation theory of a reductive group $\bG$ over a $p$-adic field $F$. 
The base field $k$ above is then the residue field of a maximal unramified extension 
$L$ of $F$. We assume $\bG$ splits over a tame extension $E$ of $L$.  
Then the Galois group $\Gal(E/L)$ is cyclic and acts on the root datum of $\bG$ via a pinned automorphism $\vt$. 
The grading corresponds to a point $x$ in the Bruhat-Tits building of $\bG(L)$, 
the group $G_0$ turns out to be the reductive quotient of 
the parahoric subgroup $\bG(L)_x$ fixing $x$, 
and the summands $\fg_i$ are quotients in the Moy-Prasad filtration of $\bG(L)_x$. 
As we will show elsewhere, the classification of positive rank gradings leads to a classification of non-degenerate $K$-types, a long outstanding problem in the representation theory of $\bG(F)$, 
and stable $G_0$-orbits in the dual of $\fg_1$ give rise to supercuspidal representations of $\bG(F)$ attached to elliptic $\bz$-regular elements of the Weyl group. These generalize the ``simple supercuspidal representations" constructed in \cite{gross-reeder}, which correspond to the Coxeter element. 

After the first version of this paper was written, we learned from A. Elashvili that $25$ years ago he, D. Panyushev and E. Vinberg had also calculated, by completely different methods,  all the positive rank gradings and little Weyl groups in types $E_{6,7,8}$ (for $k=\bc$) but they had never published their results. We thank them for comparing their tables with ours. For other aspects of positive-rank gradings on exceptional Lie algebras, see \cite{degraaf-yakimova}.

\section{Kac coordinates}
Kac \cite[chap. 8]{kac:bluebook} showed how conjugacy classes of torsion automorphisms of simple Lie algebras $\fg$ (over $\bc$) 
can be parametrized by certain labelled affine Dynkin diagrams, called {\bf Kac coordinates}. 
If we choose a root of unity $\zeta\in\bc^\times$ of order $m$,  then any automorphism  $\theta\in\fg$ of order $m$ gives a grading $\fg=\oplus_{i\in\bz/m}\ \fg_i$, where $\fg_i$ is the $\zeta^i$-eigenspace of $\theta$. 
This grading depends on the choice of $\zeta$ and if we replace $\bc$ by another ground field $k$, we are forced to assume that $m$ is  invertible in $k$. As in \cite{levy:thetap}, this  assumption will be required for our classification of positive-rank automorphisms. 
 
However, at the level of classifying {\it all} torsion automorphisms,  Serre has remarked (see \cite{serre:kac})  that, at least in the inner case, one can avoid the choice of $\zeta$ and restrictions on $k$  by replacing an automorphism $\theta$ of order $m$  with an embedding
$\bfmu_m\hra \Aut(\fg)^\circ$ of group schemes over $k$, where $\bfmu_m$ is the group scheme of $m^{th}$ roots of unity.

In this section we give an elementary treatment of Kac coordinates in Serre's more general setting, and we extend his approach to embeddings $\bfmu_m\hra\Aut(\fg)$.
In the outer case, where the image of $\bfmu_m$ does not lie in $\Aut(\fg)^\circ$,  we still find it necessary to assume the characteristic $p$ of $k$ does not divide the order of the projection of $\bfmu_m$ to 
the component group of $\Aut(\fg)$.  
Our approach differs from  \cite{kac:bluebook} in that we avoid infinite dimensional Lie algebras (cf. \cite{reeder:torsion}). 

We then discuss a family of examples, the principal embeddings of $\bfmu_m$, which play an important role in gradings of positive rank.

\subsection{Based automorphisms and affine root systems}\label{affine}
For background on finite and affine root systems  see \cite{bour456} and \cite{macdonald:affine}.
Let $R$ be an irreducible reduced finite root system spanning a real vector space $V$. 
The automorphism group of $R$ is the subgroup of $\GL(V)$ preserving $R$:
$$\Aut(R)=\{\sigma\in\GL(V):\ \sigma(R)=R\}.$$
We say an automorphism $\sigma\in\Aut(R)$ is {\bf based} if $\sigma$ preserves a base of $R$. 
If we choose a base $\De$ of $R$ then we have a splitting 
$$\Aut(R)=W\rtimes\Th,$$
where $W$ is the Weyl group of $R$ and $\Th=\{\sigma\in\Aut(R):\ \sigma(\De)=(\De)\}$.  
Since $R$ is irreducible, the group $\Th$ is isomorphic to a symmetric group $S_n$ for $n=1,2$ or $3$. 

In this section we will associate to any based automorphism $\vt\in\Aut(R)$ an affine root system $\Psi(R,\vt)$ whose isomorphism class will depend only on the order $e$ of $\vt$.  

We first establish more  notation to be used throughout the paper. 
Let $X=\bz R$ be the lattice in $V$ spanned by $R$ and let $\check X=\Hom(X,\bz)$ be the dual lattice. We denote the canonical pairing between $X$ and $\check X$ by 
$\la \lam,\check \om\ra$, for $\lam\in X$ and $\check\om\in\check X$. 

Fix a base $\De=\{\al_1,\dots,\al_\ell\}$ of $R$, where $\ell$ is the rank of $R$, 
and let $\check R\subset\check X$ be the co-root system with base 
$\check\De=\{\check\al_1,\dots,\check\al_\ell\}$, where $\check\al_i$ is the co-root corresponding to $\al_i$. 
The pairing $\la\ ,\ \ra$ extends linearly to the real vector spaces $V=\br\otimes X $ and $\check V:=\br\otimes \check X $. Thus, a root $\al\in R$ can be regarded as the linear functional $\check v\mapsto \la \al,\check v\ra$ on $\check V$, and  by duality $\Aut(R)$ can be regarded as a subgroup of $\GL(\check V)$. In this viewpoint the Weyl group $W$ is the subgroup of $\GL(\check V)$ generated by the reflections 
$s_\al:\check v\mapsto \check v-\la \al,\check v\ra\check\al$ for $\al\in R$.

Let $\check\rho$ be one-half the sum of those co-roots $\check\al\in\check R$ which are non-negative integral combinations of elements of $\check\De$. We also have 
$$\check\rho=\check\om_1+\check\om_2+\cdots+\check\om_\ell,$$
where $\{\check\om_i\}$ are the fundamental co-weights dual to $\De$, that is, $\la\al_i,\check\om_i\ra=1$ and $\la\al_i,\check\om_j\ra=0$ if $i\neq j$. 

Let $\check V^\vt=\{\check v\in \check V:\ \vt(\check v)=\check v\}$ be the subspace of $\vt$-fixed vectors in $\check V$ and let 
$R_\vt=\{\al\vert_{\check V^\vt}:\ \al\in R\}$ be the set of restrictions to 
$\check V^\vt$ of roots in $R$. By duality $\Th$ permutes the fundamental co-weights $\{\check\om_i\}$, so the vector $\check\rho$ lies in $\check V^\vt$. 
And since $\la\al,\check\rho\ra=1$ for all $\al\in \De$, it follows that no root vanishes on $\check V^\vt$. Moreover two roots $\al,\al'\in R$ have the same restriction to $\check V^\vt$  if and only if they lie in the same $\la\vt\ra$-orbit in $R$.  Hence we have
$$R_\vt=\{\be_a:\ a\in R/\vt\},$$
where $R/\vt$ is the set of $\la\vt\ra$-orbits in $R$ and 
$\be_a=\al\vert_{\check V^\vt}$ for any $\al\in a$. 

For $a\in R/\vt$, we define $\check \be_a\in \check V^\vt$ by 
\begin{equation}\label{check}
\check \be_a=
\begin{cases} 
\ \ \ \sum_{\al\in a}\check \al & \quad\text{if}\quad 2\be_a\notin R_\vt\\
2\sum_{\al\in a}\check \al & \quad\text{if}\quad 2\be_a\in R_\vt,
\end{cases}
\end{equation}
and we set $\check R_\vt=\{\check\be_a:\ a\in R/\vt\}$. Then  
$\la \be_a,\check\be_a\ra=2$ and $\la \be_a,\check\be_b\ra\in\bz$ 
 for all $a,b\in R/\vt$. 

Note that $2\be_a\notin R_\vt$ precisely when $a$ consists of ``orthogonal" roots;
that is, when $a=\{\ga_1,\dots, \ga_k\}$ with $\la\ga_i,\check \ga_j\ra=0$ for $i\neq j$. 
In this case, the element 
$$s_a:=s_{\ga_1}s_{\ga_2}\cdots s_{\ga_k}\in W$$
has order two, is independent of the order of the product and is centralized by $\vt$. 
If $2\be_a\in R_\vt$ we have $a=\{\ga_1,\ga_2\}$ where $\ga_1+\ga_2\in R$.  
In this case we define $s_a=s_{\ga_1+\ga_2}$, noting this $s_a$ is also centralized by $\vt$. 
A short calculation shows that 
$$s_a(\be_b)=\be_b-\la \be_b,\check\be_a\ra\be_a,$$
in all cases. 
On the other hand, if $\be\in b$, then $s_a(\be_b)=s_a(\be)\vert_{\check V^\vt}$, 
since $s_a$ is centralized by $\vt$. It follows that 
$\be_b-\la \be_b,\check\be_a\ra\be_a\in R_\vt$. 
These involutions $s_a$, for $a\in R/\vt$, 
generate the centralizer $W^\vt=\{w\in W:\ \vt w=w\vt\}$ 
\cite[2.3]{steinberg:varchev}. 
Thus,  $R_\vt$ is a root system (possibly non-reduced) whose Weyl group is $W^\vt$.  
The rank $\ell_\vt$ of $R_\vt$ equals the number of $\vt$-orbits in $\De$. 


Let $\sca^\vt$ be an affine space for the vector space $\check V^\vt$. 
We denote the action by $(v,x)\mapsto v+x$ for $v\in\check V^\vt$ and $x\in \sca^{\vt}$ and 
for $x,y\in\sca^{\vt}$ we let $y-x\in \check V^\vt$ be the unique vector such that $(y-x)+x=y$. 
For any affine function $\psi:\sca^\vt\to \br$  we let $\dot\psi:\check V^\vt\to \br$ be the unique linear functional such that $\psi(x+v)=\psi(x)+\la\dot\psi,v\ra$ for all 
$v\in\check V^\vt$. 

Choose a basepoint $x_0\in\sca^{\vt}$. For each linear functional $\lam:\check V^\vt\to\br$ define an affine function 
$\widetilde\lam:\sca^{\vt}\to \br$ by $\widetilde\lam(x)=\la \lam,x-x_0\ra$.  In particular, each root 
$\be_a\in R_\vt$ gives an affine function $\widetilde\be_a$ on $\sca^\vt$. 

For each orbit 
$a\in R/\vt$, set $u_a=1/|a|$.
If $ \be_a\notin 2R_\vt$, define 
$$\Psi_a=\{\widetilde\be_a+nu_a:\ n\in \bz\}.$$
If $\be_a\in 2R_\vt$, define 
$$\Psi_a=\{\widetilde\be_a+(n+\half)u_a:\ n\in \bz\}.$$ 
The resulting collection 
$$\Psi(R,\vt):=\bigcup_{a\in R/\vt}\Psi_a$$
of affine functions on $\sca^\vt$ is a reduced, irreducible affine root system 
(in the sense of \cite[1.2]{macdonald:affine}) and $x_0\in\sca^\vt$ is a special point for $\Psi(R,\vt)$. 

An {\bf alcove} in $\sca^\vt$ is a connected component of the open subset of points in $\sca$ on which no affine function in $\Psi(R,\vt)$ vanishes. 
There is a unique alcove $C\subset\sca^{\vt}$  containing $x_0$ in its closure and on which $\tilde\be_a>0$ for every $\vt$-orbit $a\subset\De$.  
The walls of $C$ are hyperplanes $\psi_i=0$, 
$i=0,1,\dots,\ell_\vt=\dim \sca^\vt$, and $\{\psi_0,\psi_1,\dots,\psi_{\ell_\vt}\}$ is a base of the affine root system $\Psi(R,\vt)$. 
The point $x_0$ lies in all but one of these walls; we choose the numbering so that 
$\psi_0(x_0)\neq 0$. 
 There are unique relatively prime positive integers $b_i$ such that $\sum b_i\dot\psi_i=0$. 
 We have $b_0=1$ and the affine function 
 $\sum_{i=0}^{\ell_\vt}b_i\psi_i$ is constant, equal to $1/e$, where $e=|\vt|$. 
  The reflections $r_i$ about the hyperplanes $\psi_i=0$ for $i=0,1,\dots,\ell_\vt$ generate an irreducible affine Coxeter group $W_\aff(R,\vt)$ which acts simply-transitively on alcoves in $\sca^\vt$.

If $\vt=1$ we recover the affine root system attached to $R$ as in \cite{bour456} and 
$W_\aff(R):=W_\aff(R,1)$ is the affine Weyl group of $R$. 

For an example with nontrivial $\vt$, take $R$ of type $A_2$ and $\vt$ of order two. 
We have $\check V=\{(x,y,z)\in\br^3:\ x+y+z=0\}$, and
$$\al_1=x-y,\quad\al_2=y-z,\quad \check\al_1=(1,-1,0),\quad \check\al_2=(0,1,-1),\quad\check\rho=(1,0,-1).$$
The nontrivial automorphism $\vt\in\Aut(R)$ permuting $\{\al_1,\al_2\}$ acts on $\check V$ by $\vt(x,y,z)=(-z,-y,-x)$. We identify  
$\check V^\vt=\{(x,0,-x):\ x\in\br\}$ with $\br$ via projection onto the first component. 
The $\la\vt\ra$-orbits in the positive roots are $a=\{\al_1,\al_2\}$ and $b=\{\al_1+\al_2\}$, 
so $\be_a=x$ and $\be_{b}=2x$.
If we identify $\sca^\vt=\br$ and take $x_0=0$, then 
$$\Psi_a=\{x+\tfrac{n}{2}:\ n\in \bz\},\qquad \Psi_b=\{2x+n+\tfrac{1}{2}:\ n\in \bz\}.$$
The alcove $C$ is the open interval $(0,\tfrac{1}{4})$ in $\br$. The walls of $C$ are defined by the vanishing of the affine roots
$$\psi_0=\tfrac{1}{2}-2x,\qquad \psi_1=x$$
which satisfy the relation $\psi_0+2\psi_1=\frac{1}{2}$, so $b_0=1$ and $b_1=2$. 
The group $W_\aff(R,\vt)$ is infinite dihedral, generated by the reflections of $\br$ about $0$ and $\tfrac{1}{4}$. 

We list the affine root systems for nontrivial $\vt$ in Table 1. 
As the structure of $\Psi(R,\vt)$ depends only on $R$ and the order $e$ of $\vt$, the pair $(R,\vt)$ is indicated by the symbol ${^eR}$, called the {\it type} of $(R,\vt)$. 
Information about $\Psi(R,\vt)$ is encoded in a {\it twisted affine diagram} $D({^eR})$ which is a graph with vertices indexed by $i\in\{0,1,\dots,\ell_\vt\}$, labelled by the integers $b_i$. 
The number $m_{ij}$ of bonds between vertices $i$ and $j$ is determined as follows. 
Choose a $W^\vt$-invariant inner product $(\ ,\ )$ on $V^\vt$ and suppose that 
$(\dot\psi_j,\dot\psi_j)\geq (\dot\psi_i,\dot\psi_i)$. Then 
$$m_{ij}=\frac{(\dot\psi_j,\dot\psi_j)}{(\dot\psi_i,\dot\psi_i)}.$$
If $m_{ij}>1$ we put an arrow pointing from vertex $j$ to vertex $i$. 

 Removing the labels and arrows from the twisted affine diagram $D({^eR})$ gives the 
Coxeter diagram $D({^eR})_\cox$ of $W_\aff(R,\vt)$ (except in type ${^2A_2}$ the four bonds should be interpreted as $r_0r_1$ having infinite order).   Table 1 gives the twisted affine diagrams for $e>1$  (their analogues for $e=1$ being well-known).
 For each type we also give the {\it twisted Coxeter number}, 
which is the sum
\begin{equation}\label{twistedcoxeter}
h_\vt=e\cdot(b_0+b_1+\cdots+b_{\ell_\vt}),
\end{equation}
 whose importance will be seen later. 
The node $i=0$ is indicated by $\bullet$. 

\begin{center}
{\small Table 1: Twisted Affine diagrams and twisted Coxeter numbers}
\begin{equation}\label{coxnumber}
{\renewcommand{\arraystretch}{1.3}
\begin{array}{cccc}
\hline
{^eR} &D({^eR})& \ell_\vt & h_\vt\\
\hline 
{^2\!A_{2}} & \twoAtwo & 1 & 6\\
 
{^2\!A_{2n}} & 
 \overset{1}\bullet\!\!\Longrightarrow\!\!\overset{2}\circ \text{----}\!\! \overset{2}\circ\!\text{--} \cdots\text{--} 
 \overset{2}\circ\!\!\Longrightarrow\!\!\overset{2}\circ & n & 4n+2\\
  
{^2\!A_{2n-1}}&\begin{matrix}
 \overset{1}\circ \text{----}\!\!\!\!\! &\overset{2}\circ&
 \!\!\!\!\!\text{--} \cdots\text{--}
 \overset{2}\circ\!\!\Longleftarrow\!\overset{1}\circ\\
&\underset{1}{\text{\rotatebox{270}{\!\!\!\!\!\!\!\!\!----$\bullet$}}}&
 \end{matrix}&n&4n-2\\
  
{^2\!D_{n}} &\overset{1}\bullet\!\!\Longleftarrow\!\!
\overset{1}\circ \text{----}\!\! \overset{1}\circ\!\text{--} \cdots\text{--} 
 \overset{1}\circ\!\!\Longrightarrow\!\!\overset{1}\circ&n-1&2n\\
  
 {^3\!D_{4}}
 &\overset{1}\bullet\text{----}\overset{2}\circ\!\Lleftarrow\!\overset{1}\circ &2&12\\

{^2\!E_{6}} &\overset{1}\bullet\!\text{----}\!
\overset{2}\circ\! \text{----}\!
\overset{3}\circ\!\!\Longleftarrow\!\!
\overset{2}\circ\text{----}\!
\overset{1}\circ&4&18\\
\hline
\end{array}}
\end{equation}

\end{center}

 {\bf Remark:\ }
Let $\scr$ be the set of pairs $(R,e)$, where $R$ is an irreducible reduced finite root system and $e$ is a divisor of $|\Th|$.  Let $\scr_{\aff}$ be the set of irreducible reduced affine root systems, as in \cite{macdonald:affine}, up to isomorphism. 
Let $\scd$ be the set of pairs 
$(D,o)$, where $D$ is the Coxeter diagram of an irreducible affine Coxeter group and $o$ is a choice of orientation of each multiple edge of $D$. The classification of reduced irreducible affine root systems \cite[1.3]{macdonald:affine} shows that the assignments 
$(R,e)\mapsto {^eR}\mapsto D({^eR})$ give bijections 
$$\scr\overset\sim\lra\scr_{\aff}\overset\sim\lra\scd.$$

 \subsection{Torsion points, Kac coordinates and the normalization algorithm}\label{kac-coord}
 Retain the notation of the previous section. 
 Let $\sca^{\vt}_{\bq}$ be the set of points in $\sca^\vt$ on which the affine roots in $\Psi(R,\vt)$ take rational values. 
 The {\it order} of a point $x\in \sca^{\vt}_{\bq}$ is the smallest positive integer $m$ such that $\psi(x)\in\frac{1}{m}\bz$ for every $\psi\in\Psi(R,\vt)$. In this case there are integers $s_i$ such that $\psi_i(x)=s_i/m$, and  
 $\gcd(s_0,\dots,s_{\ell_\vt})=1$. 
Moreover, since $b_0\psi_0+\cdots+b_{\ell_\vt}\psi_{\ell_\vt}$ is constant, equal to $ 1/e$, (recall that $e$ is the order of $\vt$) it follows that 
 $$e\cdot \sum_{i=0}^{\ell_\vt} b_i s_i=m.$$
 In particular, the order $m$ is divisible by $e$. 
 We call integer vector $(s_0,s_1,\dots, s_{\ell_\vt})$ the (un-normalized) {\bf Kac coordinates} of $x$. 
 
 The point $x$ lies in $\overline C$ precisely when all $s_i$ are non-negative; in this case we refer to the vector $(s_i)$ as {\bf normalized Kac coordinates}. The action of the affine Weyl group 
 $W_{\aff}(R,\vt)$ on $\sca_\bq^\vt$ can be visualized as an action on Kac coordinates, as follows. 
 The reflection $r_j$ about the wall $\psi_j=0$ sends the Kac coordinates $(s_i)$ to $(s_i')$, where 
$$s_i'= s_i-\la \be_i,\check \be_j\ra s_j.$$

Un-normalized Kac coordinates may have some $s_j<0$. If we apply $r_j$ and repeat this process by selecting negative nodes and applying the corresponding reflections, we will eventually obtain 
normalized Kac coordinates $(s_i')$. Geometrically, this {\bf normalization algorithm} amounts to moving a given point $x\in \sca_\bq^\Th$ into the fundamental alcove $\overline C$ by a sequence of reflections about walls, see \cite[Sec. 3.2]{reeder:torsion}.
We have implemented the normalization  algorithm on a computer and used it extensively to construct the tables in sections \ref{E678} and \ref{2E6}. 

The image of the projection $e^{-1}\sum_{i=0}^{e-1}\vt^i:\check X\to  V^\vt$ is a lattice $Y_\vt$ in $V^\vt$ which is preserved by $W^\vt$. 
The extended affine Weyl group 
$$\wtW_{\aff}(R,\vt):=W^\vt\rtimes Y_\vt$$
contains $W_{\aff}(R,\vt)$ as a normal subgroup of finite index and 
the quotient may be identified with a group of symmetries of
the oriented diagram $D({^eR})$. 
We  regard two  normalized Kac diagrams as equivalent if one is obtained from the other by a symmetry of the oriented diagram $D({^eR})$ coming from $\wtW_{\aff}(R,\vt)$. 
For $R=E_6, E_7, E_8$ and $e=1$ these diagram symmetries are: rotation of order three, reflection of order two and trivial, respectively. In type ${^2E_6}$ these diagram symmetries are trivial (see table \eqref{coxnumber}).

 \subsection{$\mu_m$-actions on Lie algebras}\label{mum}
 Let $k$ be an algebraically closed field.  All $k$-algebras are understood to be commutative with $1$,  
and in this section all group schemes are affine over $k$, and are regarded as representable functors from the category of finitely generated $k$-algebras to the category of groups. We refer to \cite{waterhouse} for more details on affine group schemes. 
 
Every finitely generated $k$-algebra $A$ is a direct product of $k$-algebras
$A=\prod_{\iota\in I(A)}A_\iota,$
where $I(A)$ indexes the connected components $\Spec(A_\iota)$ of $\Spec(A)$ and each $A_\iota$ is a $k$-algebra with no non-trivial idempotents. This decomposition is to be understood when we describe the $A$-valued points in various group schemes below. 
Each finite (abstract) group $\varGamma$ is regarded a constant group scheme, given by 
$\varGamma(A)=\prod_{\iota\in  I(A)}\varGamma(A_\iota),$
where $\varGamma(A_\iota)=\varGamma$. In other words, an element $\gamma\in\varGamma(A)$ is a function $(\iota\mapsto\gamma_\iota)$ from $I(A)$ to $\varGamma$. 

Let $\bfmu_m$ denote the group scheme of $m^{th}$ roots of unity, 
whose $A$-valued points are given by 
$$\bfmu_m(A)=\{a\in A:\ a^m=1\}=\prod_{\iota\in  I(A)}\bfmu_m(A_\iota).$$
If $m$ is nonzero in $k$ then $\bfmu_m(A_\iota)=\bfmu_m(k)$ for every $\iota\in  I(A)$,
so $\bfmu_m$ is a constant group scheme and we have
$$\bfmu_m(A)=\prod_{\iota\in  I(A)}\bfmu_m(k).$$
If $m$ is zero in $k$ then $\bfmu_m$ is not a constant group scheme.

A $k$-vector space $V$ can be regarded as a $k$-scheme such that $V(A)=A\otimes_kV$. 
To give a grading $V=\sum_{i\in \bz/m\bz}V_i$ as $k$-schemes
 is to give a morphism $\varrho:\bfmu_m\to\GL(V)$, where 
 $\GL(V)(A)$ is the automorphism group of the free $A$-module $V(A)$. 
Indeed, $\bz/m$ is canonically  isomorphic to the  Cartier dual $\Hom(\bfmu_m,\bG_m)$, 
so a morphism $\varrho:\bfmu_m\to\GL(V)$ gives a grading $V(A)=\oplus_{i\in \bz/m}V_i(A)$ where $V_i(A)=\{v\in V(A):\ \varrho(\zeta)v=\zeta^i v\quad\forall \zeta\in \bfmu_m(A)\}$.

Now let $R$ be an irreducible root system as before, with base $\De$
and group of based automorphisms $\Th$.
Set $X=\bz R$ and $\check X=\Hom(X,\bz)$. 
Then $(X,R,\check X,\check R)$ is the root datum of a connected simple algebraic group scheme $G$ over $k$ of adjoint type. 
Let $\fg$ be the Lie algebra of $G$ and let $T\subset B$ be a maximal torus contained in a Borel subgroup of $G$. We identify $R$ with the set of roots of $T$ in $\fg$, and $\De$ with the set of simple roots of $T$ in the Lie algebra of $B$.  Choose a root vector $E_i$ for  each simple root 
$\al_i\in\De$. The data $(X,R,\check X,\check R,\{E_i\})$ is called a {\bf pinning} of $G$. 

Fix an element $\vt\in\Th$. Assume the order $e$ of $\vt$ is nonzero in $k$, 
so that $\bfmu_e$ and $\la \vt\ra$ are isomorphic constant group schemes over $k$,  
and  choose  an isomorphism $\tau:\bfmu_e\to\la\vt\ra$.

By our choice of pinning $(X,R,\check X,\check R,\{E_i\})$, the group $\la\vt\ra$ may also be regarded as a subgroup of $\Aut(\fg)$ permuting the root vectors $E_i$ in the same way $\vt$ permutes the roots $\al_i$, and we have a semidirect product
$$G\rtimes\la\vt\ra\subset\Aut(\fg),$$
where the cyclic group $\la\vt\ra$ is now viewed as a constant subgroup scheme of automorphisms of $\fg$, whose points in each $k$-algebra $A$ consist of vectors $(\vt^{n_\iota})$ acting on $\fg(A)=\prod_\iota\fg(A_\iota)$, with 
$\vt^{n_\iota}$ acting on the factor $\fg(A_\iota)$.

Now let $m$ be a positive integer divisible by $e$ (but $m$ could be zero in $k$). 
Let $m/e:\bfmu_m\to \bfmu_e$ be the morphism sending $\zeta\in\bfmu_m(A)$ to 
$\zeta^{m/e}\in\bfmu_e(A)$ for every $k$-algebra $A$. 

Finally, for each rational point $x\in\sca^\vt_\bq$ of order $m$ we shall now 
define a morphism 
$$\varrho_x:\bfmu_m\to T^\vt\times\la\vt\ra,$$
where $T^\vt$ is the subscheme of $\vt$-fixed points in $T$. 
We have $x=\tfrac{1}{m}\check \lam+ x_0$, for some 
$\check \lam\in\check  X^\vt$. 
The co-character $\check \lam$  restricts to a morphism 
$\check \lam_m:\bfmu_m\to T^\vt$ and 
we define $\varrho_x$ 
on $A$-valued points by 
$$\varrho_x(\zeta)=\check \lam_m(\zeta)\times\tau(\zeta^{m/e}),\qquad\text{for}\qquad
\zeta\in\bfmu_m(A).
$$
Since
$$\Hom(\bfmu_m,T^\vt)=\check X^\vt/m\check X^\vt\simeq 
\tfrac{1}{m}\check X^\vt/\check X^\vt,
$$
we see that $\check \lam_m$ corresponds precisely 
to an orbit of $x$ under  translation by $\check X^\vt$ on 
$\sca^\vt_\bq$.  
The condition that $x$ has order $m$ means that $\check \lam_m$ does not factor through $\bfmu_d$ for any proper divisor 
$d\mid m$. 

Let $\widetilde w\in \wtW_{\aff}(R,\vt)$ have projection $w\in W^\vt$ and denote the canonical action of 
$W^\vt$ on  $T^\vt$ by $w\cdot t$, for $t\in T^\vt(A)$. 
Then we have
$$\varrho_{\widetilde w\cdot x}(\zeta)=w\cdot \varrho_x(\zeta) $$
for all $\zeta\in \bfmu_m(A)$. 
One can check (cf. \cite[section 3]{reeder:torsion}) that two points $x,y\in \sca_\bq^\vt$ of order $m$ give $G$-conjugate embeddings $\varrho_x, \varrho_y:\bfmu_m\hra T^\vt\times \vt$ if and only if $x$ and $y$ are conjugate under $\wtW_{\aff}(R,\vt)$. 
The morphism $\varrho_x$ is thus determined by the Kac coordinates $(s_0,s_1,\dots, s_{\ell_\vt})$ of $x$ and the $G$-conjugacy class of $\varrho_x$ is determined by the normalized Kac coordinates of the $\wtW_{\aff}(R,\vt)$-orbit of $x$.

\subsection{Principal $\mu_m$-actions}\label{principalmu}
We continue with the notation of section \ref{mum}. Recall that $\check\rho\in \check X^\vt$ is the sum of the fundamental co-weights $\check\om_i$.
For every positive integer $m$ divisible by $e$, we have a {\bf principal} point 
$$x_m:=x_0+\tfrac{1}{m}\check\rho\in\sca^\vt_\bq$$
of order $m$. It corresponds to the {\bf principal embedding} 
$$\varrho_m=\varrho_{x_m}:\ \bfmu_m\lra T^\vt\times\la\vt\ra,\qquad\text{given by}\qquad
\varrho_m(\zeta)=\check\rho(\zeta)\times\tau(\zeta^{m/e}).
$$

The Kac coordinates of $x_m$ and $\varrho_m$ are given as follows. 
If $1\leq i\leq\ell_\vt$ we have $\psi_i=\tilde\be_i$ for some $\be_i\in R_\vt$ which is the restriction to $\check V^\Th$ of a simple root $\al_i\in \De$. Since $\la \al_i,\check\rho\ra=1$, it follows that $\la\psi,x_m\ra=1/m$ so $s_i=1$, and  we have
$$m=e\cdot \sum_{i=0}^{\ell_\vt}b_is_i=es_0+e\cdot\sum_{i=1}^{\ell_\vt}b_i=es_0+h_\vt-e,$$
where $h_\vt=e\cdot\sum_{i=0}^{\ell_\vt}b_i$ is the twisted Coxeter number of $R_\vt$ (see \eqref{twistedcoxeter}). Hence the remaining Kac-coordinate of the principal point $x_m$ is 
$$s_0=1+\frac{m-h_\vt}{e}.$$
This is negative if $m<h_\vt-e$, in which case we can apply the normalization algorithm of section \ref{kac-coord} to obtain the normalized Kac coordinates of $x_m$. Examples are found in the tables of section \ref{exceptional}. 

We will be especially interested in the points $x_m$ where $m$ is the order of an elliptic $\bz$-regular automorphism in $W\vt$ (defined in the next section). The twisted Coxeter number $h_\vt$ is one of these special values of $m$, corresponding to $s_0=1$  (cf. section \ref{stableclassification} below).

\section{$\bz$-regular automorphisms of root systems}\label{Zregular}

We continue with the notation of section \ref{affine}: $R$ is an irreducible finite reduced root system with a chosen base $\De$ and automorphism group  $\Aut(R)=W\rtimes\Th$, 
where $W$ is the Weyl group of $R$ and $\Th$ is the subgroup of $\Aut(R)$ preserving $\De$. 

\begin{definition}\label{def:regular} An automorphism $\sigma\in \Aut(R)$ is 
{$\boldmath\bz$}-{\bf regular} if the group generated by $\sigma$ acts freely on $R$. 
\end{definition}
This is nearly equivalent to Springer's notion of a regularity (over $\bc$) \cite{springer:regular}. 
In this section we will reconcile our definition with that of Springer. 

Let $X=\bz R$ be the root lattice of $R$ 
and let $\check X=\Hom(X,\bz)$ be the co-weight lattice. 
We say that a vector $\check v\in k\otimes \check X$ is {\boldmath $k$}-{\bf regular} if $\la \al,\check v\ra\neq 0$ for every $\al\in R$. We say also that an automorphism 
$\sigma\in\Aut(R)$ is {\boldmath $k$}-{\bf regular} if $\sigma$ has a $k$-regular eigenvector in $k\otimes\check X$. 
Taking $k=\bc$  we recover Springer's definition of regularity \cite{springer:regular}. 

At first glance it appears that $\sigma$ could be $k$-regular for some fields $k$ but not others. This is why we have defined regularity over $\bz$, as in Def. \ref{def:regular}. Of course the definition of $\bz$-regularity seems quite different from that of $k$-regularity. An argument due to  Kostant for the Coxeter element (cf.  \cite[Cor. 8.2]{kostant:betti}) shows that a $k$-regular automorphism is $\bz$-regular (see
\cite[Prop. 4.10]{springer:regular}). 
The converse is almost true but requires an additional condition. 
We will prove: 

\begin{prop}\label{prop:regular} An automorphism $\sigma\in \Aut(R)$
 is $\bz$-regular if and only if for every algebraically closed field $k$ in which the order $m$ of $\sigma$  is nonzero there is  $k$-regular eigenvector for $\sigma$ in $k\otimes\check X$ whose eigenvalue has order $m$. 
\end{prop}

Suppose $\sigma=w\vt$ where $w\in W$ and $\vt\in\Th$ is a based
 automorphism of order $e$. 
If $\sigma$ has order $m$ and has a $k$-regular eigenvalue $\lam$ of order $d$, 
then $m=\lcm(d,e)$. Indeed, it is clear that $m$ is divisible by 
$n:=\lcm(d,e)$. Conversely,  we have  $\lam^n=1$ so 
$\sigma^n$ fixes a regular vector, but $\sigma^n\in W$, so in fact $\sigma^n=1$ and $m\mid n$. 
Hence the notions of 
$\bz$-regularity  and $k$-regularity coincide precisely when $e\mid d$. 
In particular they coincide if $\vt=1$, that is, if $\sigma\in W$. 
However, if $\vt$ has order $e>1$ and we take $\sigma=\vt$,  then $\sigma$ fixes the $k$-regular vector $\check\rho$ so $\sigma$ is $k$-regular (if $e\neq 0$ in $k$). However 
$\sigma$ fixes the highest root, so $\sigma$ is not $\bz$-regular. And if $\zeta\in k^\times$ has order $e$ there are no $k$-regular vectors in the $\zeta$-eigenspace of $\sigma$.

The proof of Prop. \ref{prop:regular} will be given after some preliminary lemmas.

\begin{lemma}\label{based} An automorphism $\sigma\in\Aut(R)$ is based if and only if no root of $R$ vanishes on $\check X^\sigma$. 
\end{lemma}
\proof Assume that $\sigma\in \Aut(R)$ preserves a base $\De'\subset R$. Then $\sigma$ preserves the set $R^+$  of roots in $R$ which are non-negative integral linear combinations of roots in $\De'$. The vector $\sum_{\be\in R^+}\check\be$
belongs to $\check X^\sigma$ and no root vanishes on it. 

Conversely, let $\check v\in \check X^\sigma$ be a vector on which no root in $R$ vanishes. 
Then $v$ defines a chamber $\scc$ in the real vector space $\br\otimes X$, namely,
$$\scc=\{\lam\in \br\otimes X:\ \la\lam,\check v\ra>0\}.$$
As $\sigma$ fixes $\check v$, the chamber $\scc$ is preserved $\sigma$, 
so $\sigma$ permutes the walls of $\scc$. 
The set of roots $\al$ for which $\ker\check \al$ is a wall of $\scc$ is therefore a base of $R$ preserved by $\sigma$.
\qed

Next, we say that $\sigma\in \Aut(R)$ is {\bf primitive} if $\sigma$ preserves no proper root subsystem of $R$.

\begin{lemma}\label{primitive} 
If $\sigma\in \Aut(R)$ is primitive, then its characteristic polynomial 
on $V$ is irreducible over $\bq$. That is, we have 
$\det(tI_V-\sigma\vert_V)=\Phi_m(t)$, 
where $m$ is the order of $\sigma$ and $\Phi_m(t)\in\bz[t]$ is the cyclotomic polynomial whose roots are the primitive $m^{th}$ roots of unity. 
\end{lemma}
\proof 
In this proof we change notation slightly and let $V=\bq\otimes X$ denote the {\it rational} span of $X$ and let $\overline\bq$ be an algebraic closure of $\bq$. 

For $\al\in R$, let $V_\al\subset V$ be the rational span of  the $\sigma$-orbit of $\al$. Since $V_\al$ is spanned by roots, it follows from \cite[VI.1]{bour456} that 
$R\cap V_\al$ is a root subsystem of $R$. As it is preserved by the primitive automorphism 
$\sigma$, we must have $R\subset V_\al$, so $V_\al=V$. Hence the map $\bq[t]\to V$ given by sending $f(t)\mapsto f(\sigma)\al$ is surjective, and its kernel is the ideal in $\bq[t]$ generated by the minimal polynomial $M(t)$ of $\sigma$ on $V$. Hence $\deg M(t)=\dim V$ so we have $M(t)=\det(tI_V-\sigma\vert_V)$.

We must show that $M(t)$ is irreducible over $\bq$. If not, then $M(t)$ is divisible by $\Phi_d(t)$ for some proper divisor $d\mid m$. This means $\sigma$ has an eigenvalue of order $d$ on $\overline\bq\otimes V$, implying that $\sigma^d$ has nonzero fixed-point space $\check X^{\sigma^d}$. The set of roots vanishing on  $\check X^{\sigma^d}$ is a root subsystem not equal to the whole of $R$, and therefore is empty, again using the primitivity of $\sigma$. 

By Lemma \ref{based}, $\sigma^d$ is a nontrivial  automorphism preserving a base $\De'$ of $R$. As in the proof of that lemma, the sum of the positive roots for $\De'$ is a nonzero $\overline\bq$-regular vector in $V$ fixed by $\sigma^d$. 
Hence  the nontrivial subgroup $\la\sigma^d\ra$ has trivial intersection with $W$. 
If $\sigma\in W$ this is a contradiction and the lemma is proved in this case. 

Assume that $\sigma\notin W$.  
Since $R$ is irreducible and we have shown that the projection 
$\Aut(R)\to \Th$ is injective on $\la\sigma^d\ra$, 
it follows that $\sigma^d$ has order $e\in\{2,3\}$. 
We must also have $(e,d)=1$ and $m=ed$. 
As $e$ is determined by the projection of $\sigma$ to $\Th$, it follows that  $d$ is the {\it unique} proper divisor of $m$ such that $\Phi_d(t)$ divides $M(t)$. 
Since the roots of $M(t)$ are $m^{th}$ roots of unity (because $\sigma^m=1$) 
and are distinct (since $\sigma$ is diagonalizable on $\overline\bq\otimes V$) 
and $M(t)\neq \Phi_d(t)$ by assumption, it follows that $M(t)=\Phi_m(t)\cdot \Phi_d(t)$. 

If $e=2$ then $-\sigma\in W$ is also primitive, with reducible minimal polynomial $M(-t)=
\Phi_m(-t)\cdot \Phi_d(-t)$, contradicting the case of the lemma previously proved. 
If $e=3$, then $\Phi$ has type $D_4$, so $m=3d$ and
$$4=\deg M=\phi(3d)+\phi(d)=\phi(d)[\phi(3)+1]=3\phi(d),$$
which is also impossible. The lemma is now proved in all cases. 
\qed

Now let $\sigma\in \Aut(R)$ be a $\bz$-regular automorphism of order $m$. 
Recall from Def. \ref{def:regular} that this means the group $\la\sigma\ra$ generated by $\sigma$ acts freely on $R$. 
For each $\al\in R$, let $V_\al\subset \bq\otimes X$ denote the 
$\bq$-span of the $\la \sigma\ra$-orbit of $\al$
and let $M_\al(t)$ be the minimal polynomial of $\sigma$ on $V_\al$.

\begin{lemma}\label{free}
 If $\sigma$ is $\bz$-regular of order $m$  then $\Phi_m(t)$ divides $M_\al(t)$ in $\bz[t]$, for all $\al\in R$. 
\end{lemma}
\proof Let $\zeta\in \overline\bq^\times$ be a root of unity of order $m$ and let $\al\in R$. 
It suffices to show that $\zeta$ is an eigenvalue of $\sigma$ in $\overline\bq\otimes V_\al$. 
Let $R'$ be a minimal (nonempty) $\sigma$-stable root subsystem of $R\cap V_\al$, 
and let 
$$R'=R'_0\cup R'_1\cup\cdots\cup R'_{k-1}$$
be the decomposition of $R'$ into irreducible components.
These are permuted transitively by $\sigma$; we index them so that $\sigma^iR'_0=R'_i$ for $i\in \bz/k$. The stabilizer of 
$R'_0$ in $\la \sigma\ra$ is generated by $\sigma^k$.
Correspondingly, the rational span $V'$ of $R'$ is a direct sum 
$$V'=V'_0\oplus V'_1\oplus\cdots\oplus V'_{k-1}\subset V_\al$$
where $V'_i$ is the rational span of $R'_i$. 

Suppose that $\eta:=\zeta^k$ is an eigenvalue of $\tau:=\sigma^k$ in $\overline\bq\otimes V'_0$, 
afforded by the vector $v\in \overline\bq\otimes V'_0$.  
Let $S$ and $T$ denote the group algebras over $\overline\bq$ 
of $\la\sigma\ra$ and $\la \tau\ra$, respectively, 
and let $\overline\bq_\eta$ be the  $T$-module with underlying vector space $\overline\bq$ 
on which $\tau$ acts as multiplication by $\eta$. 
There is a unique map of $S$-modules 
$$f:S\otimes_T \overline\bq_\eta\lra V'$$
such that $f(1\otimes 1)=v\in V'_0$. 
As $f(\sigma^i\otimes 1)=\sigma^iv\in \overline\bq\otimes V'_i$, 
and the spaces 
$V'_0,  V'_1, \dots, V'_{k-1}$ are linearly independent, 
it follows that $f$ is injective. Frobenius reciprocity implies that $\zeta$ appears as an eigenvalue of $\sigma$ in $\overline\bq\otimes V'$, hence also in 
$\overline\bq\otimes V_\al$. 

It therefore suffices to prove that $\eta$ appears as an eigenvalue of $\tau$ on $\overline\bq\otimes V'_0$. 
Since $\sigma$ acts freely on $R$, it follows that $\tau$ acts freely on $R'_0$ and has order $n:=m/k$ on $R'_0$. We claim that $\tau$ is primitive on $R'_0$. 
For if $R''\subset R'_0$ is a root subsystem preserved by $\tau$ then 
$R''\cup \sigma R''\cup\cdots\cup \sigma^{k-1} R''$ is a root subsystem preserved by $\sigma$ which must equal $R'$ (by minimality), so that $R''=R'_0$. Hence $\tau$ is indeed primitive on $R'_0$. 
By Lemma \ref{primitive} the characteristic polynomial of $\tau$ on $V_0'$ is the cyclotomic polynomial $\Phi_n(t)$,
which has the root $\zeta^{m/n}=\zeta^k=\eta$. Therefore $\eta$ appears as an eigenvalue of $\tau$ on $\overline\bq\otimes V'_0$, as desired. 
\qed

We are now ready to prove Prop. \ref{prop:regular}. 
Let $k$ be an algebraically closed field and set $V_k:=k\otimes X$, $\check V_k:=k\otimes \check X$. Recall that a $k$-regular vector $\check v\in \check V_k$ is one for which 
$\la\al,\check v\ra\neq 0$ for all $\al\in R$. 

For completeness we  recall the proof of the easy direction of Prop. \ref{prop:regular} (cf. \cite[4.10]{springer:regular}). 
Assume that $\sigma\in\Aut(R)$ is $k$-regular, and let $\check v\in \check V_k$ be a $k$-regular eigenvector of $\sigma$ with eigenvalue $\zeta\in k^\times$ of order $m$ equal to the order of $\sigma$. 
Suppose $\sigma^d\al=\al$ for some $\al\in R$. Then 
$$
0\neq \la \al,\check v\ra=\la \sigma^d\al,\check v\ra=\la \al,\sigma^{-d}\check v\ra
=\zeta^{-d}\la\al,\check v\ra.
$$
It follows that $\zeta^d=1$. Since $\sigma$ and $\zeta$ have the same order, it follows that 
$\sigma^d=1$. Hence $\la \sigma\ra$ acts freely on $R$, so $\sigma$ is $\bz$-regular.  

Assume now that $\sigma$ is $\bz$-regular, so that $\la\sigma\ra$ acts freely on $R$. 
Let $\bar\Phi_m(t)$ denote the image, under the map $\bz[t]\to k[t]$ induced by the canonical map $\bz\to k$, of the cyclotomic polynomial $\Phi_m(t)$. Since $m$ is nonzero in $k$, it follows that
all roots of $\bar \Phi_m(t)$ in $k$ have order $m$. 
Let $\zeta\in k^\times$ be one of them. 

Let $\al\in R$ and let $X_\al$ be the subgroup of $X$ generated by the 
$\la\sigma\ra$-orbit of $\al$. 
Then $X_\al$ is a lattice in $V_\al=\bq\otimes X_\al$ and 
$\Phi_m(t)$ divides the characteristic polynomial $\det(tI-\sigma\vert_{X_\al})$ 
in $\bz[t]$, by Lemma \ref{free}. 
Hence $\bar\Phi_m(t)$ divides $\det(tI-\sigma\vert_{k\otimes X_\al})$ in $k[t]$. 
In particular $\zeta^{-1}$ is an eigenvalue of $\sigma$ on $k\otimes X_\al$. 

The operator $ P_\zeta\in \End(V_k)$ given by 
$$P_\zeta
=1+\zeta \sigma+\zeta^{2} \sigma^2+\cdots+\zeta^{m-1} \sigma^{m-1}
$$
preserves $k\otimes X_\al$ and $P_\zeta(k\otimes X_\al)$
is the $\zeta^{-1}$-eigenspace of $\sigma$ in $k\otimes X_\al$. 
As $X_\al$ is spanned by roots $\sigma^i\al$ and 
$P_\zeta(\sigma^i\al)=\sigma^{-i}P_\zeta(\al)$, it follows that $P_\zeta(\al)\neq 0$. 

As $\al\in R$ was arbitrary, we have that $P_\zeta(\al)\neq 0$  for all $\al\in R$. Since $k$ is infinite, there exists $\check v\in \check V_k$ such that $\la P_\zeta(\al),\check v\ra\neq 0$ for all $\al\in R$. 

The dual projection
$$\check P_\zeta
=1+\zeta^{-1}\sigma+\zeta^{-2}\sigma^2+\cdots+\zeta^{1-m} \sigma^{m-1}\in \End(\check V_k)$$
satisfies 
$$\la \al, \check P_\zeta (\check v)\ra=\la  P_\zeta(\al), \check v\ra\neq 0, $$
for all $\al\in R$. Therefore $\check P_\zeta (\check v)$ is a $k$-regular eigenvector of $\sigma$ in $\check V_k$ whose eigenvalue $\zeta$ has order $m$. 
This completes the proof of Prop. \ref{prop:regular}.
\qed

\section{Positive rank gradings}\label{sec:posrank}
Let $\fg$ be the Lie algebra of a connected simple algebraic group $G$ of adjoint type over an algebraically closed field $k$ whose characteristic is not a torsion prime for $G$. 
Then $G=\Aut(\fg)^\circ$ is the identity component of $\Aut(\fg)$. 
We fix a Cartan subalgebra $\ft$ of $\fg$ with corresponding maximal torus $T=C_G(\ft)$ and let $R$ be the set of roots of  $\ft$ of $\fg$. Let $N=N_G(T)$ be the normalizer of $T$, so that $W=N/T$ is the Weyl group of $R$.

From now on we only consider gradings $\fg=\oplus_{i\in\bz/m}\ \fg_i$ whose 
period $m$ is nonzero in $k$. By choosing an $m^{th}$ root of unity in $k^\times$, 
we get an automorphism $\theta\in\Aut(\fg)$ of order $m$, 
such that $\theta$ acts on $\fg_i$ by the scalar $\zeta^i$.

In this section we show how all such gradings of positive rank may be effectively found by computing lifts to $\Aut(\fg)$ of 
automorphisms $\sigma\in\Aut(R)$.

 
 \subsection{A canonical Cartan subalgebra}\label{csa}
Given any Cartan subalgebra $\fs$ of $\fg$ with centralizer $S=C_G(\fs)$, 
let 
$$\Aut(\fg,\fs)=\{\theta\in\Aut(\fg):\ \theta(\fs)=\fs\}.$$
We have an isomorphism (obtained by conjugating $\fs$ to our fixed Cartan subalgebra $\ft$) 
$$\Aut(\fg,\fs)/S\simeq \Aut(R)$$
which is unique up to conjugacy in $\Aut(R)$.  Thus any element of $\Aut(\fg,\fs)$ gives a well-defined conjugacy class in $\Aut(R)$. 
However, an automorphism $\theta\in\Aut(\fg)$ may normalize various Cartan subalgebras $\fs$, giving rise to various classes in $\Aut(R)$. 
We will define a canonical $\theta$-stable Cartan subalgebra, which will allow us associate to $\theta$ a well-defined conjugacy class in $\Aut(R)$.

For each $\theta\in\Aut(\fg)$ whose order is nonzero in $k$ we define a canonical $\theta$-stable Cartan subalgebra $\fs$ of $\fg$ as follows. 
Let $\fc\subset\fg_1$ be a Cartan subspace. 
The centralizer $\fm=\fz_\fg(\fc)$ is a $\theta$-stable Levi subalgebra of $\fg$
and we have $\fm=\oplus\fm_i$ where $\fm_i=\fm\cap\fg_i$. 
Choose a Cartan subalgebra $\fs_0$ of $\fm_0$. 
Then $\fs_0$ contains regular elements of $\fm$ \cite[Lemma 1.3]{levy:thetap}, so the centralizer 
$$\fs:=\fz_\fm(\fs_0)$$ 
is a $\theta$-stable Cartan subalgebra of $\fm$, and $\fs$ is also a Cartan subalgebra of $\fg$. 
We have $\fs\cap \fg_0=\fs_0$ (so our notation is consistent) and $\fs\cap \fg_1=\fc$. 
Since $G_0$ is transitive on Cartan subspaces in $\fg_1$ \cite[Thm. 2.5]{levy:thetap} and 
$C_{G_0}(\fc)^\circ$ is transitive on Cartan subalgebras of its Lie algebra $\fm_0$, 
the Cartan subalgebra $\fs$ is unique up to $G_0$-conjugacy.

\subsection{A relation between $\Aut(\fg)$ and $\Aut(R)$}\label{vdash}
For $\theta\in\Aut(\fg)$ and $\sigma\in \Aut(R)$ we write 
$$\theta\vdash \sigma$$ 
if the following two conditions are fulfilled:
\begin{itemize}
\item $\theta$ and $\sigma$ have the same order;
\item $\theta$ is  $G$-conjugate to an automorphism $\theta'\in\Aut(\fg,\ft)$  such that $\theta'\vert_\ft=\sigma$. 
\end{itemize}

Assume that $\theta\vdash \sigma$ and that the common order $m$ of $\theta$ and $\sigma$ is nonzero in $k$. Choose a root of unity $\zeta\in k^\times$ of order $m$, giving a grading $\fg=\oplus_{i\in\bz/m}\ \fg_i$. 
Recall that $\rank(\theta)$ is the dimension of a Cartan subspace $\fc\subset\fg_1$ for $\theta$. 
Likewise, for $\sigma\in \Aut(R)$, let $\rank(\sigma)$ be the multiplicity of $\zeta$ as a root of the characteristic polynomial of $\sigma$ on $V$. Since $\ft$ consists of semisimple elements, it follows that 
$\rank(\theta)\geq \rank(\sigma)$. 

\begin{prop}\label{posrank} Let $\theta\in\Aut(\fg)$ be an automorphism of positive rank whose order $m$ is  nonzero in $k$. Then 
$$\rank(\theta)=\max\{\rank(\sigma):\ \theta\vdash \sigma\}.$$
\end{prop}
\proof
It suffices to show that there exists $\sigma\in \Aut(\fg,\ft)$ such that $\theta\vdash \sigma$ and $\rank(\theta)=\rank(\sigma)$. 

Replacing $\theta$ by a $G$-conjugate,  we may assume that $\ft$ is the canonical Cartan subalgebra for $\theta$ (section \ref{csa}) so that $\theta\in\Aut(\fg,\ft)$, 
and $\fc=\ft_1$ is a Cartan subspace contained in $\ft$. 
Then $\fc$ is the $\zeta$-eigenspace of  $\sigma:=\theta\vert_\ft\in\Aut(R)$. 
Since $\theta$ has order $m$, it follows that the order of $\sigma$ divides $m$. 
But  $\sigma$ has an eigenvalue of order $m$, so the order of $\sigma$ is exactly $m$. We therefore have $\theta\vdash \sigma$ and 
$\rank(\theta)=\dim\fc=\rank(\sigma)$. 
\qed

Given $\sigma\in \Aut(R)$ let $\Kac(\sigma)$ denote the set of normalized Kac diagrams of automorphisms $\theta\in \Aut(\fg,\ft)$ for which $\theta\vdash \sigma$.  Since there are only finitely many Kac diagrams of a given order, each set $\Kac(\sigma)$ is finite. 
From Prop. \ref{posrank} it follows that  
the Kac coordinates of all positive rank automorphisms of $\fg$ 
are contained in the union
\begin{equation}\label{positiverank}
\bigcup_{\sigma\in\Aut(R)/\sim}\Kac(\sigma),
\end{equation}
taken over representatives of  the $W$-conjugacy classes in $\Aut(R)$. 
Moreover $\rank(\theta)$ is the maximal $\rank(\sigma)$ for which the Kac coordinates of $\theta$ appear in $\Kac(\sigma)$. 

\subsection{Inner automorphisms}
If $\theta\in G=\Aut(\fg)^\circ$ is inner then its Kac diagram 
will belong to $\Kac(w)$ for some $w\in W$. 
In this section we refine the union \eqref{positiverank} to reduce the number of classes of $w$ to consider, and we show how to compute $\Kac(w)$ directly from $w$, for these classes. 

A subset $J\subset\{1,\dots,\ell\}$ is {\bf irreducible} if the root system $R_J$ spanned by $\{\al_j:\ j\in J\}$ is irreducible. Two subsets $J,J'$ are {\bf orthogonal} if $R_J$ and $R_{J'}$ are orthogonal. 

An element $w\in W$ is {\bf $m$-admissible} if $w$ has order $m$ and $w$ can be expressed as a product
\begin{equation}\label{admissible}
w=w_1 w_2\cdots w_d,
\end{equation} where 
each $w_i$ is contained in $W_{J_i}$ for irreducible mutually orthogonal subsets $J_1,\dots, J_d$ of $\{1,2,\dots,\ell\}$ and on the reflection representation of $W_{J_i}$ each $w_i$ has
an eigenvalue of order $m$ but no eigenvalue equal to $1$ (so $w_i$ is elliptic in $W_{J_i}$). We call \eqref{admissible} an {\bf admissible factorization} of $w$. Note that each $w_i$ also has order 
$m$, that  $\rank(w)=\sum_i\rank(w_i)$, 
and $\rank(w_i)>0$ for $1\leq i\leq d$. 

Let $G_i$ be the Levi subgroup of $G$ containing $T$ and the roots from $J_i$, 
and let $G_i'$ be the derived group of $G_i$. 
Each $w_i\in W_{J_i}$ has a lift $\dot w_i\in G_i'\cap N$ and all such lifts are conjugate by $T\cap G_i$, hence the  normalized Kac-coordinates of $\Ad(\dot w_i)$ in $\Ad(G_i')$ are well-defined. 

Given an $m$-admissible element $w=w_1\cdots w_d$ as  in \eqref{admissible}, 
let   $\Kac(w)_{\un}$ be the set of 
un-normalized Kac coordinates $(s_0,s_1,\dots,s_\ell)$ such that 
\begin{itemize}
\item For $j\in J_i$ the coordinate $s_j$ is the corresponding normalized Kac coordinate of $w_i$ in $G_i'$. 
\item For $i\in\{0,1,\dots,\ell\}-J$, the coordinate $s_i$ ranges over a set of representatives for $\bz/m$. 
\item $\sum_{i=0}^\ell a_i s_i=m$. 
\end{itemize}

If $w$ is any automorphism of $T$ we set $(1-w)T:=\{t\cdot w(t)^{-1}:\ t\in T\}$. 

\begin{lemma}\label{lem:Kac(w)} If $w$ is $m$-admissible, then $\Kac(w)$ is the set of Kac diagrams obtained by applying the normalization algorithm of section \ref{kac-coord} to the elements of 
$\Kac(w)_{\un}$. 
\end{lemma}
\proof Each Kac diagram in $\Kac(w)_{\un}$ is that of a lift of $w$ in $N$ of order $m$. 
Hence the normalization of this diagram lies in $\Kac(w)$. Conversely, suppose
$(s_i)$ are normalized Kac coordinates lying in $\Kac(w)$. 
By definition, there is an inner automorphism $\theta\vdash w$ (notation of section \ref{vdash}) of order $m$ with these normalized Kac-coordinates, and we may assume that $\theta=\Ad(n)$ for some $n\in N$, 
a lift of $w$. Then 
$$n=\dot w_1\dot w_2\cdots \dot w_d\cdot t$$
where each $\dot w_i$ is a lift of $w_i$ and $t\in T$. 
Let $Z$ be the maximal torus in the center of $G_1'\cdot G_2'\cdots G_d'\cdot T$. 
Then $T=Z\cdot (1-w)T$, so we may conjugate $n$ by $T$ to arrange that $t\in Z$. 
Next, we conjugate each $\dot w_i$ in $G_i'$ to an element $t_i\in T\cap G_i'$, 
thus conjugating $n$ to 
$$n'=t_1\cdot t_2\cdots t_d\cdot t\in T.$$
Since $n'$ has order $m$ there exists $\check\lam\in\check X$ such that 
$n'=\check\lam(\zeta)$. As in section \ref{kac-coord}, the point $x=x_0+\frac{1}{m}\check\lam\in\sca_\bq$ has order $m$ and the simple affine roots $\psi_i$ take values $\psi_i(x)=s_i'/m$, where $s_i'$ are the Kac coordinates of $n'$ and $\sum_{i=0}^\ell a_i s'_i=m$.
If $j\in J_i$ then  $s'_j$ is a Kac coordinate of the $G'_i$-conjugate 
$\dot w_i$ of $t_i$, and if 
$i\in\{0,1,\dots,\ell\}-J$ we have $\al_i(n')=\zeta^{s'_i}$, so the class of $s_i'$ in $\bz/m$ is determined.  
Hence the  Kac coordinates $(s'_i)$ lie in $\Kac(w)_{\un}$ and their normalization is $(s_i)$. 
\qed

\begin{prop}\label{Kac(w)} 
Let $\theta\in\Aut(\fg)^\circ$ be an inner automorphism of order $m$ nonzero in $k$ with  
$\rank(\theta)>0$. Then there exists an $m$-admissible element $w\in W$ such that $\theta\vdash w$, and the rank of $\theta$ is given by
$$\rank(\theta)=\max\{\rank(w):\ \theta\vdash w\},$$
where the maximum is taken over all $W$-conjugacy classes of $m$-admissible elements $w\in W$ such that $\theta\vdash w$. 
\end{prop}
\proof We may assume that  $\ft$ is the canonical Cartan subalgebra for $\theta$, so that 
$\theta=\Ad(n)$ for some $n\in N$. The element  $w=nT\in N/T=W$ has order $m$ and $\theta\vdash w$. Recall that the canonical Cartan subalgebra has the property that $\ft_1$ is a 
Cartan subspace for $\theta$. Hence $\rank(\theta)=\rank(w)>0$ .

Assume first that $\ft_0=0$, that is, $w$ is elliptic. 
Then $w$ is $m$-admissible and its admissible factorization 
\eqref{admissible} is $w=w_1$, with $d=1$, 
 so the proposition is proved in this case. 

Assume now that $\ft_0\neq 0$. 
Let $R_0$ be the set of roots in $R$ vanishing on $\ft_0$. Since $R_0$ is the root system of a Levi subgroup of $G$, 
there is a basis $\De=\{\al_1,\al_2,\dots,\al_\ell\}$ of $R$ such that 
$\De_0:=\De\cap R_0$ is a basis of $R_0$. 
We have $\De_0=\{\al_j:\ j\in J\}$ for some subset $J\subset\{1,2,\cdots,\ell\}$. 
Decomposing $R_0$ into irreducible root systems $R_0^{i}$, we have 
corresponding decompositions
\begin{displaymath}
\begin{split}
R_0&=R_0^1\cup R_0^2\cup\cdots\cup R_0^n,\\
\De_0&=\De_0^1\cup \De_0^2\cup\cdots\cup \De_0^n,\\
J&= J_1\cup J_2\cup\cdots\cup J_n,\\
W_J&=W_{J_1}\times W_{J_2}\times\cdots\times W_{J_n},\\
w&=w_1\cdot w_2\cdots\cdot w_n.
\end{split}
\end{displaymath}
By construction,  $w$ is elliptic in $W_J$ and has an eigenvalue of order $m$ on the reflection representation of $W_J$. 
Therefore, each $w_i$ is elliptic in $W_{J_i}$ and has eigenvalues of order dividing $m$. 
And since $\rank(w)>0$ there is some number $d\geq 1$ of $w_i$'s having an eigenvalue of order exactly $m$. 
Let the factors be numbered so that $w_i$ has an eigenvalue of order $m$ for $i\leq d$, 
and $w_i$ has no eigenvalue of order $m$ for $i>d$. The element 
$$w'=w_1w_2\cdots w_d$$
is $m$-admissible. 

As before, let  $G_i$ be the Levi subgroup of $G$ containing $T$ and the root subgroups from $J_i$,
and let $G_i'$ be the derived subgroup of $G_i$. The derived group of $C_G(\ft_0)$ is a 
commuting product 
$G_1'\cdot G_2'\cdots G_n'$. 

Each $w_i$ has a lift 
$\dot w_i\in N\cap G_i'$; such a lift is unique up to conjugacy by $T\cap G_i'$ 
and we have 
$$\theta=\dot w_1\dot w_2\cdots \dot w_n \cdot t$$
for some $t\in T$. For $i>d$ we conjugate $\dot w_i$ in $G'_i$ to an element $t_i\in T$, 
obtaining a conjugate $\theta'$ of $\theta$ having the form
$$\theta'=\dot w_1\dot w_2\cdots \dot w_d \cdot t'.$$
Therefore $\theta\vdash w'$ and $w'$ is $m$-admissible of the same rank as $\theta$. 
The proposition is proved. 
\qed

\section{Principal and stable gradings}\label{principalstable}
Retain the set-up of section \ref{sec:posrank}.
Let $B$ be a Borel subgroup of $G=\Aut(\fg)^\circ$ containing our fixed maximal torus $T$.  The algebraic group $G$ has root datum 
$(X,R,\check X, \check R)$, where $X=X^\ast(T)$ (resp. $\check X=X_\ast(T)$) are the lattices of weights (resp. co-weights) of $T$, and $R$ (resp. $\check R$) are the sets of  roots (resp. co-roots) of $T$ in $G$. 
The base $\De$ of $R$ is the set of simple roots of $T$ in $B$. 
As before, we choose a pinning 
$(X,R,\check X, \check R,\{E_i\})$, where $E_i\in\fg$ is a root vector for the simple root $\al_i\in\De$. This choice
gives an isomorphism  from $\Aut(R,\De)$ to the group 
$\Th=\{\vt\in\Aut(\fg,\ft):\ \vt\{E_i\}=\{E_i\}\ \}$ of pinned automorphisms, and we have a splitting
$$\Aut(\fg)=G\rtimes\Th.$$

\subsection{Principal gradings}\label{principalgradings}
For each positive integer $m$ and pinned automorphism $\vt\in\Aut(R,\De)$, we have a principal grading $\fg=\oplus_{i\in\bz/m}\ \fg_i$ given (as in section \ref{principalmu}) by the point  $x_m:=\frac{1}{m}\check\rho+x_0$ 
(Recall that $\check \rho$ is the sum of the fundamental co-weights dual to the simple roots $\al_i\in\De$.)
The  normalized Kac diagram of $x_m$  may be obtained via the algorithm described in section \ref{principalmu}. (These Kac diagrams may also be found in \cite{degraaf:niltheta}.)

Note that $\fg_1$ contains the regular nilpotent element 
$E:=E_1+E_2+\cdots+E_\ell$ associated to our pinning. 
If $m$ is nonzero in $k$ and we choose a root of unity $\zeta\in k^\times$ of order $m$, then $\fg_i$ is the $\zeta^i$-eigenspace for the automorphism 
$$\theta_m:=\check\rho(\zeta)\vt.$$ 
Note that the $\zeta$-eigenspace $\fg_1$ for $\theta_m$ contains the regular nilpotent element $E:=E_1+E_2+\cdots+E_\ell$ associated to our pinning. Conversely if 
$\theta=\check\lam(\zeta)\vt$ is an automorphism of order $m$ whose $\fg_1$ contains a regular nilpotent element then $\theta$ is principal. 
If the characteristic $p$ of $k$ is zero or sufficiently large, the element
$\check\rho(\zeta)$ is the image of $\begin{bmatrix}\zeta&0\\0&1\end{bmatrix}$ under the principal embedding $\PGL_2\hra G$ associated by the Jacobson-Morozov theorem to $E$. Elsewhere in the literature a principal automorphism is called ``$N$-regular". 

The first aim of this section is to show that  lifts to $\Aut(\fg)$ of $\bz$-regular elliptic automorphisms 
$\sigma\in\Aut(R)$ are principal. (Recall that an automorphism $\sigma\in\Aut(R)$ is called {\bf elliptic} if $X^\sigma=0$.) 

More precisely, let $\sigma=w\vt\in W\vt$ be an elliptic $\bz$-regular automorphism of $R$ (Def. \ref{def:regular}). 
 Let $n\in N$ be a lift of $w$. Since $\sigma$ is elliptic the fixed-point group $T^\sigma$ is finite, so the coset $nT\vt\subset G\vt$ consists of a single $T$-orbit under conjugation. It follows that the $G$-conjugacy class $C_\sigma$ of $n\vt$ in $G\vt$ depends only on $\sigma$. In this section we will prove the following. 

\begin{prop}\label{clift} Assume $\sigma\in W\vt$ is elliptic and $\bz$-regular and that the order $m$ of $\sigma$ is nonzero in $k$. 
Then the conjugacy class $C_\sigma$ contains $\check\rho(\zeta)\vt$ for every $\zeta\in k^\times$ of order $m$. 
\end{prop}


The second aim of this section is to characterize the principal gradings which arise from 
 elliptic $\bz$-regular automorphisms of $R$ in terms of stability (see section \ref{sec:stable}).

\subsection{Conjugacy results}
If $\sigma$ is an automorphism of an abelian group $A$, we set 
\[
(1-\sigma)A:=\{a\cdot\sigma(a)^{-1}:\ a\in A\}.
\]
Let $N^\vt, W^\vt$ denote the fixed-point subgroups of $\vt$ in $N, W$ respectively,
and let $N_\vt=\{n\in N:\ \vt(n)\equiv n\mod T\}$. 
It is known (see \cite{steinberg:endo}) that $N_\vt=N^\vt\cdot T$. This group acts on the coset $T\vt$ by conjugation. 
Meanwhile the fixed-point group $W^\vt$ acts on the quotient torus 
$$T_\vt=T/(1-\vt)T$$
whose character and cocharacter groups 
$X^\ast(T_\vt)=X^\vt$ and $ X_\ast(T_\vt)=\check X/(1-\vt)\check X$ are the invariants and coinvariants of $\vt$ in $X$ and $\check X$, respectively. 

We now recall some conjugacy results from \cite{borel:corvallis} and \cite{reeder:torsion} which are stated over $\bc$ but whose proofs are unchanged if $\bc$ is replaced by any algebraically closed field $k$. 
First, we have \cite[6.4]{borel:corvallis}:
\begin{lemma}\label{nu} The natural projection $\nu:T\to T_\vt$ induces a bijection
$$T\vt/N_\vt \lra T_\vt/W^\vt,$$
sending $t\vt\mod N_\vt\mapsto \nu(t)\mod W^{\vt}$. 
\end{lemma}
From \cite[Lemma 3.2]{reeder:torsion} each semisimple element $g\vt\in G\vt$ is $G$-conjugate to an element of $t\vt$ with $t\in T^\vt$. Now \cite[6.5]{borel:corvallis} shows that sending $g\vt$ to the class of $\nu(t)$ modulo $W^\vt$ gives a bijection between the set of semisimple $G$-conjugacy classes in $G\vt$ and the orbit space $T_\vt/W^\vt$.


Now the affine variety $T_\vt/W^\vt$ has a canonical $\bz$-form,
namely the ring $\bz[X^\vt]^{W^\vt}$ of $W^{\vt}$-invariants in the integral group ring of the character group $X^\vt$ of $T_\vt$.
Indeed, let $X^{\vt}_+$ be the set of dominant weights in $X^\vt$ and for each $\lam\in X^\vt_+$,
let  $\eta_\lam$  be the sum in $\bz[X^\vt]$ over the $W^\vt$-orbit of $\lam$, and let
$\eta^k_\lam$ be the same sum in the group ring $k[X^\vt]$.
Then $\{\eta_\lam:\ \lam\in X_+^{\vt}\}$ and   $\{\eta_\lam^k:\ \lam\in X_+^{\vt}\}$
are bases of $\bz[X^\vt]^{W^\vt}$ and $k[X^\vt]^{W^\vt}$ respectively,
and $\{1\otimes\eta_\lam:\ \lam\in X^\vt_+\}$ is a $k$-basis of
$k\otimes_\bz(\bz[X^\vt]^{W^\vt})$. It follows that the canonical mapping
$\bz[X^\vt]^{W^\vt}\lra k[X^\vt]^{W^\vt}$ induces an isomorphism
\begin{equation}\label{XZform}
k\otimes_\bz(\bz[X^\vt]^{W^\vt})\overset\sim\lra k[X^\vt]^{W^\vt}.
\end{equation}

The torus $T_\vt$ is a maximal torus in a connected reductive group $G_\vt$ with Weyl group $W^\vt$,
so $\bz[X^\vt]^{W^\vt}$ has another $\bz$-basis, $\{\chi_\lam:\ \lam\in X^\vt_+\}$, where
$$\chi_\lam=\sum_{\mu\in X^\vt} m_\lam^\mu\mu,$$
and $m_\lam^\mu$ is the multiplicity of the weight $\mu$ in the irreducible representation
of highest weight $\lam$ of the complex group with the same root datum as $G_\vt$. Therefore
$k[X^\vt]^{W^\vt}$ has another $k$-basis, $\{\chi_\lam^k:\ \lam\in X^\vt_+\}$,
where $\chi_\lam^k\in k[X^\vt]^{W^\vt}$ is the image of $1\otimes \chi_\lam$ under the isomorphism \eqref{XZform}.

We now regard $G$ as a Chevalley group scheme over $\bz$, writing $G(A)$ for the group of $A$-valued points in a commutative ring $A$. 
The group heretofore denoted by $G$ is now $G(k)$. Likewise $T$ and $N$ are now group schemes over $\bz$. 

Let $\lam\in X^\vt_+$ and let $V$ be the irreducible representation of $G(\bc)$ of highest weight $\lam$. Since $\vt\lam=\lam$ it follows that 
 $V$ extends uniquely to a representation of 
$G(\bc)\cdot\la\vt\ra$ such that $\vt$ acts trivially on the highest weight space $V(\lam)$. 
 
Choose a $G(\bz)$-stable lattice $M$ in $V$ such that $M\cap V(\mu)$ 
spans each weight space $V(\mu)$ in $V$ and  $\vt M=M$. 
For example, we could take  $M$ to be the admissible $\bz$-form of $V$ constructed by Kostant in \cite{kostant:zform}. 
We get a representation of 
$G(k)\cdot\la\vt\ra$ on $V_k:=k\otimes M$ which may be reducible and which depends on $M$. However, since $M$ contains a basis of $V$, the traces on $V_k$ of elements of $G(k)\cdot\la\vt\ra$ are independent of the choice of $M$. 

Let $A=\bz[\zeta]\subset\bc$ be the cyclotomic ring generated by a root of unity 
$\zeta\in \bc^\times$ of order $m$. Assume that $k$ is algebraically closed and $m$ is nonzero in $k$.  Choose $\zeta_k\in k^\times$ a root of unity of order $m$. 
We have ring homomorphisms 
$$\bc\overset{\iota}\hookleftarrow A\overset{\pi}\lra k,$$
where $\iota$ is the inclusion and $\pi(\zeta)=\zeta_k$. We use the same letters to denote maps on groups of points, e.g., 
$$G(\bc)\overset{\iota}\hookleftarrow G(A)\overset{\pi}\lra G(k),$$
and similarly for $T$ and $N$.

\begin{lemma}\label{kconj1} Let $s,t\in T(k)^\vt$ be elements of order $m$ such that 
$\tr(s\vt, V_k)=\tr(t\vt,V_k)$ for all irreducible representations $V$ of $G(\bc)$ whose highest weight belongs to $X_+^\vt$.  Then $s\vt$ and $t\vt$ are $G(k)$-conjugate. 
\end{lemma}
\proof
Let $V'$ be the representation of $G_\vt(\bc)$ with the same highest weight as $V$. 
And choose a lattice $M'\subset V'$ analogous to $M$ above. Since $s$ has order $m$ there is a co-weight $\check\om\in \check X$ such that 
$$s=\check\om(\zeta_k)=\pi\check\om(\zeta).$$
For each $\mu\in X^\vt$ let $M(\mu)=M\cap V(\mu)$ and likewise set 
$M'(\mu)=M'\cap V'(\mu)$. We have
\begin{equation*}
\begin{split}
\tr(s\vt,V_k)
&=\sum_{\mu\in X^\vt}\mu(s)\cdot \tr(\vt,k\otimes M(\mu))
=\sum_{\mu\in X^\vt}\zeta_k^{\la\mu,\check\om\ra}\cdot 
\pi\left( \tr(\vt,M(\mu))\right)\\
&=\pi\left(\sum_{\mu\in X^\vt}\zeta^{\la\mu,\check\om\ra}\cdot 
\tr(\vt,M(\mu))\right).
\end{split}
\end{equation*}
By a result of Jantzen (see for example \cite{kumar-lusztig-prasad}) we have 
$$
\sum_{\mu\in X^\vt}\zeta^{\la\mu,\check\om\ra}\cdot \tr(\vt,M(\mu))
=
\sum_{\mu\in X^\vt}\zeta^{\la\mu,\check\om\ra}\cdot \dim M'(\mu).
$$
It follows that 
\begin{equation*}
\tr(s\vt,V_k)=\pi\left(\sum_{\mu\in X^\vt}\zeta^{\la\mu,\check\om\ra}\cdot \dim M'(\mu)\right)=\tr(\nu(s),V'_k).
\end{equation*}
Applying this identity to $t\vt$ as well, we find that 
$$\tr(\nu(s),V'_k)=\tr(\nu(t),V'_k).$$
Therefore $\chi_\lam^k(\nu(s))=\chi_\lam^k(\nu(t))$ for every $\lam\in X_+^\vt$.
Since these $\chi_\lam^k$ are a basis of $k[X^\vt]^{W^{\vt}}$, it follows from \cite[Cor. 6.6]{steinberg:regular}
that  $\nu(s)\equiv\nu(t)\mod W^\vt$.
By Lemma \ref{nu} we have that $s\vt$ and $t\vt$ are $G(k)$-conjugate, as claimed. 
\qed

Now suppose $g\in G(\bz)$ and $g\vt$ is semisimple of order $m$. 
Let $s\in T(\bc)^\vt$ and $t\in T(k)^\vt$ be such that $\iota(g)\vt$ is $G(\bc)$-conjugate to $s\vt$ and $\pi(g)\vt$ is $G(k)$-conjugate to $t\vt$. 

\begin{lemma}\label{k-conj} In the situation just described, we have $s\in T(A)$ and $\pi(s)\vt$ is $G(k)$-conjugate to $t\vt$. 
\end{lemma}
\proof
As above we have $s=\check\om(\zeta)$ for some co-weight $\check\om\in\check X$. 
It follows that $s\in T(A)$. Moreover, $g\vt$ preserves the lattice $M$, so we have
$$\tr(\iota(g)\vt,M)=\tr(s\vt, V)=\sum_{\mu\in X^\vt}\zeta^{\la \mu,\check\om\ra}\cdot\tr(\vt,M(\mu)).
$$
Applying $\pi$ to both sides we get 
\begin{equation}\label{s}
\pi\left(\tr(\iota(g)\vt,M)\right)
=\sum_{\mu\in X^\vt}\zeta_k^{\la \mu,\check\om\ra}\cdot\tr(\vt,V_k(\mu))
=\tr(\pi(s)\vt,V_k).
\end{equation}
On the other hand, we can first apply $\pi:G(A)\to G(k)$ and then take traces. 
This gives 
\begin{equation}\label{t}
\pi\left(\tr(\iota(g)\vt,M)\right)=\tr(\pi(g)\vt,V_k)=\tr(t\vt,V_k).
\end{equation}
Comparing the  expressions \eqref{s} and \eqref{t} and using Lemma \ref{kconj1} we see that 
$\pi(s)\vt$ and $t\vt$ are $G(k)$-conjugate as claimed. 
\qed

We are ready to prove Prop. \ref{clift}.
Recall that $w\vt\in W\vt$ is an elliptic $\bz$-regular automorphism of $R$  whose order $m$ is nonzero in the algebraically closed field $k$. Let $\zeta\in k^\times$  be a root of unity of order $m$. 
Recall that  $\check\rho$ is the sum of the fundamental co-weights arising from our chosen pinning.
We have $\check\rho\in\check X^\vt$  and  $\check\rho(\zeta)\in T(k)^\vt$. 
We now prove Prop. \ref{clift} in the following form. 
\begin{prop}\label{clift2}   
For any lift $n\in N(k)$ of $w$, the element
$n\vt\in G(k)\vt$ is $G(k)$-conjugate to $\check\rho(\zeta)\vt$. 
\end{prop}

\proof
Assume first that $k$ has characteristic zero. In this case the proof relies on
 \cite[Thm. 3.3]{panyushev:theta} and is similar to the proof of  
 \cite[Thm. 4.2 (iii)]{panyushev:theta}. 
The automorphism $\tau:=\check\rho(\zeta)\vt\in\Aut(\fg)$ has order $m$ and gives a grading $\fg=\oplus_{i\in\bz/m}\ \fg'_i$, where $\fg'_i$ is the $\zeta^i$-eigenspace of $\tau$. 
The sum
 $E=\sum_{i=1}^\ell E_i$ of the simple root vectors in our pinning belongs to $\fg_1'$.  
By \cite[Thm. 3.3(v)]{panyushev:theta}, the dimension of a Cartan subspace  $\fc\subset\fg_1'$ may be computed as follows. 
Let $f_1,\dots,f_\ell\in k[\ft]$ be homogeneous generators for the algebra of $W$-invariant polynomials on $\ft$. Assume, as we may, that each $f_i$ is an eigenvector for $\vt$, with eigenvalue denoted $\vep_i$, and set $d_i=\deg f_i$. 
The integer
$$a(m,\vt):=|\{i:\ 1\leq i\leq\ell,\ \vep_i\zeta^{d_i}=1\}|$$
depends only on $m$ and $\vt$, and we have 
 $$\dim\fc=a(m,\vt).$$
 Let $\fs$ be a canonical Cartan subalgebra for $\tau$ (section \ref{csa}). 
 There exists $g\in G$ such that $\ft=\Ad(g)\fs$, and we set $\theta'=g\tau g^{-1}$. 
Since $\theta'$ normalizes $\ft$ and belongs to $G\vt$ we have $\theta'\in N\vt$. 
Let $w'\vt\in W\vt$ be the projection of $\theta'$. Then $\Ad(g)\fc$ is the 
 $\zeta$-eigenspace $\ft(w'\vt,\zeta)$ of $w'\vt$ in $\ft$, so 
 $$\dim \ft(w'\vt,\zeta)=a(m,\vt).$$
 Since $w\vt$ is $\bz$-regular and therefore $k$-regular (by Prop. \ref{prop:regular}), 
 it follows from \cite[Prop. 3.6]{springer:regular} that we also have 
$\dim \ft(w\vt,\zeta)=a(m,\vt)$, and therefore
$$\dim\ft(w\vt,\zeta)=\dim\ft(w'\vt,\zeta).$$
By  \cite[Thm. 6.4 (iv)]{springer:regular} the elements  $w\vt, w'\vt\in W\vt$ are  conjugate under $W$. 
It follows that $n\vt$ is $N$-conjugate to an element of $T\theta'$. 
As $w'\vt$ is also elliptic, it follows that $n\vt$ is actually conjugate to $\theta'$, 
and hence to $\tau=\check\rho(\zeta)\vt$, as claimed.

Now assume that $k$ has positive characteristic not dividing $m$. 
Let $A$ be the cyclotomic subring of $\bc$ generated by $z=e^{2\pi i/m}$ and 
let $\pi:A\to k$ be the ring homomorphism mapping $z\mapsto \zeta$. 
By ellipticity, all lifts of $w\vt$ to $N(k)\vt$ are $T(k)$-conjugate, so we may choose our lift to be of the form $\pi(n)$ with $n\in N(\bz)$. From the characteristic zero case just proved, we have that $\iota(n)\vt$ is $G(\bc)$-conjugate to $\check\rho(z)\vt$. By Lemma \ref{k-conj} it follows that 
$\pi(n)\vt$ is $G(k)$-conjugate to $\check\rho(\zeta)\vt$, as claimed. 
\qed

\subsection{Stable gradings} \label{sec:stable}
Let $H$ be a connected reductive $k$-group  acting on a $k$-vector space $V$. 
A vector $v\in V$ is called {\bf $H$-stable} (in the sense of Geometric Invariant Theory, see 
\cite{mumford:stable}) if the $H$-orbit of $v$ is closed and the stabilizer of $v$ in $H$ is finite. 
The second condition means  that the stabilizer $H_v$ is a finite algebraic group:
it has only finitely many points over the algebraically closed field $k$.

Recall we are assuming  the characteristic  of $k$ is not a torsion prime for $G$
and that the period $m$ of the grading $\fg=\oplus_{i\in\bz/m}\fg_i$ is nonzero in $k$. 
We have chosen a  root of unity $\zeta\in k^\times$ of order $m$, and $\theta\in\Aut(\fg)$ is the automorphism of order $m$ whose $\zeta^i$-eigenspace is $\fg_i$. 

We say the grading $\fg=\oplus_{i\in \bz/m}\ \fg_i$ (or the automorphism $\theta$) 
is {\bf stable} if there are $G_0$-stable vectors in $\fg_1$. In this section we will show that stable gradings are closely related to elliptic $\bz$-regular automorphisms of the root system $R$.  

\begin{lemma}\label{stabless} A vector $v\in\fg_1$  is stable if and only if $v$ is a  regular semisimple element of $\fg$ and the action of $\theta$ on the Cartan subalgebra centralizing $v$ is elliptic. 
\end{lemma}
\proof Vinberg showed (\cite[Prop. 3]{vinberg:graded}) that the $G_0$-orbit of $v$ is closed in $\fg_1$ if and only if $v$ is semisimple in $\fg$. His proof works also in positive characteristic  (see \cite[2.12-3]{levy:thetap}). If $v$ is semisimple its centralizer $C_G(v)$ is connected (since $p$ is not a torsion prime, by \cite[Thm. 3.14]{steinberg:torsion}) and reductive with semisimple derived subgroup $H$ \cite[13.19, 14.2]{borel:linear}. As $v$ is an eigenvector for $\theta$ we have $\theta(H)=H$. 
If $v$ is stable then $H^\theta$ is finite. 
Let $\pi:H_{sc}\to H$ be the simply-connected covering of $H$. 
We lift $\theta$ to an automorphism of $H_{sc}$, denoting it again by $\theta$. 
Now $H_{sc}^\theta$ is connected \cite[chap. 8]{steinberg:endo} 
so $\pi(H_{sc}^\theta)\subset\left(H^\theta\right)^\circ$ is trivial.  
Since $\ker\pi$ is finite, we must have $H_{sc}^\theta=1$. 
This implies that $H_{sc}=1$. For otherwise, by \cite[chap. 8]{steinberg:endo}, 
there would be a maximal torus $T'$ contained in a Borel subgroup $B'$ of $H_{sc}$ such that $\theta(T')=T'$ and $\theta(B')=B'$, and $H_{sc}^\theta$ would have rank equal to the number of $\theta$-orbits on the set of simple roots of $T'$ in $B'$. Therefore $H_{sc}=1$, so $H=1$ and $C_G(v)$ is a torus. This means that $v$ is regular in $\fg$. 
The reverse implication is clear.
\qed

Prop. \ref{clift} and Lemma \ref{stabless} have the following corollaries. 

\begin{cor} \label{stable} Let $\theta\in G\vt$ have order $m$ nonzero in $k$.  The following are equivalent.
\begin{enumerate}
\item The grading on $\fg$ given by $\theta$ is stable;
\item The action of $\theta$ on its canonical Cartan subalgebra induces 
an elliptic $\bz$-regular automorphism of $R$;
\item 
$\theta$ is principal and $m$ is the order of an elliptic $\bz$-regular element of $W\vt$. 
\end{enumerate}
\end{cor}

\begin{cor}\label{classifiction} The map sending a stable automorphism $\theta\in \Aut(\fg)$ to the automorphism of $R$ induced by the action of $\theta$ on its canonical Cartan subalgebra gives a bijection between the $G$-conjugacy classes of stable automorphisms of $\fg$ and the $W$-conjugacy classes of elliptic $\bz$-regular automorphisms of $R$.
\end{cor}

\section{Affine-pinned automorphisms}\label{affine-pinned}
In this section we construct certain automorphisms of $\fg$ arising from symmetries of the affine Dynkin diagram. These will be used to study outer automorphisms of $E_6$.

Assume $\fg$ is a simple Lie algebra over $\bc$ with adjoint group $G=\Aut(\fg)^\circ$. 
Let $N, T$ be the normalizer and centralizer of a Cartan subalgebra $\ft$ of $\fg$ and let 
 $W=N/T$. 
Let $R$ be the set of roots of $T$ in $\fg$ and choose a base 
$\De=\{\al_1,\dots,\al_\ell\}$ of $R$. 
Let $\al_0$ be the lowest root of $R$ with respect to $\De$ and set $\Pi=\{\al_i:\ i\in I\}$, where $I=\{0,1,\dots,\ell\}$.  The subgroup of $W$ preserving $\Pi$, 
$$W_\Pi=\{w\in W:\ w\Pi=\Pi\}$$
is isomorphic to the fundamental group of $G$.
Each element $w\in W_\Pi$ determines a permutation $\sigma$ of $I$ such that 
$$w\cdot \al_i=\al_{\sigma (i)}.$$ 

Choose a Chevalley lattice $\fg_\bz\subset\fg$ spanned by a lattice in $\ft$ and root vectors for $T$. 
An {\it affine pinning}  is a set 
$\widetilde\Pi=\{E_0,E_1,\cdots,E_\ell\}$
consisting of nonzero root vectors $E_i\in \fg_{\al_i}\cap \fg(\bz)$ for each $i\in I$. 
Let $N(\bz)$ be the stabilizer of $\fg(\bz)$ in $N$, 
and consider the subgroup 
$$N_{\widetilde\Pi}=\{n\in N(\bz):\ n\widetilde\Pi=\widetilde\Pi\}.$$

\begin{lemma} Let  $\widetilde\Pi$ be an affine pinning. Then the projection $N\to W$ restricts to an isomorphism $f:N_{\widetilde\Pi}\overset\sim\lra W_\Pi$. 
\end{lemma}
\proof It is clear that $f(N_{\widetilde\Pi})\subset W_\Pi$. 
An element in $\ker f$  lies in $T$ and fixes each root vector $E_i$, hence lies in the center of $G$, which is trivial since $G$ is adjoint. Hence $f$ is injective. 

Let $w\in W_\Pi$. Since the projection $N\to W$ is surjective on $N(\bz)$ 
\cite[Lemma 22]{steinberg:yale}, 
there is a lift $n'$ of $w$ such that $n'\in N(\bz)$. 
For each $i\in I$ we have  
$n'\cdot E_i=c_i E_{\sigma(i)}$, for some $c_i=\pm 1$. 

Let  $\check\om_1,\dots,\check\om_\ell\in X_\ast(T)$ be the fundamental coweights of $T$ dual to 
$\al_1,\dots, \al_\ell$. 
The element $t=\prod_{i=1}^\ell \check\om_i(c_i)$ lies in $T(\bz)$ and the new lift  $n=n't$ of $w$ satisfies $n\cdot E_i=E_{\sigma( i)}$ for $1\leq i\leq \ell$. 

Let $d$ be the order of $w$. Then $\sigma^d=1$ so $n^d$ fixes $E_i$ for each $1\leq i\leq \ell$. Hence $n^d\in T$ and belongs to the kernel of each simple root $\al_i$. 
Since $G$ is adjoint, it follows that $n^d=1$. 

Let $i=\sigma(0)$. It follows from \cite[VI.2.2]{bour456} that $\sigma^j(0)\neq 0$ for $1\leq j<d$.
By what has been proved, we have 
$$n^{-1}\cdot E_i=n^{d-1}\cdot E_i=E_{\sigma^{d-1}(i)}=E_{\sigma^{-1}(i)}=E_0.$$
It follows that $n\cdot E_0=E_i$, so $n$ is a lift of $w$ in $N_{\widetilde\Pi}$. 
\qed

Now let $k$ be an algebraically closed field of characteristic not equal to two, 
and view $G$ as a group scheme over $\bz$, via the lattice $\fg_\bz$. 
Take  $w\in W_\Pi$ of order two. Again from \cite[VI.2.2]{bour456} there exists a unique minuscule coweight $\check\om_j$ such that $w\check\om_j=-\check\om_j$. 
Since $2\neq 0$ in $k$, the natural map $T(\bz)\to T(k)$ is injective, which implies that
the map $N(\bz)\to N(k)$ is injective. We now let $n$ be the image in $N(k)$ of 
the unique lift of $w$ in $N_{\widetilde\Pi}$. 
\begin{prop}\label{kacminuscule}
There exists an affine pinning $\widetilde\Pi$ such that $n$ is $G(k)$-conjugate to $\check \om_j(-1)$. The Kac coordinates of $\Ad(n)$ are given by:
$$
s_i=
\begin{cases} 
1&\quad\text{for}\quad i\in\{0,j\}\\
0&\quad\text{for}\quad i\notin\{0,j\}.
\end{cases}
$$
These labels give the unique $w$-invariant Kac-diagram of order two having $s_0\neq 0$.
\end{prop}
\proof
By \cite[Lemma 5]{carter:weyl} there are mutually orthogonal roots $\gamma_1,\dots,\gamma_m\in R$ with corresponding reflections $r_1,\dots, r_m\in W$, such that 
\begin{equation}\label{product}
w=r_1r_2\cdots r_m.
\end{equation}
Since $\check \om_j$ is minuscule we have $\la \al,\check\om_j\ra\in\{-1,0,1\}$ 
for each $\al\in R$. The positive roots made negative by $w$ are those for which 
$\la \al,\check\om_j\ra\neq 0$. Since $w\gamma_i=-\gamma_i$ for each $i$, 
we may choose the sign of each $\gamma_i$ so that $\la \gamma_i,\check \om_j\ra=1$. And since
$$-\check \om_j=w\cdot\check \om_j=\check \om_j-\sum_{i=1}^m\la \ga_i,\check \om_j\ra\check \ga_i,$$
it then follows that 
\begin{equation}\label{sum}
\check \gamma_1+\check \gamma_2+\cdots +\check \gamma_m=2\check \om_j.
\end{equation}
For each $i=1,\dots,m$  there exists a morphism 
$\vp_i: SL_2\to G$ 
over $\bz$ whose restriction to the diagonal subgroup is given by
$$\vp_i\left(\begin{bmatrix}t&0\\0&t^{-1}\end{bmatrix}\right)=\check \gamma_i(t)$$ 
and such that $\vp_i\left(\begin{bmatrix}0&-1\\1&0\end{bmatrix}\right)\in N(\bz)$ and is a representative of $r_i$. 

Since the roots $\gamma_i$ are mutually orthogonal, the images of these homomorphisms $\vp_i$ commute with one another. Hence we have a $\bz$-morphism
$$\vp:SL_2\to G,\qquad \text{given by}\qquad 
\vp\left(\begin{bmatrix}a&b\\c&d\end{bmatrix}\right)=\prod_{i=1}^m\vp_i\left(\begin{bmatrix}a&b\\c&d\end{bmatrix}\right).
$$
By equation  \eqref{product} the element 
\begin{equation}\label{npinned}
n:=\vp\left(\begin{bmatrix}0&-1\\1&0\end{bmatrix}\right)
\end{equation}
belongs to $N(\bz)$ and represents $w$. 
Equation \eqref{sum} implies that 
$$\vp\left(\begin{bmatrix}t&0\\0&t^{-1}\end{bmatrix}\right)=\check \om_j(t)^2,$$ 
which in turn implies that $n$ has order two. 
Since the matrices 
$\begin{bmatrix}0&-1\\1&0\end{bmatrix}$ and 
$\begin{bmatrix}\sqrt{-1}&0\\0&-\sqrt{-1}\end{bmatrix}$ 
 are conjugate in $\SL_2$, it follows that $n$ is conjugate  to $\check \om_j(-1)$ in $G$, 
 and that $\Ad(n)$ has the asserted Kac-coordinates.
  
 We construct an affine pinning stable under $n$ as follows. 
Choose representatives $\al_i$ of the $w$-orbits in $\Pi$, and choose arbitrary nonzero root vectors $E_i\in \fg(\bz)$ for these roots. Let $\sigma$ be the permutation of $I$ induced by $w$. If $w\cdot\al_i\neq\al_i$, let $E_{\sigma(i)}=n\cdot E_i$. Since $n$ has order two, we have $n\cdot E_{\sigma(i)}=E_i$. If $w\cdot\al_i=\al_i$ then $\al_i$ is orthogonal to each of the roots $\ga_1,\dots,\ga_m$, since the latter are negated by $w$. It follows that the image of each homomorphism $\vp_1,\dots,\vp_m$ centralizes the root space $\fg_{\al_i}$, so  any nonzero vector $E_i\in\fg_{\al_s}\cap\fg(\bz)$ is fixed by $n$. 
The collection $\widetilde\Pi=\{E_i\}$ of vectors thus defined is an affine pinning stable under $n$. \qed

The following lemma will also be useful. 

\begin{lemma}\label{nfixed} Let $S=(T^n)^\circ$ be the identity component of the subgroup of $T$ centralized by $n$. Then $S$ is centralized by the entire group $\vp(\SL_2)$. 
\end{lemma}
\proof
Since $2\check\om_j$ is a simple co-weight in $\vp(\SL_2)$ and $\check\om_j$ is minuscule, we have that $\la \al,2\check\om_j\ra\in \{-2,0,2\}$ for every root $\al\in R$. 
Hence $\vp(\SL_2)$ acts on $\fg$ as a sum of copies of the trivial and adjoint representations. Applying the element 
\[\begin{bmatrix} 1&0\\-t&1\end{bmatrix}\cdot
\begin{bmatrix} 0&-1/t\\t&0\end{bmatrix}=
\begin{bmatrix} 1&-1/t\\0&1\end{bmatrix}\cdot
\begin{bmatrix} 1&0\\t&1\end{bmatrix}
\]
to a vector in the zero weight space and comparing components in the $-2$ weight space, we find (since the characteristic of $k$ is not two) that any vector in $\fg$ invariant under the normalizer of 
$2\check\om_j(k^\times)$ in $\vp(\SL_2)$ is invariant under all of $\vp(\SL_2)$. 
Since the Lie algebra of $S$ consists of such vectors, the lemma is proved. 

\qed

\section{Little Weyl groups}\label{littleweyl}

Let $\theta$ be an automorphism of $\fg$ whose order $m$ is invertible in $k$. 
Choose a root of unity $\zeta\in k^\times$ of order $m$ and let 
 $\fg=\oplus_{i\in\bz/m}\ \fg_i$ be the grading of $\fg$ into $\zeta^i$-eigenspaces of $\theta$. 
Choose a Cartan subspace $\fc$ in  $\fg_1$ and assume the rank 
 $r=\dim \fc$ is positive. 
The little Weyl group is defined as
$$W(\fc,\theta)=N_{G_0}(\fc)/Z_{G_0}(\fc),$$
where $G_0=(G^\theta)^\circ$ is the connected subgroup of $G$ with Lie algebra $\fg_0$.  When it is necessary to specify $G$ in the little Weyl group we will write $W_G(\fc,\theta)$. 

It is clear from the definition that $W(\fc,\theta)$ acts faithfully on $\fc$. 
From \cite{vinberg:graded} and \cite{levy:thetap}, 
it is known that the action of $W(\fc,\theta)$ on $\fc$ is generated by transformations fixing a hyperplane in $\fc$,  that
the restriction map $k[\fg_1]^{G_0}\to k[\fc]^{W(\fc,\theta)}$ is an isomorphism, 
and that this ring is a polynomial ring with homogeneous generators $f_1,\dots,f_r$, 
such that 
$$|W(\fc,\theta)|=\prod_{i=1}^r\deg(f_i).$$

\subsection{Upper bounds on the little Weyl group}
Recall we have fixed a Cartan subalgebra $\ft$ in $\fg$, with normalizer and centralizer $N$ and $T$ in $G$ and we have identified $W=N/T$. 

Replacing $\theta$ by a $G$-conjugate if necessary, we may assume $\ft$ is the canonical Cartan subalgebra for $\theta$ (see \ref{csa}).
In particular $\fc$ is the $\zeta$-eigenspace of $\theta$ in $\ft$. 
Then $\theta$ normalizes $N$ and $T$ in $\Aut(\fg)$, giving an action of $\theta$ on $W$; 
let $W^\theta=\{y\in W:\ \theta(y)=y\}$ be the fixed point subgroup of $\theta$ in $W$.

Elements in $W^\theta$ commute with the action of $\theta$ on $\ft$, so $W^\theta$ acts on the eigenspace $\fc$. 
Let 
\begin{equation}\label{W1}
W_1^\theta:=W^\theta/C_W(\fc)^\theta
\end{equation}
be the quotient acting faithfully on $\fc$. 
Since $\ft$ is a Cartan subalgebra in the Levi subalgebra $\fm=\fz_\fg(\fc)$, it follows that 
every element of $W(\fc, \theta)$ has  a representative in $N$ and 
that $W(\fc,\theta)$ may be viewed as a subgroup of $W_1^\theta$. Thus, 
we have an embedding 
$$W(\fc,\theta)\hra W_1^\theta.$$
Note that $W(\fc,\theta)$  is more subtle than $W_1^\theta$.  For it can happen that two automorphisms $\theta$ and $\theta'$ of the same order  agree  on $\ft$ and $W$, so they have the same Cartan subspace $\fc$ and $W_1^{\theta}=W_1^{\theta'}$, but nevertheless $W(\theta,\fc)\neq W(\theta',\fc)$ (e.g. cases $4_a$ and $4_b$ in $E_6$; these examples are also used in 
\cite[4.5]{panyushev:theta} to illustrate other subtleties). 

A still coarser group,  depending  only on $\fc$ and not on $\theta$ is 
$$W(\fc):= N_{W}(\fc)/C_{W}(\fc). 
$$
As subgroups of $\GL(\fc)$, we have containments
$$W(\fc,\theta)\subset W_1^\theta\subset W(\fc).$$
Under certain circumstances one or both of these containments is an equality. 

\begin{lemma}\label{littleweylregular} Suppose $\fc$ contains a regular element of $\fg$. Then 
$$W_1^\theta=W^\theta=W(\fc).$$
\end{lemma}
\proof By regularity it is clear that $W_1^\theta=W^\theta$ and that $W(\fc)=N_W(\fc)$. 
And any $y\in N_W(\fc)$ commutes with the scalar action of $\theta$ on $\fc$ so the commutator $[y,\theta]$ is trivial in $W$, again by regularity. 

\qed

Panyushev \cite[Thm. 4.7]{panyushev:theta} has shown that both containments above are equalities if $\theta$ is principal:

\begin{prop}[Panyushev]\label{pan} If $\theta$ is principal  then 
$W(\fc,\theta)=W_1^\theta=W(\fc).$
\end{prop}
We note that Panyushev works in characteristic zero, but his geometric proof works equally well in good characteristic $p\nmid m$, using the invariant theoretic results of \cite{levy:thetap}. 

\begin{cor}\label{regprinc} If $\theta$ is principal and the restriction of $\theta$ to $\ft$ induces a $\bz$-regular automorphism of $R$ then $W(\fc,\theta)=W^\theta$. 
\end{cor}
\proof By Prop. \ref{prop:regular}, $\bz$-regularity implies $k$-regularity, so $W(\fc)=W^\theta_1$ is just $W^\theta$.  \qed

This sharpens the first result in this direction, which was proved in  Vinberg's original work
\cite[Prop. 19]{vinberg:graded}:
\begin{cor}[Vinberg]\label{vinberg:stable} If $\theta$ gives a stable grading on $\fg$ then 
$W(\fc,\theta)=W^\theta$. 
\end{cor}

\subsection{Little Weyl groups for inner gradings}
Assume now that $\theta$ is inner, and let the restriction 
of $\theta$ to $\ft$ be given by the element $w\in W$. 
In this section we give upper and lower bounds for $W(\fc,\theta)$ depending only on $w$, 
under certain conditions; these will suffice to compute almost all little Weyl groups in type $E_n$. 
The fixed-point group
$$W^\theta=C_W(w),$$
is now the centralizer of $w$ in $W$, which acts on the $\zeta$-eigenspace  $\fc$
of $w$ in $\ft$. 
The quotient by the kernel of this action is the group $W_1^\theta$. 
Simple upper and lower bounds for  $W(\fc,\theta)$ can be obtained as follows. 

\begin{lemma}\label{simplebound} If $U$ is any subgroup of $C_W(w)$ acting trivially on $\fc$ then we have the inequalities
$$m\leq |W(\fc,\theta)|\leq \frac{|C_W(w)|}{|U|}.$$
\end{lemma}
\proof
Since $\theta$ is semisimple it lies in the identity component $G_0$ of its centralizer in $G$. Hence the cyclic group $\la\theta\ra$ embeds in $W(\fc,\theta)$, whence the lower bound. 
The upper bound follows from \eqref{W1}.
\qed

Information about $C_W(w)$, including its order, is given in \cite{carter:weyl}. 
Using the tables therein, one can often find a fairly large subgroup $U\subset C_W(w)$ 
as in Lemma \ref{simplebound}.

{\bf Example 1:\ } In type $E_8$ there are eight cases  (namely $12_b$ through $12_i$ in the tables below) where $w$ is a Coxeter element in 
$W(E_6)$.  
From \cite{carter:weyl}  we have $|C_W(w)|=144$. Hence the centralizer is given by 
$$C_W(w)=\la w\ra\times \la -w^6\ra\times W(A_2),$$
where $A_2$ is orthogonal to the $E_6$. 
Since $\fc$ lives in the $E_6$ Levi subalgebra and $w^6$ acts by $-1$ on $\fc$, 
the inequalities of Lemma \ref{simplebound} become equalities for 
 $U=\la -w^6\ra\times W(A_2)$. 
Hence $W(\fc,\theta)\simeq \mu_{12}$ in these eight cases. 

{\bf Example 2:\ } In type $E_8$ there are four cases ($6_h$ through $6_k$) where $w$ is a Coxeter element in  $W(D_4)$. Let  $\De_4=\{\be_1,\dots,\be_4\}$  be a base of the corresponding root subsystem of type $D_4$. The subgroup of $W(E_8)$ permuting $\De_4$ is a symmetric group $S_3$. We may choose the Coxeter element $w$ to be centralized by this $S_3$, and $\fc$ is a line in the span of the co-root vectors $\{d\check\be_i(1)\}$. 
 The roots of $E_8$ orthogonal to $\De_4$ form another system of type $D_4$, 
 hence there is a subgroup $W_2\simeq W(D_4)$ fixing each root in $\De_4$ and therefore 
 acting trivially on $\fc$. Since $S_3$ normalizes $\De_4$ it also normalizes $W_2$. 
From \cite{carter:weyl}  we have $|C_W(w)|=6\cdot 6\cdot 192$, 
so the inequalities of Lemma \ref{simplebound} hold for $U\simeq S_3\ltimes W(D_4)$. 
Hence $W(\fc,\theta)\simeq \mu_{6}$ in these four cases.

{\bf Example 3:\ } In type $E_7$ there are two cases ($9_a$ and $9_b$) 
where $w$ is the square of a Coxeter element 
and we have $C_W(w)=\la -w\ra\simeq\mu_{18}$. 
Since $w$ is $\bz$-regular,  Lemma \ref{simplebound} only gives the inequalities
$$9\leq |W(\fc,\theta)|\leq 18.$$
In fact, we have  $W(\fc,\theta)\simeq\mu_{18}$ and $\mu_9$  in cases $9_a$ and $9_b$, respectively. This shows that, in general, 
$W(\fc,\theta)$ depends on $\theta$, and not just on $w$. 
We will return to this example after sharpening our lower bound, as follows. 

For any subset $J\subset\{1,\dots,\ell\}$ let $R_J$ be the root subsystem generated by 
$\{\al_j:\ j\in J\}$, let $W_J$ be Weyl group of $R_J$ and let  $\fg_J$ be the subalgebra of $\fg$ generated by the root spaces $\fg_\al$ for $\al\in R_J$. If the action of $\theta$ on $\ft$ is given by an element $w\in W_J$ then $\theta$ induces an automorphism $\theta_J$ of $\fg_J$.

\begin{lemma}\label{Jlowerbound} Suppose $\theta$ normalizes the Cartan subalgebra $\ft$ and has image 
$w\in W_J$ for some  subset $J\subset\{1,\dots,\ell\}$ 
such that the following conditions hold. 
\begin{enumerate}
\item $\theta$ is conjugate to an automorphism $\theta'=\Ad(t)$ where $t\in T$ satisfies $\al_j(t)=\zeta$ for all $j\in J$; 
\item The rank of $w$ on $\ft$ is equal to the rank of $\theta$; 
\item The principal automorphisms of $\fg_J$ of order $m$ have rank equal to the rank of $\theta$. 
\item $w$ is $\bz$-regular in $W_J$;
\end{enumerate}
Then there is an embedding $C_{W_J}(w)\hra W(\fc,\theta)$.
\end{lemma}
\proof 
Condition 1 means there is $g\in G$ such that the automorphism 
$$\theta'=g\theta g^{-1}=\Ad(t),$$
where $t\in T$ satisfies $\al_j(t)=\zeta$ for all $j\in J$. 
We have $t=\check\rho_J(\zeta)z$ where $\check\rho_J$ is half the sum of the positive co-roots of $R_J$ (with respect to $\De_J$) and $z\in \ker\al_j$ for all $j\in J$. 

Condition 2 means that the eigenspace $\fc:=\ft(w,\zeta)$ is a Cartan subspace for $\theta$. 
Note that $\fc\subset \fg_J$. 
Let $\fc_J$ be a Cartan subspace for the automorphism 
$$\theta'_J:=\theta'\vert_{\fg_J}=\Ad(\check\rho_J(\zeta))\in G_J,$$
where $G_J=\Aut(\fg_J)^\circ$. 

As $\theta'_J$ is principal of order $m$, we have $\dim\fc_J=\dim\fc$, by condition 3. 

Now $\fc':=\Ad(g)\fc$ is a Cartan subspace for $\theta'$ in 
$\fg(\theta',\zeta)$, and the latter subspace contains 
$\fg_J(\theta'_J,\zeta)$, which in turn contains $\fc_J$. 
Thus $\fc'$ and $\fc_J$ are two Cartan subspaces in $\fg(\theta',\zeta)$, 
so there is $h\in G^{\theta'}$ such that $\Ad(hg)\fc=\Ad(h)\fc'=\fc_J$ 
\cite[Thm. 2.5]{levy:thetap}. 
Conjugation by $hg$ gives an isomorphism 
$$W_G(\fc,\theta)\overset\sim\lra W_G(\fc_J,\theta').$$
Since the latter group contains $W_{G_J}(\fc_J,\theta'_J)$, we have an embedding 
$$W_{G_J}(\fc_J,\theta'_J)\hra W_G(\fc,\theta).$$

Let $\ft_J=\ft\cap\fg_J$ and let $\ft_J'$ be a $\theta'_J$-stable Cartan subalgebra of $\fg_J$ containing $\fc_J$. 
Then there is $b\in G_J$ such that $\Ad(b)\ft_J'\subset \ft_J$, 
so $b\theta'_Jb^{-1}$ normalizes $\ft_J$ and $\fc_J':=\Ad(b)\fc_J$ is a Cartan subspace for 
$b\theta'_Jb^{-1}$ contained in $\ft_J$. 
Let $w'\in W_J$ be the element induced by $b\theta'_Jb^{-1}$. 
We now have two elements $w,w'\in W_J$ having equidimensional $\zeta$-eigenspaces $\fc$ and $\fc_J'$ in $\ft_J$.

The one-parameter subgroups of $G_J$ which centralize $\ft_J$ form a lattice giving a
$\bz$-form $\check X_J$  of $\ft_J$. Let $A$ be the cyclotomic subring of $\bc$  generated by $z=e^{2pi i/m}$ and let $\pi:A\to k$ be the ring homomorphism sending $z\mapsto \zeta$. 
Since the map $\pi:\mu_m(\bc^\times)\to \mu_m(k^\times)$ is an isomorphism, 
it follows that the $z$-eigenspaces of $w$ and $w'$ in $\check X_J\otimes\bc$ have the same dimension. 

Now  $w$ is $k$-regular on $\ft_J=k\otimes\check X_J$, by condition 4. 
Hence $w$ is $\bc$-regular on $\bc\otimes\check X_J$, by Prop. \ref{prop:regular}. 
By \cite[6.4]{springer:regular}, the elements $w$ and $w'$ are conjugate in $W_J$, so $w'$ is $k$-regular on 
$\ft_J$. Hence the principal automorphism  $b\theta'_Jb^{-1}$ of $\fg_J$ has regular vectors in 
$\Ad(b)\fc_J$, so the principal automorphism $\theta'_J$ has regular vectors in $\fc_J$. 
It now follows from Cor. \ref{regprinc} that 
$W_{G_J}(\fc_J,\theta'_J)\simeq C_{W_J}(w')\simeq C_{W_J}(w)$. 
\qed

{\bf Remarks:\ } 1.\ In practice, condition 1 means the  normalized Kac diagram of $\theta$ can be conjugated under the affine Weyl group $W_{\aff}(R)$ to a (usually un-normalized) Kac diagram with $1$ on each node for $j\in J$. We will see that condition 1 is verified as a byproduct of the normalization algorithm.  

2.\ The element $w$ is usually elliptic in $W_J$. When this holds, condition 3 is implied by conditions 2 and 4, as follows  from Prop. \ref{clift2}. 

3.\ Recall that the order of $C_{W_J}(w)$ is the product of those degrees of $W_J$ which are divisible by the order $m$ of $w$. Thus the lower bound in Prop. \ref{Jlowerbound} is completely explicit. 

{\bf Example 3 revisited:\ } Recall that $G$ has type $E_7$ and $w$ is the square of a Coxeter element. We give the normalized Kac diagram for each $\theta$, the un-normalized diagram for each $\theta'$, whose subdiagram of $1$'s determines $J$. 
$$
\begin{array}{cccc}
\text{No.} & \theta& \theta' & J\\
\hline
9_a&\EVII{0}{ 1}{ 0}{ 0}{ 1}{ 0}{ 1}{ 1} &\EVII{-8}{1}{1}{1}{1}{1}{1}{1}  & E_7\\
9_b&\EVII{1}{ 0}{ 1}{ 1}{ 0}{ 0}{ 1}{ 1}   & \EVII{-7}{1}{1}{1}{1}{1}{1}{0}& E_6\\
\hline
\end{array}
$$
Lemma  \ref {Jlowerbound} shows that $9_a$ has little Weyl group $W(\fc,\theta)\simeq\mu_{18}$, but does not decide case $9_b$, which we treat using invariant theory (see section \ref{littleweylE}).

\subsection{Stable isotropy groups}\label{stableisotropy}
Assume that $\theta\in\Aut(\fg)$ gives a stable grading $\fg=\oplus_{i\in\bz/m}\ \fg_i$. By definition there is a regular semisimple element $v\in\fg_1$ whose isotropy subgroup in $G_0$ is finite. Fix a Cartan subspace $\fc\subset\fg_1$ and let $S$ be the unique maximal torus in $G$ centralizing $\fc$. In the proof of Lemma \ref{stabless} we saw that $C_G(v)$ is a torus, so we must have $C_G(v)=S$.  It follows that  all stable vectors in $\fc$ have the same isotropy group in $G_0$, equal to
$$S_0:=S\cap G_0.$$
We now give a more explicit description of $S_0$. 

First,  $S_0$ is contained in the fixed-point subgroup $S^\theta$, which is finite of order 
$$|S^\theta|=\det(1-\theta\vert_{X^\ast(S)}).$$
Let $N(S)$ be the normalizer of $S$ in $G$. 
Then $N(S)^\theta$ meets all components of $G^\theta$, and it follows from Cor. \ref{vinberg:stable}  that 
the inclusion $S^\theta\hra G^\theta$ induces an isomorphism 
$$S^\theta/S_0\simeq G^\theta/G_0.$$
This quotient depends only on the image $\vt$ of $\theta$ in the component group of $\Aut(\fg)$. To see this, let
$$G_{sc}\overset\pi\lra G$$ 
be the simply-connected covering of $G$ and set $Z=\ker\pi$. 
Then $\theta$ and $\vt$ lift to automorphisms of $G_{sc}$
which we again denote by $\theta$ and $\vt$.
Since $G_{sc}^\theta$ is connected and $\theta=\vt$ on $Z$, 
we have an exact sequence 
$$1\lra Z^\vt\lra G_{sc}^\theta\lra G_0\lra 1,$$
which restricts to an exact sequence 
$$1\lra Z^\vt\lra S_{sc}^\theta\lra S_0\lra 1,$$
where $S_{sc}=\pi^{-1}(S)$. Since 
$$|S^\theta|=|S_{sc}^\theta|,$$
it follows that we have another exact sequence 
$$1\lra S_0\lra S^\theta\lra Z/(1-\vt)Z\lra 1.$$
On the other hand, $Z/(1-\vt)Z$ is isomorphic to the subgroup  $\Om_\vt\subset \wtW_{\aff}(R,\vt)$ stabilizing the alcove $C$. 
 The group $\Om_\vt$ acts as symmetries of the twisted affine Dynkin diagram 
 $D({^eR})$. These groups are well-known if $e=1$; for $e>1$, $\Om_\vt$ is the full symmetry group of $D({^eR})$ and has order $1$ or $2$. 
 It follows that if $\theta$ is stable then 
 the isotropy group $S_0$ fits into an exact sequence 
 \begin{equation}\label{isotropy2}
 1\lra S_0\lra S^\theta\lra \Om_\vt\lra 1.
 \end{equation}
 The groups $S_0$ are tabulated for exceptional groups in Sect. \ref{exceptional}.
 
\subsection{Stable orbits and elliptic curves}

Certain remarkable stable gradings have appeared in recent work of Barghava and Shankar on the average rank of elliptic curves 
(\cite{bhargava-shankar1}, \cite{bhargava-shankar2}). These gradings have periods $m=2,3,4,5$ and are of types ${^2\!A_2},\ {^3\!D_4},\ {^2\!E_6},\ E_8$ respectively, as tabulated below. Here  $\mathbf d$ stands for the natural representation of $\SL_d$. 

\begin{center}
$$
{\renewcommand{\arraystretch}{1.2}
\begin{array}{cccccc}
\hline
m&\text{Kac coord.}& W(\fc,\theta)&\text{degrees}&G_0&\fg_1\\
\hline 
&&&&&\\
2&\twoAtwoo &\SL_2(\bz/2)&2,3
&\SL_2/\bfmu_2&\Sym^4(\mathbf 2)\\
&&&&&\\
 
3&0\ 0\Lleftarrow 1&\SL_2(\bz/3)&4,6
&\SL_3/\bfmu_3&\Sym^3(\mathbf 3)\\
&&&&&\\
 
4&0\ 0\ 0\Leftarrow 1\ 0  &\bfmu_2\times\SL_2(\bz/4)&8,12
&(\SL_2\times\SL_4)/\bfmu_4&\mathbf 2\boxtimes\Sym^2(\mathbf 4)\\
 &&&&&\\
 
5&\E{0}{ 0}{ 0}{ 0}{ 0}{ 1}{ 0}{ 0}{ 0} &\bfmu_5\times\SL_2(\bz/5)&20,30
&(\SL_5\times\SL_5)/\bfmu_5&\mathbf 5\boxtimes\Lambda^2\mathbf 5\\
&&&&&\\
\hline
\end{array}}
$$
\end{center}
For each $m=2,3,4,5$ the isotropy subgroup $S_0$ is isomorphic to  $\bfmu_m\times\bfmu_m$ and the little Weyl group  $W(\fc,\theta)$ is isomorphic to the group $W_m$ with presentation 
$$W_m=\la s,t:\ s^m=t^m=1, \quad sts=tst\ra.$$
(Note that $W_m$ is infinite  for $m> 5$.) The exact sequence 
$$
1\lra S_0\lra N_{G_0}(\fc)\lra W(\fc,\theta)\lra 1
$$
gives a homomorphism $W(\fc,\theta)\to \Aut(S_0)=\GL_2(\bz/m\bz)$ 
with image $\SL_2(\bz/m\bz)$ and split kernel $\la\theta^e\ra\simeq \bfmu_{m/e}$, 
as tabulated above (see also \cite{reeder:cyclotomic}). 

In each case the number $|R|$ of roots is equal to $m\cdot (m-1)\cdot (12/b),$
where $b=4,3,2,1$ is the maximal number of bonds between two nodes in the twisted affine diagram $D({^eR})$. 
We have $\dim G_0=|R|/m$ and 
the degrees $d_1<d_2$ have the property that $3d_1=2d_2=|R|/(m-1)$. 
Let $I, J\in k[\fc]^{W(\fc,\theta)}$ be homogeneous generators of degrees $d_1, d_2$. 
The discriminant on $\ft$ (product of all the roots in $R$) has restriction to $\fc$ given by 
$D^{m-1}$ (up to nonzero scalar),  where $D=-4I^3-27J^2$. The stable vectors $v\in \fc$ are those where $D(v)\neq 0$, and each stable vector $v$ corresponds to an elliptic curve 
$E_v$ with equation
$$y^2=x^3+I(v)\cdot x+J(v)$$
whose $m$-torsion group $E_v[m]$ is isomorphic (as an algebraic group over $k$) to $S_0$. 
For more information, along with some generalizations to hyperelliptic curves, see \cite{gross:manjulrep}. 

\section{Classification of stable gradings}\label{stableclassification}

Let $\theta\in G\vt$ be an automorphism of $\fg$ whose order $m$ is invertible in $k$, 
associated to the grading $\fg=\oplus_{i\in\bz/m}\ \fg_i$. 
After conjugating $\theta$ by an element of $G$ we may assume that $\ft$ is the canonical Cartan subalgebra of $\theta$. Then  $\theta\vert_\ft=w\vt$, for some $w\in W$. 
In section \ref{principalstable} we have seen that $\theta$ is stable if and only if $w\vt$ is an elliptic $\bz$-regular automorphism of $R$ , in which case $\theta$ is $G$-conjugate to $\check\rho(\zeta)\vt$ for some/any root of unity $\zeta\in k^\times$ of order $m$. Moreover,  the $G$-conjugacy class of $\theta$ is completely determined by its order $m$. 
The values of $m$ which can arise  are the orders of elliptic $\bz$-regular automorphisms of $R$ in $W\vt$; these are classified in \cite{springer:regular}. 

For example, the elliptic $\bz$-regular elements in $W\vt$ of maximal order are the 
{\bf $\vt$-Coxeter elements}, whose order is the 
$\vt$-Coxeter number 
$$h_\vt=e\cdot (b_1+b_2+\cdots+b_{\ell_\vt})$$
(see \eqref{twistedcoxeter}). 
These form a single $W$-conjugacy class in $W\vt$, 
representatives of which include elements of the form $w\vt$, where 
$w$ is the product, in any order, of one reflection $r_i$ taken from each of the $\vt$-orbits on simple reflections.  

For any algebraically closed field $k$ in which $h_\vt$ is invertible 
and any $\zeta\in k^\times$ of order $h_\vt$, 
the automorphism 
$$\theta_{\cox}=\Ad(\check\rho(\zeta))\vt\in\Aut(\fg)$$
is stable of order $h_\vt$ and acts on its canonical Cartan subalgebra via a $\vt$-Coxeter element. The Kac coordinates  of $\theta_\cox$ have 
$s_i=1$ for all $i\in\{0,\dots,\ell_\vt\}$ and are already normalized. 

For $m<h_\vt$ the automorphism $\check\rho(\zeta)\vt$ corresponds to a point in $\sca_\bq^\vt$ with un-normalized coordinates $s_i=1$ for $i\neq 0$ and $s_0=1+(m-h_\vt)/e$ (see \ref{principalmu}). Here we must apply the normalization algorithm to obtain normalized Kac coordinates.  By \eqref{isotropy2} these normalized Kac diagrams will be invariant under the symmetry group of the diagram 
$D({^eR})$. 
The resulting classification of the stable gradings in all types is tabulated for exceptional Lie algebras in section \ref{exceptional} and for classical Lie algebras in section 
\ref{classical}.

\subsection{Stable gradings of exceptional Lie algebras}
\label{exceptional}
Here we tabulate the stable gradings for exceptional Lie algebras, 
along with the corresponding elliptic $\bz$-regular element $w\vt\in W\vt$ 
and the isotropy group $S_0$ (see section \ref{stableisotropy}). 
The column labelled $A$ will be explained in section \ref{distinguished}.

\begin{center}
{\small Table 2: The stable gradings  for ${E_6}$}
$$
{\renewcommand{\arraystretch}{1.2}
\begin{array}{cccccc}
\hline
m&\text{un-normalized}&\text{normalized}& w&S_0&A\\
\hline 
12=h_\vt&\EVI{1}{1}{1}{1}{1}{1}{1} &\EVI{1}{1}{1}{1}{1}{1}{1}&E_6
&1&E_6\\
&&&&&\\
 
9&\EVI{-2}{1}{1}{1}{1}{1}{1}&\EVI{1}{1}{1}{1}{0}{1}{1} &E_6(a_1) 
&1&E_6(a_1)\\
&&&&&\\
 
6&\EVI{-5}{1}{1}{1}{1}{1}{1}&\EVI{1}{ 1}{ 0}{ 0}{ 1}{ 0}{ 1}  &E_6(a_2)
&1&E_6(a_3)\\
 &&&&&\\
3&\EVI{-8}{1}{1}{1}{1}{1}{1}&\EVI{0}{ 0}{ 0}{ 0}{ 1}{ 0}{ 0} &3A_2 
&\bfmu_3\times\bfmu_3&-\\
\hline
\end{array}}
$$
\end{center}

\vskip25pt
\begin{center}
{\small Table 3: The stable gradings  for ${^2\!E_6}$}
$$
{\renewcommand{\arraystretch}{1.3}
\begin{array}{ccccc}
\hline
m&\text{un-normalized}&\text{normalized}& w\vt &S_0\\
\hline
18=h_\vt&\outEVI{1}{1}{1}{1}{1} &\outEVI{1}{1}{1}{1}{1}&-E_6(a_1)
&1 \\
 
12&\outEVI{-\!2}{1}{1}{1}{1} &\outEVI{1}{1}{0}{1}{1}&-E_6 
&1 \\
 
6&\outEVI{-\!5}{1}{1}{1}{1}&\outEVI{1}{0}{0}{1}{0}&-(3A_2)
&1  \\
 4&\outEVI{-\!6}{1}{1}{1}{1}&\outEVI{0}{0}{0}{1}{0}&-D_4(a_1) 
&\bfmu_4\times\bfmu_4\\
 
2&\outEVI{-\!7}{1}{1}{1}{1} &\outEVI{0}{0}{0}{0}{1}&-1
&\bfmu_2^6 \\
\hline
\end{array}}
$$
\end{center}

\vskip25pt
\begin{center}
{\small Table 4: The stable gradings  for ${E_7}$}
$$
{\renewcommand{\arraystretch}{1.3}
\begin{array}{cccccc}
\hline
m&\text{un-normalized}&\text{normalized}& w&S_0&A\\
\hline 
18=h_\vt&\EVII{1}{1}{1}{1}{1}{1}{1}{1}  &\EVII{1}{1}{1}{1}{1}{1}{1}{1} &E_7
&1&E_7\\
 
14&\EVII{-3}{1}{1}{1}{1}{1}{1}{1} &\EVII{1}{1}{1}{1}{0}{1}{1}{1}&E_7(a_1) 
&1&E_7(a_1)\\
 
6&\EVII{-11}{1}{1}{1}{1}{1}{1}{1}&\EVII{1}{ 0}{ 0}{ 0}{ 1}{ 0}{ 0}{ 1} &E_7(a_4)
&1 &E_7(a_5)\\
 
2&\EVII{-15}{1}{1}{1}{1}{1}{1}{1} &\EVII{0}{ 0}{ 1}{ 0}{ 0}{ 0}{ 0}{ 0}&7A_1
&\bfmu_2^6&-\\
\hline
\end{array}}
$$
\end{center}

\begin{center}
{\small Table 5: The stable gradings  for ${E_8}$}
$$
{\renewcommand{\arraystretch}{1.3}
\begin{array}{cccccc}
\hline
m&\text{un-normalized}&\text{normalized}& w&S_0&A\\
\hline
30=h_\vt&\E{1}{ 1}{ 1}{ 1}{ 1}{ 1}{ 1}{ 1}{ 1}\quad &\E{1}{ 1}{ 1}{ 1}{ 1}{ 1}{ 1}{ 1}{ 1} &E_8
&1&E_8\\

24&\E{-5}{ 1}{ 1}{ 1}{ 1}{ 1}{ 1}{ 1}{ 1} \quad &\E{1}{ 1}{ 1}{ 1}{ 0}{ 1}{ 1}{ 1}{ 1}&E_8(a_1) &1&E_8(a_1)\\

20&\E{-9}{ 1}{ 1}{ 1}{ 1}{ 1}{ 1}{ 1}{1}\quad   &\E{1}{ 1}{ 1}{ 1}{ 0}{ 1}{ 0}{ 1}{ 1}&E_8(a_2) &1&E_8(a_2)\\

15&\E{-14}{ 1}{ 1}{ 1}{ 1}{ 1}{ 1}{ 1}{ 1}\quad&\E{1}{ 1}{ 0}{ 0}{ 1}{ 0}{ 1}{ 0}{ 1}&E_8(a_5)
&1&E_8(a_4)\\

12&\E{-17}{ 1}{ 1}{ 1}{ 1}{ 1}{ 1}{ 1}{ 1}\quad&\E{1}{ 1}{ 0}{ 0}{ 1}{ 0}{ 0}{ 1}{ 0}&E_8(a_3)
&1&E_8(a_5)\\

10&\E{-19}{ 1}{ 1}{ 1}{ 1}{ 1}{ 1}{ 1}{ 1}\quad&\E{1}{ 0}{ 0}{ 0}{ 1}{ 0}{ 0}{ 1}{ 0}&E_8(a_6)=-2A_4
&1&E_8(a_6)\\

8&\E{-21}{ 1}{ 1}{ 1}{ 1}{ 1}{ 1}{ 1}{ 1}\quad &\E{0}{ 0}{ 0}{ 0}{ 1}{ 0}{ 0}{ 0}{ 1}&D_8(a_3) &\bfmu_2\times\bfmu_2&-\\

6&\E{-23}{ 1}{ 1}{ 1}{ 1}{ 1}{ 1}{ 1}{ 1}\quad &\E{1}{ 0}{ 0}{ 0}{ 0}{ 1}{ 0}{ 0}{ 0}&E_8(a_8)=-4A_2
&1&E_8(a_7)\\

5&\E{-24}{ 1}{ 1}{ 1}{ 1}{ 1}{ 1}{ 1}{ 1}\quad &\E{0}{ 0}{ 0}{ 0}{ 0}{ 1}{ 0}{ 0}{ 0}&2A_4
&\bfmu_5\times\bfmu_5&-\\

4&\E{-25}{ 1}{ 1}{ 1}{ 1}{ 1}{ 1}{ 1}{ 1}\quad&\E{0}{ 0}{ 0}{ 0}{ 0}{ 0}{ 1}{ 0}{ 0}&2D_4(a_1)
&\bfmu_2^4&-\\

3&\E{-26}{ 1}{ 1}{ 1}{ 1}{ 1}{ 1}{ 1}{ 1}\quad&\E{0}{ 0}{ 1}{ 0}{ 0}{ 0}{ 0}{ 0}{ 0}&4A_2
&\bfmu_3^4&-\\

2&\E{-27}{ 1}{ 1}{ 1}{ 1}{ 1}{ 1}{ 1}{ 1}\quad  &\E{0}{ 1}{ 0}{ 0}{ 0}{ 0}{ 0}{ 0}{0}&8A_1=-1
&\bfmu_2^8&-\\
\hline
\end{array}}
$$
\end{center}

\vskip15pt
\begin{center}
{\small Table 6: The stable gradings  for ${F_4}$}
$$
{\renewcommand{\arraystretch}{1.3}
\begin{array}{cccccc}
\hline
m&\text{un-normalized}&\text{normalized}& w&S_0&A \\
\hline
12=h_\vt&\ \ \FIV{1}{1}{1}{1}{1} &\FIV{1}{1}{1}{1}{1}&F_4 &1&F_4\\
8&\FIV{-3}{1}{1}{1}{1} &\FIV{1}{1}{1}{0}{1}&B_4 
&\bfmu_2&F_4(a_1)\\
6&\FIV{-5}{1}{1}{1}{1}&\FIV{1}{0}{1}{0}{1}&F_4(a_1)
&1&F_4(a_2)\\
4&\FIV{-7}{1}{1}{1}{1}&\FIV{1}{0}{1}{0}{0}&D_4(a_1) 
&\bfmu_2\times\bfmu_2&F_4(a_3)\\
3&\FIV{-8}{1}{1}{1}{1}&\FIV{0}{0}{1}{0}{0}&A_2+\tilde A_2
&\bfmu_3\times\bfmu_3&-\\
2&\FIV{-9}{1}{1}{1}{1} &\FIV{0}{1}{0}{0}{0}&4A_1 
&\bfmu_2^4&-\\
\hline
\end{array}}
$$
\end{center}

\vskip15pt
\begin{center}
{\small Table 7: The stable gradings  for ${G_2}$ }
$$
{\renewcommand{\arraystretch}{1.3}
\begin{array}{cccccc}
\hline
m&\text{un-normalized}&\text{normalized}& w &S_0&A\\
\hline
6=h_\vt&\ \ 1\ 1\Rrightarrow 1 &1\ 1\Rrightarrow 1&G_2 
&1&G_2\\
 
3&-2\ 1\Rrightarrow 1 &1\ 1\Rrightarrow 0&A_2
&\bfmu_3&G_2(a_1)\\
 
2&-3\ 1\Rrightarrow 1&0\ 1\Rrightarrow 0&A_1+\tilde A_1
&\bfmu_2\times\bfmu_2&-\\
\hline
\end{array}}
$$
\end{center}

\vskip25pt
\begin{center}
{\small Table 8: The stable gradings  for ${^3D_4}$ }
$$
{\renewcommand{\arraystretch}{1.3}
\begin{array}{c c c c c }
\hline
m&\text{un-normalized}&\text{normalized}& w\vt\in W(F_4)&S_0\\
\hline
12=h_\vt&1\ 1\Lleftarrow 1 &1\ 1\Lleftarrow 1&F_4 
&1\\

6&-1\ 1\Lleftarrow 1&1\ 0\Lleftarrow 1&F_4(a_1)
&1\\

3&-2\ 1\Lleftarrow 1&0\ 0\Lleftarrow 1&A_2+\tilde A_2
&\bfmu_3\times\bfmu_3\\
\hline
\end{array}}
$$
\end{center}

\subsection{Stable gradings of classical Lie algebras}
\label{classical}
Here we tabulate the stable gradings of classical Lie algebras. 
For inner type $A_n$ the only stable grading is the Coxeter one, so we omit this case.

\subsubsection{Type ${^2\!A_\ell}$}
The stable gradings in type ${^2\!A_\ell}$ correspond to divisors of $\ell$ and $\ell+1$, each having odd quotient $d=m/2$. Conjugacy classes in the symmetric group are denoted by their partitions. For example, $[d^{2k+1}]$ consists of the products of 
$2k+1$ disjoint $d$-cycles. 

\begin{center}
{\small Table 9: The stable gradings for ${^2A_2}$}
$$
{\renewcommand{\arraystretch}{1.3}
\begin{array}{ c c c c}
\hline
m=2d&\text{Kac diagram}& w\vt&S_0\\
\hline
6=h_\vt
&
\twoAtwooo
&-1\times [3]
&1\\
2 &\twoAtwoo &-[1^3]
&\bfmu_2\times\bfmu_2\\
\hline
 \end{array}}
 $$
 \end{center}

\begin{center}
{\small Table 10: The stable gradings for ${^2A_{2n}}$, $n\geq 2$}
$$
{\renewcommand{\arraystretch}{1.5}
\begin{array}{c c c c}
\hline
m=2d&\text{Kac diagram}& w\vt&S_0\\
\hline
2(2n+1)=h_\vt
&1 \Rightarrow 1\ \ 1\ \cdots\ 1\ \ 1\ \Rightarrow 1 &-1\times[2n+1]
&1\\
2&1 \Rightarrow 0\ \ 0\ \ 0\cdots\ 0\ \ 0\ \Rightarrow 0 &-1\times[1^{2n+1}]
&\bfmu_2^{2n}\\
\frac{2(2n+1)}{2{k}+1},\quad {k}>0 &
1 \Rightarrow\underset{A_{2{k}} } { \underbrace{ 0\cdots0}  }\ \ 1\ \ 
\underset{A_{2{k}} } { \underbrace{ 0\cdots0}  }\ \ 1\ \ \cdots\ \ 1\ \ 
\underset{A_{2{k}} } { \underbrace{ 0\cdots0}  }\Rightarrow 1 
&-1\times [d^{2{k}+1}]
&\bfmu_2^{2{k}}\\
\frac{2n}{k},\quad 1<\frac{n}{k}\ \text{odd}\quad &
1 \Rightarrow\underset{A_{2{k}-1} } { \underbrace{ 0\cdots0}  }\ \ 1\ \ 
\underset{A_{2{k}-1} } { \underbrace{ 0\cdots0}  }\ \ 1\ \ \cdots\ \ 1\ \ 
\underset{B_{k} } { \underbrace{ 0\cdots0 \Rightarrow 0}}
&-1\times [d^{2{k}},1]
&\bfmu_2^{2{k}}\\
\hline
 \end{array}}
 $$
 \end{center}

\begin{center}
{\small Table 11: The stable gradings for ${^2A_{2n-1}}$, $n\geq 3$}
$$
{\renewcommand{\arraystretch}{1.5}
\begin{array}{ cc c c}
\hline
m&\text{Kac diagram}&  w\vt &S_0\\
\hline
&&&\\
2(2n-1)=h_\vt&
\begin{split}
&1\\
1\ \  \ & 1\ \ \ 1\ \ \ 1\ \ \ 1\cdots 1\ \ \ 1\Leftarrow1\\
\end{split}\qquad 
& 
-1\times [2n-1]
&1\\
&&&\\
2n\quad (\text{$n$ odd})&
\begin{split}
&1\\
1\ \  \ & 0\ \ \ 1\ \ \ 0\ \ \ 1\cdots 1\ \ \ 0\Leftarrow1\\
\end{split}\qquad 
& 
-1\times [n^2]
&1\\
&&&\\
\frac{2(2n-1)}{2{k}+1},\quad {k}>0 &
\begin{split}
&\quad \ \ 0\\
&\underset{ D_{{k}+1}     } {\underbrace{0\ \ \ \ 0\ \cdots \ 0}}\ \ 1\ \ 
\underset{A_{2{k}} } {\underbrace{ 0\ \cdots \ 0}}\ \ 1\cdots
\ 1\ \ \underset{A_{2{k}} } {\underbrace{ 0\ \cdots\ 0}}
\Leftarrow 1
\end{split} 
&-1\times [d^{2{k}+1},1]&\bfmu_2^{2{k}}\\
&&&\\
\frac{2n}{k},\quad 1<\frac{n}{k}\  \text{odd}\quad
&
\begin{split}
&\quad \ \ 0\\
&\underset{ D_{{k}}     } {\underbrace{0\ \ \ \ 0\ \cdots \ 0}}\ \ 1\ \ 
\underset{A_{2{k}-1} } {\underbrace{ 0\ \cdots \ 0}}\ \ 1\cdots
\ 1\ \ \underset{A_{2{k}-1} } {\underbrace{ 0\ \cdots\ 0}}
\Leftarrow 1
\end{split} 
&-1\times [d^{2{k}}]
&\bfmu_2^{2{k}-2}\\
\hline
\end{array}}
$$
\end{center}

\newpage
\subsubsection{Types $B_n, C_n$}
The stable gradings  for type $B_n$ and $C_n$ correspond to divisors $k$ of $n$, with period $m=2n/k$. The corresponding class in $W(B_n)=W(C_n)$,
denoted $kB_{n/k}$, consists of the $k^{th}$ powers of a Coxeter element. 
\begin{center}
{\small Table 12: The stable gradings for type $B_n$}
$$
{\renewcommand{\arraystretch}{1.3}
\begin{array}{cccc}
\hline
k=\frac{2n}{m}& \text{Kac diagram}& w& S_0\\
\hline 
&&&\\
1
&
\begin{split}
&1\\
1\ \ \ & 1\ \ \ 1\ \ \ 1\ \ \ 1\cdots 1\ \ \ 1\Rightarrow1\\
\end{split}
& 
B_{n}
&1\\
 &&&\\
\underset{n\ \text{even}}{2}
&
\begin{split}
&1\\
1\ \ \ & 0\ \ \ 1\ \ \ 0\ \ \ 1\cdots 0\ \ \ 1\Rightarrow0\\
\end{split}
& 
2 B_{n/2}
&\bfmu_2\\
&&&\\
\underset{{k\ \text{even}}}{k>2}
&\begin{split}
&\quad \ 0\\
&\underset{ D_{{k}/2}     } {\underbrace{0\ \ \ 0\ \cdots \ 0}}\ \ \ 1\ \ \ 
\underset{A_{{k}-1} } {\underbrace{ 0\ \cdots \ 0}}\ \ 1\cdots
1\ \ \ \underset{A_{{k}-1} } {\underbrace{ 0\ \cdots\   0}}\ \ \ 1\ \ \ 
\underset{B_{{k}/2} } {\underbrace{ 0\cdots 0\  \Rightarrow 0}}
\end{split} 
& {k}  B_{n/{k}}
&\bfmu_2^{{k}-1}\\
&&&\\
\underset{k\ \text{odd}}{k>1}
&
  \begin{split}
&\quad \ 0\\
&\underset{ D_{({k}+1)/2}     } {\underbrace{0\ \ \ 0\ \cdots  \ 0}}\ \ \ 1\ \ \ 
\underset{A_{{k}-1} } {\underbrace{ 0\ \cdots \ \ 0}}\ \ 1\cdots
1\ \ \ \underset{A_{{k}-1} } {\underbrace{ 0\ \cdots\ \  0}}\ \ \ 1\ \ \ 
\underset{B_{({k}-1)/2}  } {\underbrace{ 0\cdots 0\  \Rightarrow 0}}
\end{split} 
&
{k}  B_{n/{k}}
&\bfmu_2^{{k}-1}\\
\hline
\end{array}}
$$
\end{center}

\vskip50pt
\begin{center}
{\small Table 13: The stable gradings for type $C_n$}
$$
{\renewcommand{\arraystretch}{1.5}
\begin{array}{cccc}
\hline
k=\frac{2n}{m}&\text{Kac diagram}&  w&S_0\\
\hline 
1
&
1 \Rightarrow 1\ \ 1 \cdots 1\ \ 1\Leftarrow 1 
& B_{n}
&1 \\
 
k>1\quad
&
1 \Rightarrow\underset{A_{{k}-1} } { \underbrace{ 0\cdots0}  }\ \ 1\ \ 
\underset{A_{{k}-1} } { \underbrace{ 0\cdots0}  }\ \ 1\ \ \cdots\ \ 1\ \ 
\underset{A_{{k}-1} } { \underbrace{ 0\cdots0}  }\Leftarrow 1 
&{k}B_{n/{k}}
&\bfmu_2^{{k}-1} \\
\hline
 \end{array}}
 $$
 \end{center}

\newpage
\subsubsection{Types $D_n$ and ${^2\!D_n}$ \quad ($n\geq 4$)}

The stable gradings  for type $D_n$ correspond to even divisors $k$ of $n$ and odd divisors $\ell$ of $n-1$. The stable gradings  for type ${^2\!D_n}$ correspond to odd divisors $\ell$ of $n$ and even divisors $k$ of $n-1$.
\begin{center}
{\small Table 14: The stable gradings for type $D_n$, $n\geq 4$}
$$
{\renewcommand{\arraystretch}{1.5}
\begin{array}{cccc}
\hline
m&\text{Kac diagram}&  w&S_0\\
\hline
2n-2=h_\vt\quad&
\begin{split}
&\ 1\qquad\qquad\  1\\
1\ \ \ \!& \ 1\ \ \ 1\cdots\  1\ \ \  1\ \ \ 1\\
\end{split}
& B_1+B_{n-1}
&1\\
&&&\\
n\quad (\text{if $n$ is even})&
\begin{split}
&1\qquad\qquad\qquad\ \  1\\
1\ \ & 0\ \ 1\ \ 0\ \ 1\cdots0\ \ 1\ \  0\ \ 1\\
\end{split}
& 2 B_{n/2}
&1\\
&&&\\
\frac{2n}{k}\quad  2<{k}\  \text{even}&
 \begin{split}
&\quad 0\qquad\qquad\qquad\qquad\qquad\qquad\qquad\qquad
\ \ \!\ \ \!0\\
&\underset{ D_{{k}/2}     } {\underbrace{0\ \ 0\cdots \  0}}\ \ 1\ \ 
\underset{A_{{k}-1} } {\underbrace{ 0\  \cdots \  0}}\ \ 1\cdots
1\ \ \underset{A_{{k}-1} } {\underbrace{ 0\ \cdots\  0}}\ \ 1\ \ 
\underset{D_{{k}/2} } {\underbrace{ 0\cdots 0\  \  \ 0}}
\end{split} 
& {k} B_{n/{k}}
&\bfmu_2^{{k}-2}\\
&&&\\
\frac{2n-2}{\ell}\quad 1< {\ell}\ \text{odd}&
 \begin{split}
&\quad 0\qquad\qquad\qquad\qquad\qquad\qquad\qquad\qquad\ \ \ \!0\\
&\underset{ D_{({\ell}+1)/2}     } {\underbrace{0\ \ 0\cdots \ 0}}\ \ 1\ \ 
\underset{A_{{\ell}-1} } {\underbrace{ 0\  \cdots \ 0}}\ \ 1\cdots
1\ \ \underset{A_{{\ell}-1} } {\underbrace{ 0\  \cdots\   0}}\ \ 1\ \ 
\underset{D_{({\ell}+1)/2} } {\underbrace{ 0\cdots 0\  \  \ 0}}
\end{split} 
&B_1+{\ell} B_{(n-1)/{\ell}}
&\bfmu_2^{{\ell}-1}\\
\hline
\end{array}}
$$
\end{center}

\vskip25pt
\begin{center}
{\small Table 15: The stable gradings for type ${^2\!D_n}$, $n\geq 3$}
$$
{\renewcommand{\arraystretch}{1.5}
\begin{array}{cccc}
\hline
m&\text{Kac diagram}&  w&S_0\\
\hline
2n=h_\vt
&
1\Leftarrow 1\ \ 1\cdots1\ \ \ 1 \Rightarrow 1
& B_n
&1\\
&&&\\
n-1\quad 
(\text{if $n$ is odd})&
0\Leftarrow 1\ \ 0\ \ 1\ \ 0\cdots1\ \ 0\ \  1 \Rightarrow 0
& B_1+2 B_{n/2}
&\mu_2\times\mu_2\\
&&&\\
\frac{2n}{\ell}\quad  2<{\ell}\  \text{odd}
&\underset{ B_{({\ell}-1)/2}     } {\underbrace{0\Leftarrow 0\cdots  0}}
\ \ 1\ \ 
\underset{A_{{\ell}-1} } {\underbrace{ 0\cdots 0}}\ \ 1\cdots
\ 1\ \underset{A_{{\ell}-1} } {\underbrace{ 0\cdots 0}}\ \ 1\ \ 
\underset{B_{({\ell}-1)/2} } {\underbrace{ 0\cdots 0 \Rightarrow 0}}
& {\ell} B_{n/{\ell}}
&\bfmu_2^{{\ell}-1}\\
&&&\\
\frac{2n-2}{k}\quad 1< {k}\ \text{even}&
\underset{ B_{{k}/2}     } {\underbrace{0\Leftarrow 0\cdots  0}}
\ \ 1\ \ 
\underset{A_{{k}-1} } {\underbrace{ 0\cdots  0}}\ \ 1\cdots
\ 1\  \underset{A_{{k}-1} } {\underbrace{ 0 \cdots 0}}\ \ 1\ \ 
\underset{B_{k/2} } {\underbrace{ 0\cdots 0 \Rightarrow 0}}
&B_1+{k} B_{(n-1)/{k}}
&\bfmu_2^{{k}}\\
\hline
\end{array}}
$$
\end{center}

\subsection{Distinguished nilpotent elements and stable gradings}
\label{distinguished}

Kac coordinates of stable gradings are of two kinds, according as $s_0=0$ or $s_0=1$. 
Expanding on section 9 of \cite{springer:regular}, we show here that all stable gradings with $s_0=1$ in exceptional Lie algebras are related to distinguished nilpotent elements. 
For simplicity, we assume in this section only that $k$ has characteristic zero. 

Let $A$ be a distinguished nilpotent element in $\fg$. 
That is, the connected centralizer $C_G(A)^\circ$ is unipotent. 
There is a  homomorphism $\check\lam:k^\times\to G$, 
such that $\Ad(\check\lam(t))A=tA$ for all 
$t\in k^\ast$.
This gives a grading 
$$\fg=\bigoplus_{j=-a}^{a}\fg(j),$$
where $\fg(j)=\{x\in \fg:\ \lam(t)x=t^j\cdot x\ \  \forall t\in k^\times\}$ and 
$a=\max\{j:\ \fg(j)\neq 0\}.$ 
Since  $A$ is distinguished the linear map $\ad(A):\fg(0)\to\fg(1)$ is a bijection. 

Set $m=a+1$, assume this is nonzero in $k$,  and choose a root of unity $\zeta\in k^\times$ of order $m$. 
The inner automorphism  $\theta_A:=\Ad(\check\lam(\zeta))\in \Aut(\fg)^\circ$ has order $m$, 
giving rise to a $\bz/m$-grading 
$$\fg=\bigoplus_{i\in\bz/m}\ \fg_i,$$
where $\fg_i$ is the $\zeta^i$-eigenspace of $\theta_A$ in $\fg$. We have
$$\fg_i=\sum_{\substack{-a\leq j\leq a\\ j\equiv i\mod m}}\fg(j),$$
so that 
$$\fg_0=\fg(0)\qquad\text{and}\qquad \fg_1=\fg(-a)\oplus\fg(1).$$ 

Choose a maximal torus $T$ in a Borel subgroup $B$ of $G$ such that $\check\lam\in X_\ast(T)$ and $\la \al,\check\lam\ra\geq 0$ for all roots $\al$ of $T$ in $B$. For each of the simple roots 
$\al_1,\dots,\al_\ell$ we have $\la \al_i,\check\lam\ra\in\{0,1\}$. We set 
$s_i=\la \al_i,\check\lam\ra$, and also put $s_0=1$. 
Since $\fg(-a)$ contains the lowest root space, it follows that $(s_0,s_1,\dots, s_\ell)$ are the normalized Kac-coordinates of $\theta_A$. 

\begin{prop}\label{9.5} The following are equivalent.
\begin{enumerate}
\item There exists $M\in \fg(-a)$ such that $M+A$ is regular semisimple.
\item There exists $M\in \fg(-a)$ such that $M+A$ is semisimple.
\item The automorphism $\theta_A$ is stable. 
\end{enumerate}
\end{prop}
\proof Implication $1\Rightarrow 2$ is obvious. 

We prove $2\Rightarrow 3$.  
Since $A$ is distinguished, the centralizer $C_{G_0}(A)$ is finite.
Since $G_0$ preserves each summand $\fg(j)$, we have 
$C_{G_0}(M+A)\subset C_{G_0}(A)$. Hence  $C_{G_0}(M+A)$ is also finite, so the $G_0$-orbit of $M+A$ in $\fg_1$  is stable. 

The implication $3\Rightarrow 1$ is proved in \cite[9.5]{springer:regular}. 
We give Springer's argument here for completeness. 
Let $F$ be a $G$-invariant polynomial on $\fg$ such that $F(x)\neq 0$ if and only if $x$ is regular semisimple. For example, we can choose $F$ corresponding, under the Chevalley isomorphism $k[\ft]^G\overset\sim\to k[\ft]^W$,  to the product of the roots.  Now assuming that $3$ holds, there are  vectors $Z\in \fg(-a)$ and $Y_0\in\fg(1)$ such that $Z+Y_0$ is semisimple and has finite stabilizer in $G_0$. The centralizer $\fm=\fz(Z+Y_0)$ is then reductive, with $\fm^\theta=0$, so $\fm$ is a Cartan subalgebra of $\fg$ and $Z+Y_0$ is in fact regular semisimple. Hence the polynomial $F_Z$ on $\fg(1)$ given by $F_Z(Y):=F(Z+Y)$ does not vanish identically.
Since $A$ is distinguished, the orbit $\Ad(G_0)A$ is dense in $\fg(1)$, so there is $g\in G_0$ such that $F_Z(\Ad(g)A)=F(\Ad(g)^{-1}Z+A)\neq 0$. It follows that $\Ad(g)^{-1}Z+A$ is regular semisimple so $1$ holds. \qed

We say that a distinguished nilpotent element $A\in\fg$ is {\bf $S$-distinguished} if the equivalent conditions of Prop. \ref{9.5} hold. 

{\bf A non-example:\ } It can happen that $\fg(-a)+\fg(1)$ contains semisimple elements, but none have the form $M+A$ with $M\in\fg(-a)$. For example, suppose $\fg=\fsp_6$ and $A$ has Jordan blocks $(4,2)$. The automorphism $\theta_A$ has Kac coordinates 
$$1\Rightarrow 1\ \  \ \!0\Leftarrow 1$$
and has rank equal to $1$. 
It corresponds to $w\in W(C_3)$ of type 
$C_2\times C_1$, which is not $\bz$-regular, so $A$ is not $S$-distinguished.

\begin{prop}\label{s0} Assume that $\fg$ is of exceptional type and that 
$\theta\in\Aut(\fg)^\circ$ is a stable inner automorphism whose Kac coordinates satisfy $s_0=1$. Then $\theta=\theta_A$ where $A$ is an $S$-distinguished nilpotent element in $\fg$. 
\end{prop}
\proof In the tables of section \ref{exceptional} we have listed, for each $\theta$ with $s_0=1$, the conjugacy class of a nilpotent element $A$ such that $\theta_A$ has the normalized Kac coordinates of $\theta$. 
\qed

{\bf Remark 1:\ } For $n$ even there is a unique $S$-distinguished non-regular nilpotent class
in $\fso_{2n}$ which is also $S$-distinguished in $\fso_{2n+1}$, having Jordan partitions $[2n+1,2n-1]$ and $[2n+1, 2n-1,1]$, respectively. For $A$ in these classes $\theta_A$ has order $n$. In these and the exceptional cases, the map $A\mapsto\theta_A$ is a  bijection from the set of $S$-distinguished nilpotent $G$-orbits in $\fg$ to the set of inner gradings on $\fg$ with $s_0=1$. 
However, Prop. \ref{9.5} is false for  $C_n$, $n\geq 2$.

{\bf Remark 2:\ } If $A$ is $S$-distinguished then $\fz(M+A)$ is a canonical Cartan subalgebra for $\theta_A$ on which $\theta_A$ acts by an element of the conjugacy class in $W$ associated to $A$ via the Kazhdan-Lusztig map \cite{kazhdan-lusztig:affinefixedpoints}. This follows from the argument in \cite[9.11]{kazhdan-lusztig:affinefixedpoints}, 
and confirms two entries in \cite[Table 1]{spaltenstein:kl} (for $A=E_8(a_6), E_8(a_7)$), listed there as conjectural.

{\bf Remark 3:\ } There are exactly three cases where  $\fg_0$ is a maximal proper Levi subalgebra in $\fg$. These occur in $G_2, F_4$ and $E_8$, for $a=2,3,5$ respectively, where $C_{G_0}(A)$ is a symmetric group $S_3, S_4, S_5$. These groups act irreducibly on the subspaces $\fg(-a)$ of  dimensions $1,2,4$, in which the stabilizers of a vector in general position are the isotropy groups $S_0=\bfmu_3, \bfmu_2\times\bfmu_2$, $1$. These are the maximal abelian normal subgroups of $C_{G_0}(A)$.

\section{Positive rank gradings for type $E_{6,7,8}$ (inner case)}\label{E678}

Assume now that $\fg$ has type $E_n$, for $n=6,7,8$. 
From Prop. \ref{Kac(w)} we have the following algorithm to find all inner automorphisms of 
$\fg$  having positive rank. For each $m\geq 1$ list 
the $W$-conjugacy classes of $m$-admissible elements in $W$. For a representative $w$ of each class, form the list $\Kac(w)_{\un}$ and apply the normalization algorithm to each element of 
$\Kac(w)_{\un}$, discarding duplicate results,  to obtain the list $\Kac(w)$ of normalized Kac coordinates. Then by Prop. \ref{Kac(w)}, 
the union of the lists $\Kac(w)$ over all conjugacy-classes of $m$-admissible $w$ gives all positive rank inner automorphisms of order $m$. 

To find the $\Kac(w)_{\un}$ when each $w_i$ is $\bz$-regular, we can use Prop. \ref{clift} to find the Kac coordinates of each $w_i$, which lead to those of $w$ via the normalization algorithm. 
It turns out that we obtain all positive rank gradings from those $m$-admissible $w$ for which each factor $w_i$ is not only elliptic but also $\bz$-regular in $W_{J_i}$. However, we do not have an {\it a priori} proof of this fact, 
so we must also compute Kac coordinates of lifts in the small number of cases 
where not all $w_i$ are $\bz$-regular. 

These non-regular cases are handled as follows. 
By induction, we assume $w=w_i$ lies in no proper reflection subgroup and we consider the powers of $w$. 
To illustrate the method, take the nonregular element $w=E_8(a_7)=-A_2E_6$ of order $12$. 
First list the $32$ normalized Kac coordinates $(s_i)$ with $s_i\in\{0,1\}$ 
and $s_0 + 2 s_1 + 3 s_2 + 4 s_3 + 6 s_4 + 5 s_5 + 4 s_6 + 3 s_7 + 2 s_8 =12$. 
We have $w^2$ and $w^3$ in the  classes $A_2E_6(a_2)$ and  $2A_2+2A_1$ 
whose lifts have Kac coordinates 
$\E{0}{0}{1}{0}{0}{0}{0}{1}{0}$ and $\E{0}{ 0}{ 0}{ 1}{ 0}{ 0}{ 0}{ 0}{ 0}$, respectively. 
Only one of the 32 elements on the list satisfies these two conditions, 
namely $\E{1}{ 0}{ 0}{ 1}{ 0}{ 1}{ 0}{ 0}{ 1}$. Therefore this is the Kac diagram 
for the lift of $w$ in the class
$E_8(a_7)$. 

\subsection{A preliminary list of Kac coordinates for positive rank gradings of inner type}
\label{preliminary}
For each possible order $m$ we  list the $m$-admissible elements in $W(E_{6,7,8})$,
 the rank $r=\rank(w)$. 
and the form of the un-normalized Kac-coordinates of the lifts of $w$ In the column $\Kac(w)_{\un}$, each $\ast$ is an independent variable integer ranging  over a set of representatives of $\bz/m$ such that the order is always $m$. For each vector of $\ast$-values we apply the normalization algorithm to obtain the normalized Kac coordinates $\Kac(w)$ in the last column. The sets $\Kac(w)$ are not disjoint. In  a second set of tables (section \ref{Eposrank}), we will select, for each $\theta$ appearing in $\cup_w\Kac(w)$, a $w$ of maximal rank for which $\Kac(w)$ contains $\theta$. 

We use Carter's notation for conjugacy classes in $W$, augmented as follows. 
If $X$ is a conjugacy class and $-1\in W$ then $-X=\{-w:\ w\in X\}$. This makes some classes easier to understand; for example, 
$E_8(a_7)=-A_2E_6$. 

\begin{center}
{\small Table 16:  $\Kac(w)_{\un}$ and $\Kac(w)$ for $m$-admissible $w$ in $W(E_6)$ }
$$
\begin{array}{|c|c|c|l|l|}
\hline
m& w & r&\Kac(w)_{\un}& \Kac(w)\\
\hline\hline
12& E_6 & 1& \EVI{1}{1}{1}{1}{1}{1}{1}& \EVI{1}{1}{1}{1}{1}{1}{1}\\
\hline
9& E_6(a_1) & 1& \EVI{1}{1}{1}{1}{0}{1}{1}& \EVI{1}{1}{1}{1}{0}{1}{1}\\
\hline
8& D_5& 1& \EVI{1}{ \ast}{ 1}{ 1}{ 1}{ 1}{ \ast}& 
\EVI{0}{ 1}{ 1}{ 1}{ 0}{ 1}{ 1}\qquad \EVI{1}{ 1}{ 1}{ 0}{ 1}{ 0}{ 1}\\
\hline
6& E_6(a_2) & 2& \EVI{1}{ 1}{ 0}{ 0}{ 1}{ 0}{ 1}& \EVI{1}{ 1}{ 0}{ 0}{ 1}{ 0}{ 1}\\
\hline
6& A_5& 1&\EVI{\ast}{ 1}{ \ast}{ 1}{ 1}{ 1}{ 1} &\EVI{1}{ 1}{ 0}{ 0}{ 1}{ 0}{ 1}
\qquad\EVI{0}{0}{1}{1}{0}{1}{0}\\
\hline
6& D_4& 1& \EVI{\ast}{ \ast}{ 1}{ 1}{ 1}{ 1}{ \ast}&\EVI{1}{ 0}{ 1}{ 0}{ 1}{ 0}{ 0} \qquad
\EVI{0}{ 1}{ 1}{ 0}{ 0}{ 1}{ 1}
       \qquad \EVI{0}{ 1}{ 1}{ 1}{ 0}{ 0}{ 1}\qquad \EVI{1}{ 1}{ 0}{ 0}{ 1}{ 0}{ 1}\\
\hline
5& A_4& 1&\EVI{\ast}{ 1}{ \ast}{ 1}{ 1}{ 1}{ \ast} &\EVI{1}{ 0}{ 0}{ 1}{ 0}{ 1}{ 0}
\qquad \EVI{0}{ 1}{ 0}{ 0}{ 1}{ 0}{ 1}
        \qquad \EVI{1}{ 1}{ 1}{ 0}{ 0}{ 0}{ 1}\\
\hline
4& D_4(a_1)& 2&\EVI{\ast}{ \ast}{ 1}{ 1}{ 0}{ 1}{ \ast} &
\EVI{1}{ 0}{ 0}{ 0}{ 1}{ 0}{ 0}\qquad \EVI{0}{ 1}{ 1}{ 0}{ 0}{ 0}{ 1}\\
\hline
4& A_3& 1&\EVI{\ast}{ \ast}{ \ast}{ 1}{ 1}{ 1}{ \ast} &
\EVI{0}{ 0}{ 0}{ 0}{ 1}{ 0}{ 1}\qquad \EVI{0}{ 0}{ 0}{ 1}{ 0}{ 1}{ 0}
\qquad \EVI{0}{ 1}{ 0}{ 0}{ 0}{ 1}{ 1}\qquad \EVI{0}{ 1}{ 0}{ 1}{ 0}{ 0}{ 1}
\qquad \EVI{1}{ 1}{ 0}{ 0}{ 0}{ 1}{ 0}\\
\hline
3& 3A_2& 3&\EVI{1}{ 1}{ 1}{ 1}{ \ast}{ 1}{ 1} &\EVI{0}{ 0}{ 0}{ 0}{ 1}{ 0}{ 0}\\
\hline
3& 2A_2& 2& \EVI{\ast}{ 1}{ \ast}{ 1}{ \ast}{ 1}{ 1} &\EVI{0}{ 0}{ 0}{ 0}{ 1}{ 0}{ 0}\qquad \EVI{1}{ 1}{ 0}{ 0}{ 0}{ 0}{ 1}\\
\hline
3& A_2& 1&\EVI{\ast}{\ast}{1}{\ast}{1}{\ast}{\ast} &
\EVI{1}{ 0}{ 1}{ 0}{ 0}{ 0}{ 0}\qquad \EVI{0}{ 0}{ 0}{ 0}{ 1}{ 0}{ 0}
\qquad \EVI{0}{ 0}{ 0}{ 1}{ 0}{ 0}{ 1}\qquad \EVI{0}{ 1}{ 0}{ 0}{ 0}{ 1}{ 0}
\qquad \EVI{1}{ 1}{ 0}{ 0}{ 0}{ 0}{ 1}\\
\hline
2& 4A_1& 4&\EVI{1}{ 1}{ \ast}{ \ast}{ 1}{ \ast}{ 1} &\EVI{0}{ 0}{ 1}{ 0}{ 0}{ 0}{ 0}\\
\hline
2& 3A_1& 3&\EVI{\ast}{ 1}{ \ast}{ \ast}{ 1}{ \ast}{ 1}  &\EVI{0}{ 0}{ 1}{ 0}{ 0}{ 0}{ 0}\\
\hline
2& 2A_1& 2&\EVI{1}{ \ast}{ \ast}{ \ast}{ 1}{ \ast}{ \ast} & \EVI{0}{ 0}{ 1}{ 0}{ 0}{ 0}{ 0}
\qquad \EVI{0}{ 1}{ 0}{ 0}{ 0}{ 0}{ 1}\\
\hline
2& A_1& 1&\EVI{\ast}{ \ast}{ \ast}{ \ast}{ 1}{ \ast}{ \ast}  &
\EVI{0}{ 0}{ 1}{ 0}{ 0}{ 0}{ 0}\qquad \EVI{0}{ 1}{ 0}{ 0}{ 0}{ 0}{ 1}\\
\hline
\end{array}
$$
\end{center}

\begin{center}
{\small Table 17: $\Kac(w)_{\un}$ and $\Kac(w)$ for $m$-admissible $w$ in $W(E_7)$ }
$$
\begin{array}{|c|c|c|l|l|}
\hline
m& w & r&\Kac(w)_{\un}& \Kac(w)\\
\hline\hline
18& E_7 & 1& \EVII{1}{1}{1}{1}{1}{1}{1}{1}& \EVII{1}{1}{1}{1}{1}{1}{1}{1}\\
\hline
14& E_7(a_1)& 1& \EVII{1}{1}{1}{1}{0}{1}{1}{1}& \EVII{1}{1}{1}{1}{0}{1}{1}{1}\\
\hline
12& E_7(a_2) & 1& \EVII{1}{ 1}{ 1}{ 0}{ 1}{ 0}{ 1}{ 1}& \EVII{1}{ 1}{ 1}{ 0}{ 1}{ 0}{ 1}{ 1}\\
\hline
12& E_6& 1& \EVII{\ast}{1}{1}{1}{1}{1}{1}{\ast}& 
\EVII{1}{0}{1}{1}{0}{1}{1}{1}\qquad \EVII{0}{ 1}{ 0}{ 0}{ 1}{ 1}{ 1}{ 1}
\qquad\EVII{1}{ 1}{ 1}{ 0}{ 1}{ 0}{ 1}{ 1}
\\
\hline
10& D_6 & 1& \EVII{\ast}{\ast}{1}{1}{1}{1}{1}{1}& 
\EVII{0}{ 1}{ 1}{ 0}{ 1}{ 0}{ 1}{ 0}\qquad\EVII {1}{ 0}{ 1}{ 1}{ 0}{ 1}{ 0}{ 1}
\qquad\EVII {1}{ 1}{ 0}{ 0}{ 1}{ 0}{ 1}{ 1}\\
\hline
9& E_6(a_1) & 1& \EVII{\ast}{1}{1}{1}{0}{1}{1}{\ast}& 
\EVII{0}{ 1}{ 0}{ 0}{ 1}{ 0}{ 1}{ 1}\qquad\EVII {1}{ 0}{ 1}{ 1}{ 0}{ 0}{ 1}{ 1} \\
\hline
8& D_5& 1& \EVII{\ast}{\ast}{1}{1}{1}{1}{1}{\ast}& 
\EVII{0}{ 0}{ 1}{ 1}{ 0}{ 0}{ 1}{ 1}\qquad\EVII{0}{ 1}{ 0}{ 0}{ 0}{ 1}{ 1}{ 1} 
\qquad\EVII{0}{ 1}{ 0}{ 0}{ 1}{ 0}{ 1}{ 0}
\qquad\EVII{0}{ 1}{ 1}{ 0}{ 0}{ 1}{ 0}{ 1}\\
&&&&\EVII{1}{ 0}{ 0}{ 0}{ 1}{ 0}{ 1}{ 1}\qquad\EVII{1}{ 0}{ 0}{ 1}{ 0}{ 1}{ 0}{ 1}
\qquad\EVII{1}{ 0}{ 1}{ 0}{ 1}{ 0}{ 0}{ 1}\qquad\EVII{1}{ 1}{ 1}{ 0}{ 0}{ 0}{ 1}{ 1}\\
\hline
8& D_6(a_1) & 1& \EVII{\ast}{\ast}{1}{1}{0}{1}{1}{1}& 
\EVII{0}{ 1}{ 1}{ 0}{ 0}{ 1}{ 0}{ 1}\qquad\EVII {1}{ 0}{ 0}{ 0}{ 1}{ 0}{ 1}{ 1}\\
\hline

8& A_7& 1& \EVII{1}{1}{\ast}{1}{1}{1}{1}{1}& 
\EVII{0}{ 1}{ 0}{ 0}{ 1}{ 0}{ 1}{ 0} \\
\hline
7& A_6 & 1& \EVII{\ast}{1}{\ast}{1}{1}{1}{1}{1}& 
\EVII{0}{ 1}{ 0}{ 0}{ 1}{ 0}{ 0}{ 1} \\
\hline
6& E_7(a_4) & 3& \EVII{1}{0}{0}{0}{1}{0}{0}{1}& \EVII{1}{0}{0}{0}{1}{0}{0}{1}\\
\hline
6& D_6(a_2) & 2& \EVII{\ast}{\ast}{1}{1}{0}{1}{0}{1}& 
\EVII{0}{ 1}{ 1}{ 0}{ 0}{ 0}{ 1}{ 0}\qquad\EVII{1}{ 0}{ 0}{ 0}{ 1}{ 0}{ 0}{ 1}\\
\hline
6& E_6(a_2) & 2& \EVII{\ast}{1}{0}{0}{1}{0}{1}{\ast}& 
\EVII{0}{ 1}{ 0}{ 0}{ 0}{ 1}{ 0}{ 1}\qquad\EVII {1}{ 0}{ 0}{ 0}{ 1}{ 0}{ 0}{ 1}\\
\hline
6& D_4& 1& \EVII{\ast}{\ast}{1}{1}{1}{1}{\ast}{\ast}& 
\EVII{1}{ 0}{ 1}{ 0}{ 0}{ 1}{ 0}{ 0}\qquad\EVII {0}{ 0}{ 0}{ 0}{ 0}{ 1}{ 1}{ 1}
\qquad\EVII{0}{ 0}{ 0}{ 1}{ 0}{ 0}{ 1}{ 1}\qquad\EVII{0}{ 0}{ 1}{ 0}{ 1}{ 0}{ 0}{ 0}
\qquad\EVII{0}{ 0}{ 1}{ 1}{ 0}{ 0}{ 0}{ 1}\\
&&&&\EVII {0}{ 1}{ 0}{ 0}{ 0}{ 1}{ 0}{ 1}
\qquad\EVII {0}{ 1}{ 1}{ 0}{ 0}{ 0}{ 1}{ 0}\qquad\EVII {1}{ 0}{ 0}{ 0}{ 1}{ 0}{ 0}{ 1}
\qquad\EVII {1}{ 0}{ 1}{ 0}{ 0}{ 0}{ 1}{ 1}\\
\hline

6& A_5'' & 1& \EVII{\ast}{1}{\ast}{1}{1}{1}{1}{\ast}& 
\EVII{0}{ 0}{ 0}{ 0}{ 1}{ 0}{ 1}{ 0}\qquad\EVII {0}{ 1}{ 0}{ 0}{ 0}{ 1}{ 0}{ 1}
\qquad\EVII{0}{ 1}{ 1}{ 0}{ 0}{ 0}{ 1}{ 0}\qquad\EVII{1}{ 0}{ 0}{ 0}{ 1}{ 0}{ 0}{ 1}\\
\hline
6& A_5' & 1& \EVII{1}{1}{1}{1}{1}{\ast}{\ast}{\ast}& 
\EVII{0}{ 0}{ 0}{ 1}{ 0}{ 1}{ 0}{ 0}\qquad\EVII {0}{ 1}{ 1}{ 0}{ 0}{ 0}{ 1}{ 0}
\qquad\EVII{1}{ 0}{ 0}{ 0}{ 1}{ 0}{ 0}{ 1}\qquad\EVII{1}{ 1}{ 0}{ 0}{ 0}{ 0}{ 1}{ 1}\\
\hline
5& A_4& 1& \EVII{\ast}{\ast}{\ast}{1}{1}{1}{1}{\ast}& 
\EVII{0}{ 0}{ 0}{ 0}{ 1}{ 0}{ 0}{ 1}\qquad\EVII {0}{ 0}{ 0}{ 1}{ 0}{ 0}{ 1}{ 0}
\qquad\EVII{0}{ 1}{ 0}{ 0}{ 0}{ 0}{ 1}{ 1}\qquad\EVII{0}{ 1}{ 1}{ 0}{ 0}{ 0}{ 0}{ 1}
\qquad\EVII{1}{ 0}{ 0}{ 0}{ 0}{ 1}{ 0}{ 1}\\
\hline

4& 2A_3& 2& \EVII{1}{1}{\ast}{1}{\ast}{1}{1}{1}& 
\EVII{0}{ 0}{ 0}{ 0}{ 1}{ 0}{ 0}{ 0} \qquad\EVII {0}{ 1}{ 0}{ 0}{ 0}{ 0}{ 1}{ 0} \\
\hline

4& D_4(a_1)& 2& \EVII{\ast}{\ast}{1}{1}{0}{1}{\ast}{\ast}& 
\EVII{0}{ 0}{ 0}{ 0}{ 0}{ 1}{ 0}{ 1}\qquad\EVII {0}{ 0}{ 0}{ 0}{ 1}{ 0}{ 0}{ 0}
\qquad\EVII{0}{ 0}{ 0}{ 1}{ 0}{ 0}{ 0}{ 1}\qquad\EVII{0}{ 1}{ 0}{ 0}{ 0}{ 0}{ 1}{ 0}
\qquad\EVII{1}{ 0}{ 1}{ 0}{ 0}{ 0}{ 0}{ 1}\\
\hline

4& A_3& 1& \EVII{\ast}{\ast}{\ast}{1}{1}{1}{\ast}{\ast}& 
\EVII{0}{ 0}{ 0}{ 0}{ 0}{ 1}{ 0}{ 1}\qquad\EVII {0}{ 0}{ 0}{ 0}{ 1}{ 0}{ 0}{ 0}
\qquad\EVII{0}{ 0}{ 0}{ 1}{ 0}{ 0}{ 0}{ 1}\qquad\EVII{0}{ 0}{ 1}{ 0}{ 0}{ 0}{ 1}{ 0}\\
&&&&\EVII{0}{ 1}{ 0}{ 0}{ 0}{ 0}{ 1}{ 0}\qquad\EVII{1}{ 0}{ 0}{ 0}{ 0}{ 0}{ 1}{ 1}
\qquad\EVII{1}{ 0}{ 1}{ 0}{ 0}{ 0}{ 0}{ 1}\\
\hline

3& 3A_2& 3& \EVII{1}{1}{1}{\ast}{1}{\ast}{1}{1}& \EVII{0}{ 0}{ 0}{ 0}{ 0}{ 1}{ 0}{ 0}\\
\hline

3& 2A_2& 2& \EVII{1}{1}{1}{\ast}{1}{\ast}{\ast}{\ast}& 
\EVII{0}{ 0}{ 0}{ 0}{ 0}{ 1}{ 0}{ 0}\qquad\EVII {0}{ 1}{ 0}{ 0}{ 0}{ 0}{ 0}{ 1}\\
\hline

3& A_2& 1& \EVII{\ast}{\ast}{1}{\ast}{1}{\ast}{\ast}{\ast}& 
\EVII{0}{ 0}{ 0}{ 0}{ 0}{ 0}{ 1}{ 1}\qquad\EVII{0}{ 0}{ 0}{ 0}{ 0}{ 1}{ 0}{ 0}
\qquad\EVII{0}{ 0}{ 1}{ 0}{ 0}{ 0}{ 0}{ 1}\qquad\EVII{0}{ 1}{ 0}{ 0}{ 0}{ 0}{ 0}{ 1}\\
\hline

2& 7A_1& 7& \EVII{0}{ 0}{ 1}{ 0}{ 0}{ 0}{ 0}{ 0}&\EVII{0}{ 0}{ 1}{ 0}{ 0}{ 0}{ 0}{ 0}\\
\hline
2& 6A_1& 6& \EVII{\ast}{\ast}{0}{0}{0}{1}{0}{0}&\EVII{0}{ 0}{ 1}{ 0}{ 0}{ 0}{ 0}{ 0}\\
\hline
2& 5A_1& 5& \EVII{1}{\ast}{1}{1}{\ast}{1}{\ast}{1}&\EVII{0}{ 0}{ 1}{ 0}{ 0}{ 0}{ 0}{ 0}\\
\hline
2& 4A_1'& 4& \EVII{1}{\ast}{1}{1}{\ast}{1}{\ast}{\ast}& 
\EVII{0}{ 0}{ 1}{ 0}{ 0}{ 0}{ 0}{ 0}\\
\hline
2& 4A_1''& 4& \EVII{\ast}{\ast}{0}{0}{1}{0}{\ast}{\ast}& 
\EVII{0}{ 0}{ 1}{ 0}{ 0}{ 0}{ 0}{ 0}\qquad\EVII {0}{ 0}{ 0}{ 0}{ 0}{ 0}{ 1}{ 0}\\
\hline
2& 3A_1'& 3& \EVII{1}{\ast}{1}{1}{\ast}{\ast}{\ast}{\ast}& 
\EVII{0}{ 0}{ 1}{ 0}{ 0}{ 0}{ 0}{ 0}\qquad\EVII {1}{ 0}{ 0}{ 0}{ 0}{ 0}{ 0}{ 1}\\
\hline
2& 3A_1''& 3& \EVII{1}{\ast}{1}{\ast}{\ast}{1}{\ast}{\ast}& 
\EVII {0}{ 0}{ 1}{ 0}{ 0}{ 0}{ 0}{ 0}\qquad\EVII{0}{ 0}{ 0}{ 0}{ 0}{ 0}{ 1}{ 0}\\
\hline
2& 2A_1& 2& \EVII{\ast}{\ast}{\ast}{1}{\ast}{1}{\ast}{\ast}& 
\EVII{0}{ 0}{ 1}{ 0}{ 0}{ 0}{ 0}{ 0}\qquad\EVII {0}{ 0}{ 0}{ 0}{ 0}{ 0}{ 1}{ 0}
\qquad\EVII {1}{ 0}{ 0}{ 0}{ 0}{ 0}{ 0}{ 1}\\
\hline
2& A_1& 1& \EVII{\ast}{\ast}{\ast}{\ast}{\ast}{1}{\ast}{\ast}& 
\EVII{0}{ 0}{ 1}{ 0}{ 0}{ 0}{ 0}{ 0}\qquad\EVII {0}{ 0}{ 0}{ 0}{ 0}{ 0}{ 1}{ 0}
\qquad\EVII {1}{ 0}{ 0}{ 0}{ 0}{ 0}{ 0}{ 1}\\
\hline

\end{array}
$$
\end{center}

\begin{center}
{\small Table 18: $\Kac(w)_{\un}$ and $\Kac(w)$ for $m$-admissible $w$ in $W(E_8)$ }
$$
\begin{array}{|c|c|c|l|l|}
\hline
m& w & r&\Kac(w)_{\un}& \Kac(w)\\
\hline\hline
30& E_8 & 1& \E{1}{ 1}{ 1}{ 1}{ 1}{ 1}{ 1}{ 1}{ 1}& \E{1}{ 1}{ 1}{ 1}{ 1}{ 1}{ 1}{ 1}{ 1}\\
\hline
24& E_8(a_1)& 1& \E{1}{ 1}{ 1}{ 1}{ 0}{ 1}{ 1}{ 1}{ 1}& \E{1}{ 1}{ 1}{ 1}{ 0}{ 1}{ 1}{ 1}{ 1}\\
\hline
20& E_8(a_2) & 1& \E{1}{ 1}{ 1}{ 1}{ 0}{ 1}{ 0}{ 1}{ 1}& \E{1}{ 1}{ 1}{ 1}{ 0}{ 1}{ 0}{ 1}{ 1}\\
\hline
18& E_8(a_4) & 1& \E{1}{ 0}{ 1}{ 1}{ 0}{ 1}{ 0}{ 1}{ 1}& \E{1}{ 0}{ 1}{ 1}{ 0}{ 1}{ 0}{ 1}{ 1}\\
\hline

18& E_7 & 1& \E{\ast}{ 1}{ 1}{ 1}{ 1}{ 1}{ 1}{ 1}{ \ast}& 
\E{0}{ 1}{ 0}{ 1}{ 1}{ 0}{ 1}{ 0}{ 1}\qquad \E{1}{ 0}{ 1}{ 1}{ 0}{ 1}{ 0}{ 1}{ 1}
\qquad \E{1}{ 1}{ 0}{ 0}{ 1}{ 0}{ 1}{ 1}{ 1}\qquad \E{1}{ 1}{ 1}{ 0}{ 1}{ 0}{ 1}{ 0}{ 1}\\
&&&&\E{1}{ 1}{ 1}{ 1}{ 0}{ 1}{ 0}{ 1}{ 0}\\
\hline

15& E_8(a_5) & 1& \E{1}{ 1}{ 0}{ 0}{ 1}{ 0}{ 1}{ 0}{ 1}& \E{1}{ 1}{ 0}{ 0}{ 1}{ 0}{ 1}{ 0}{ 1}\\
\hline

14& E_7(a_1) & 1& \E{\ast}{ 1}{ 1}{ 1}{ 0}{ 1}{ 1}{ 1}{ \ast}& 
\E{0}{ 1}{ 0}{ 0}{ 1}{ 0}{ 1}{ 0}{ 1}\qquad \E{1}{ 0}{ 1}{ 1}{ 0}{ 0}{ 1}{ 0}{ 1}
\qquad \E{1}{ 1}{ 0}{ 0}{ 1}{ 0}{ 0}{ 1}{ 1}\qquad \E{1}{ 1}{ 1}{ 0}{ 0}{ 1}{ 0}{ 1}{ 0}
\\
\hline

14& D_8 & 1& \E{1}{ \ast}{ 1}{ 1}{ 1}{ 1}{ 1}{ 1}{ 1}& \E{0}{ 1}{ 0}{ 0}{ 1}{ 0}{ 1}{ 0}{ 1}\\
\hline

12& E_8(a_3) & 2& \E{1}{ 1}{ 0}{ 0}{ 1}{ 0}{ 0}{ 1}{ 0}& \E{1}{ 1}{ 0}{ 0}{ 1}{ 0}{ 0}{ 1}{ 0}\\
\hline

12& E_8(a_7) & 1& \E{1}{ 0}{ 0}{ 1}{ 0}{ 1}{ 0}{ 0}{ 1}& \E{1}{ 0}{ 0}{ 1}{ 0}{ 1}{ 0}{ 0}{ 1}\\
\hline

12& E_7(a_2) & 1& \E{\ast}{ 1}{ 1}{ 0}{ 1}{ 0}{ 1}{ 1}{ \ast}& 
\E{0}{ 1}{ 0}{ 1}{ 0}{ 0}{ 1}{ 0}{ 1}\qquad \E{1}{ 0}{ 0}{ 1}{ 0}{ 1}{ 0}{ 0}{ 1}
\qquad \E{1}{ 1}{ 0}{ 0}{ 1}{ 0}{ 0}{ 1}{ 0}\qquad \E{1}{ 1}{ 1}{ 0}{ 0}{ 0}{ 1}{ 0}{ 1}
\\
\hline

12& D_8(a_1) & 1& \E{1}{ 1}{ 0}{ 0}{ 1}{ 0}{ 0}{ 1}{ 0}& \E{1}{ 1}{ 0}{ 0}{ 1}{ 0}{ 0}{ 1}{ 0}\\
\hline

12& D_7 & 1& \E{\ast}{ \ast}{ 1}{ 1}{ 1}{ 1}{ 1}{ 1}{ 1}& \E{0}{ 0}{ 0}{ 0}{ 1}{ 0}{ 1}{ 0}{ 1}
\qquad
\E{1}{ 1}{ 0}{ 0}{ 1}{ 0}{ 0}{ 1}{ 0}\\
\hline

12& E_6 & 1& \E{\ast}{1}{ 1}{ 1}{ 1}{ 1}{ 1}{ \ast}{ \ast}&
\E{0}{ 0}{ 1}{ 0}{ 1}{ 0}{ 0}{ 1}{ 0}\qquad\E{0}{ 1}{ 0}{ 1}{ 0}{ 0}{ 1}{ 0}{ 1}
\qquad\E{0}{ 1}{ 1}{ 0}{ 0}{ 1}{ 0}{ 0}{ 1}\qquad\E{1}{ 0}{ 0}{ 0}{ 1}{ 0}{ 0}{ 1}{ 1}\\
&&&&\E{1}{ 0}{ 0}{ 1}{ 0}{ 1}{ 0}{ 0}{ 1}\qquad\E{1}{ 0}{ 1}{ 1}{ 0}{ 0}{ 1}{ 0}{ 0}
\qquad\E{1}{ 1}{ 0}{ 0}{ 0}{ 0}{ 1}{ 1}{ 1}\qquad\E {1}{ 1}{ 0}{ 0}{ 1}{ 0}{ 0}{ 1}{ 0}\\
&&&&\E{1}{ 1}{ 1}{ 0}{ 0}{ 0}{ 1}{ 0}{ 1}\\
\hline

10& E_8(a_6) & 2& \E{1}{ 0}{ 0}{ 0}{ 1}{ 0}{ 0}{ 1}{ 0}& \E{1}{ 0}{ 0}{ 0}{ 1}{ 0}{ 0}{ 1}{ 0}\\
\hline

10& D_6 & 1& \E{\ast}{\ast}{ 1}{ 1}{ 1}{ 1}{ 1}{ 1}{ \ast}&
\E{0}{ 0}{ 0}{ 1}{ 0}{ 0}{ 1}{ 0}{ 1}\qquad\E{0}{ 1}{ 0}{ 0}{ 1}{ 0}{ 0}{ 0}{ 1}
\qquad\E{0}{ 1}{ 0}{ 1}{ 0}{ 0}{ 1}{ 0}{ 0}\qquad\E{1}{ 0}{ 0}{ 0}{ 1}{ 0}{ 0}{ 1}{ 0}\\
&&&&\E{1}{ 1}{ 0}{ 0}{ 0}{ 1}{ 0}{ 0}{ 1}\qquad\E{1}{ 1}{ 1}{ 0}{ 0}{ 0}{ 1}{ 0}{ 0}\\
\hline

9& E_6(a_1) & 1& \E{\ast}{1}{ 1}{ 1}{ 0}{ 1}{ 1}{ \ast}{ \ast}&
\E{0}{ 0}{ 0}{ 0}{ 1}{ 0}{ 0}{ 1}{ 0}\qquad\E{0}{ 1}{ 0}{ 0}{ 0}{ 1}{ 0}{ 0}{ 1}
\qquad\E{1}{ 0}{ 0}{ 0}{ 1}{ 0}{ 0}{ 0}{ 1}\qquad\E{1}{ 0}{ 0}{ 1}{ 0}{ 0}{ 1}{ 0}{ 0}\\
&&&&\E{1}{ 1}{ 0}{ 0}{ 0}{ 0}{ 1}{ 0}{ 1}\qquad\E{1}{ 1}{ 1}{ 0}{ 0}{ 0}{ 0}{ 1}{ 0}
\\
\hline

9& A_8& 1& \E{1}{ 1}{ \ast}{ 1}{ 1}{ 1}{ 1}{ 1}{ 1}& \E{0}{ 0}{ 0}{ 0}{ 1}{ 0}{ 0}{ 1}{ 0}\\
\hline
8& D_8(a_3) & 2& \E{0}{ 0}{ 0}{ 0}{ 1}{ 0}{ 0}{ 0}{ 1}& \E{0}{ 0}{ 0}{ 0}{ 1}{ 0}{ 0}{ 0}{ 1}\\
\hline

8& D_6(a_1) & 1& \E{\ast}{\ast}{ 1}{ 1}{ 0}{ 1}{ 1}{ 1}{ \ast}&
\E{0}{ 0}{ 0}{ 0}{ 1}{ 0}{ 0}{ 0}{ 1}\qquad\E{0}{ 1}{ 1}{ 0}{ 0}{ 0}{ 0}{ 1}{ 0}
\qquad\E{1}{ 0}{ 0}{ 0}{ 0}{ 1}{ 0}{ 0}{ 1}\qquad\E{1}{ 0}{ 0}{ 1}{ 0}{ 0}{ 0}{ 1}{ 0}\\
&&&&\E{1}{ 1}{ 0}{ 0}{ 0}{ 1}{ 0}{ 0}{ 0}
\\
\hline

8& D_5 & 1& \E{\ast}{\ast}{ 1}{ 1}{ 1}{ 1}{ 1}{ \ast}{ \ast}&
\E{0}{ 0}{ 0}{ 0}{ 1}{ 0}{ 0}{ 0}{ 1}\qquad\E{0}{ 0}{ 0}{ 1}{ 0}{ 0}{ 1}{ 0}{ 0}
\qquad\E{0}{ 1}{ 0}{0}{ 0}{ 0}{ 1}{ 0}{ 1}\qquad\E{0}{ 1}{ 0}{ 1}{ 0}{ 0}{ 0}{ 0}{ 1}\\
&&&&\E{0}{ 1}{ 1}{ 0}{ 0}{ 0}{ 0}{ 1}{ 0}\qquad\E{1}{ 0}{ 0}{ 0}{ 0}{ 1}{ 0}{ 0}{ 1}
\qquad\E{1}{ 0}{ 0}{ 1}{ 0}{ 0}{ 0}{ 1}{ 0}\qquad\E{1}{0}{ 1}{ 0}{ 0}{ 0}{ 1}{ 0}{ 0}\\
&&&&\E{1}{ 1}{ 0}{ 0}{ 0}{ 0}{ 0}{ 1}{ 1}\qquad\E{1}{ 1}{ 0}{ 0}{ 0}{ 1}{ 0}{ 0}{ 0}
\qquad\E{1}{ 1}{ 1}{ 0}{ 0}{ 0}{ 0}{ 0}{ 1}
\\
\hline

8& A_7' & 1& \E{1}{\ast}{ 1}{ \ast}{ 1}{1}{1}{1}{1}& 
\E{0}{ 0}{ 0}{ 0}{ 1}{ 0}{ 0}{ 0}{ 1}\qquad \E{0}{ 0}{ 0}{ 1}{ 0}{ 0}{ 1}{ 0}{ 0}
\qquad\E{0}{ 1}{ 0}{ 0}{ 0}{ 0}{ 1}{ 0}{ 1}.
\\
\hline

8& A_7'' & 1& \E{\ast}{ 1}{\ast}{1}{1}{1}{1}{1}{1}& 
\E{0}{ 0}{ 0}{ 0}{ 1}{ 0}{ 0}{ 0}{ 1}
\\
\hline

7& A_6 & 1& \E{\ast}{ 1}{\ast}{1}{1}{1}{1}{1}{\ast}& 
\E{0}{ 0}{ 0}{ 0}{ 0}{ 1}{ 0}{ 0}{ 1}\qquad\E{0}{ 0}{ 0}{ 1}{ 0}{ 0}{ 0}{ 1}{ 0}
\qquad\E{1}{ 0}{ 0}{ 0}{ 1}{ 0}{ 0}{ 0}{ 0}\qquad\E{1}{ 1}{ 0}{ 0}{ 0}{ 0}{ 1}{ 0}{ 0}
\\
\hline
\end{array}
$$
\end{center}


\begin{center}
{\small Table 18 continued:  $\Kac(w)_{\un}$ and $\Kac(w)$ for $m$-admissible $w$ in $W(E_8)$ }
$$
\begin{array}{|c|c|c|l|l|}
\hline
m& w & r&\Kac(w)_{\un}& \Kac(w)\\
\hline\hline

6& E_8(a_8) & 4& \E{1}{ 0}{ 0}{ 0}{ 0}{ 1}{ 0}{ 0}{ 0}& \E{1}{ 0}{ 0}{ 0}{ 0}{ 1}{ 0}{ 0}{ 0}\\
\hline

6& E_7(a_4) & 3& \E{\ast}{ 0}{ 0}{ 0}{ 1}{ 0}{ 0}{ 1}{ \ast}& 
\E{0}{ 0}{ 0}{ 1}{ 0}{ 0}{ 0}{ 0}{ 1}\qquad \E{1}{ 0}{ 0}{ 0}{ 0}{ 1}{ 0}{ 0}{ 0}\\
\hline

6& E_6(a_2) & 2& \E{\ast}{ 1}{ 0}{ 0}{ 1}{ 0}{ 1}{ \ast}{ \ast}& 
\E{0}{ 0}{ 0}{ 1}{ 0}{ 0}{ 0}{ 0}{ 1}\qquad\E{0}{ 0}{ 1}{ 0}{ 0}{ 0}{ 0}{ 1}{ 0}
\qquad\E{1}{ 0}{ 0}{ 0}{ 0}{ 1}{ 0}{ 0}{ 0}\qquad\E{1}{ 1}{ 0}{ 0}{ 0}{ 0}{ 0}{ 1}{ 0}\\
\hline

6& D_6(a_2) & 2& \E{\ast}{ \ast}{ 1}{ 1}{ 0}{ 1}{ 0}{ 1}{ \ast}& 
\E{0}{ 0}{ 0}{ 1}{ 0}{ 0}{ 0}{ 0}{ 1}\qquad\E{0}{ 1}{ 0}{ 0}{ 0}{ 0}{ 1}{ 0}{ 0}
\qquad\E{1}{ 0}{ 0}{ 0}{ 0}{ 1}{ 0}{ 0}{ 0}\\
\hline

6& D_4 & 1& \E{\ast}{ \ast}{ 1}{1}{ 1}{ 1}{\ast}{ \ast}{ \ast}& 
\E{0}{ 0}{ 0}{ 1}{ 0}{ 0}{ 0}{ 0}{ 1}\qquad\E{0}{ 0}{ 1}{ 0}{ 0}{ 0}{ 0}{ 1}{ 0}
\qquad\E{0}{ 1}{ 0}{ 0}{ 0}{ 0}{ 1}{ 0}{ 0}\qquad\E{0}{ 1}{ 0}{ 1}{ 0}{ 0}{ 0}{ 0}{ 0}\\
&&&&\E{1}{ 0}{ 0}{ 0}{ 0}{ 0}{ 0}{ 1}{ 1}\qquad\E{1}{ 0}{ 0}{ 0}{ 0}{ 1}{ 0}{ 0}{ 0}
\qquad\E{1}{ 0}{ 1}{ 0}{ 0}{ 0}{ 0}{ 0}{ 1}\qquad\E{1}{ 1}{ 0}{ 0}{ 0}{ 0}{ 0}{ 1}{ 0}\\
&&&&\E{1}{ 1}{ 1}{ 0}{ 0}{ 0}{ 0}{ 0}{ 0}\\
\hline

6& A_5 & 1& \E{\ast}{1}{ \ast}{1}{ 1}{ 1}{1}{ \ast}{ \ast}& 
\E{0}{ 0}{ 0}{ 0}{ 0}{ 0}{ 1}{ 0}{ 1}\qquad\E{0}{ 0}{ 0}{ 0}{ 1}{ 0}{ 0}{ 0}{ 0}
\qquad\E{0}{ 0}{ 0}{ 1}{ 0}{ 0}{ 0}{ 0}{ 1}\qquad\E{0}{ 0}{ 1}{ 0}{ 0}{ 0}{ 0}{ 1}{ 0}\\
&&&&\E{0}{ 1}{ 0}{ 0}{ 0}{ 0}{ 1}{ 0}{ 0}\qquad\E{1}{ 0}{ 0}{ 0}{ 0}{ 1}{ 0}{ 0}{ 0}
\qquad\E{1}{ 1}{ 0}{ 0}{ 0}{ 0}{ 0}{ 1}{ 0}\\
\hline

5& 2A_4 & 2& \E{1}{1}{1}{1}{1}{\ast}{1}{1}{1}& \E{0}{ 0}{ 0}{ 0}{ 0}{ 1}{ 0}{ 0}{ 0}\\
\hline

5& A_4 & 1& \E{\ast}{1}{ \ast}{1}{ 1}{ 1}{\ast}{ \ast}{ \ast}& 
\E{0}{ 0}{ 0}{ 0}{ 0}{ 1}{ 0}{ 0}{ 0}\qquad\E{0}{ 0}{ 1}{ 0}{ 0}{ 0}{ 0}{ 0}{ 1}
\qquad\E{0}{ 1}{ 0}{ 0}{ 0}{ 0}{ 0}{ 1}{ 0}\qquad\E{1}{ 0}{ 0}{ 0}{ 0}{ 0}{ 1}{ 0}{ 0}\\
&&&&\E{1}{ 0}{ 0}{ 1}{ 0}{ 0}{ 0}{ 0}{ 0}\qquad\E{1}{ 1}{ 0}{ 0}{ 0}{ 0}{ 0}{ 0}{ 1}\\
\hline

4& 2D_4(a_1) & 4& \E{1}{ \ast}{ 1}{1}{ 0}{ 1}{\ast}{ 1}{ 0}& 
\E{0}{ 0}{ 0}{ 0}{ 0}{ 0}{ 1}{ 0}{ 0}\qquad\E{1}{ 0}{ 1}{ 0}{ 0}{ 0}{ 0}{ 0}{ 0}\\
\hline

4& D_4(a_1) & 2& \E{\ast}{ \ast}{ 1}{1}{ 0}{ 1}{\ast}{ \ast}{ \ast}& 
\E{0}{ 0}{ 0}{ 0}{ 0}{ 0}{ 1}{ 0}{ 0}\qquad\E{0}{ 0}{ 0}{ 1}{ 0}{ 0}{ 0}{ 0}{ 0}
\qquad\E{0}{ 1}{ 0}{ 0}{ 0}{ 0}{ 0}{ 0}{ 1}\qquad\E{1}{ 0}{ 0}{ 0}{ 0}{ 0}{ 0}{ 1}{ 0}\\
&&&&\E{1}{ 0}{ 1}{ 0}{ 0}{ 0}{ 0}{ 0}{ 0}\\
\hline

4& 2A_3' & 2& \E{\ast}{ 1}{ \ast}{1}{ \ast}{ 1}{1}{ 1}{ \ast}& 
\E{0}{ 0}{ 0}{ 0}{ 0}{ 0}{ 1}{ 0}{ 0}\qquad\E{0}{ 0}{ 0}{ 1}{ 0}{ 0}{ 0}{ 0}{ 0} 
\qquad\E{0}{ 1}{ 0}{ 0}{ 0}{ 0}{ 0}{ 0}{ 1} \qquad\E{1}{ 0}{ 1}{ 0}{ 0}{ 0}{ 0}{ 0}{ 0}\\
\hline

4& 2A_3''& 2& \E{\ast}{ 1}{ \ast}{1}{ 1}{ \ast}{1}{ 1}{1}& 
\E{0}{ 0}{ 0}{ 0}{ 0}{ 0}{ 1}{ 0}{ 0}\\
\hline

4& A_3& 1& \E{\ast}{ 1}{ \ast}{1}{ 1}{ \ast}{\ast}{ \ast}{\ast}& 
\E{0}{ 0}{ 0}{ 0}{ 0}{ 0}{ 1}{ 0}{ 0}\qquad\E{0}{ 0}{ 0}{ 1}{ 0}{ 0}{ 0}{ 0}{ 0}
\qquad\E{0}{ 1}{ 0}{ 0}{ 0}{ 0}{ 0}{ 0}{ 1}\qquad\E{1}{ 0}{ 0}{ 0}{ 0}{ 0}{ 0}{ 1}{ 0}\\
&&&&\E{1}{ 0}{ 1}{ 0}{ 0}{ 0}{ 0}{ 0}{ 0}\\
\hline

3& 4A_2& 4& \E{1}{ 1}{ 1}{1}{ \ast}{ 1}{1}{ \ast}{1}& 
\E{0}{ 0}{ 1}{ 0}{ 0}{ 0}{ 0}{ 0}{ 0}\\
\hline

3& 3A_2& 3& \E{1}{ 1}{ \ast}{1}{ \ast}{ 1}{1}{ \ast}{1}& 
\E{0}{ 0}{ 0}{ 0}{ 0}{ 0}{ 0}{ 1}{ 0}\qquad\E{0}{ 0}{ 1}{ 0}{ 0}{ 0}{ 0}{ 0}{ 0}\\
\hline

3& 2A_2& 2& \E{\ast}{ 1}{ \ast}{1}{ \ast}{ 1}{1}{ \ast}{\ast}& 
\E{0}{ 0}{ 0}{ 0}{ 0}{ 0}{ 0}{ 1}{ 0}\qquad\E{0}{ 0}{ 1}{ 0}{ 0}{ 0}{ 0}{ 0}{ 0}
\qquad\E{1}{ 1}{ 0}{ 0}{ 0}{ 0}{ 0}{ 0}{ 0}\\
\hline

3& A_2& 1&\E{\ast}{ 1}{ \ast}{1}{ \ast}{ \ast}{\ast}{ \ast}{\ast}& 
\E{0}{ 0}{ 0}{ 0}{ 0}{ 0}{ 0}{ 1}{ 0}\qquad\E{0}{ 0}{ 1}{ 0}{ 0}{ 0}{ 0}{ 0}{ 0}
\qquad\E{1}{ 0}{ 0}{ 0}{ 0}{ 0}{ 0}{ 0}{ 1}\qquad\E{1}{ 1}{ 0}{ 0}{ 0}{ 0}{ 0}{ 0}{ 0}\\
\hline

2& 8A_1& 7& \E{0}{ 1}{ 0}{ 0}{ 0}{ 0}{ 0}{ 0}{0}&\E{0}{ 1}{ 0}{ 0}{ 0}{ 0}{ 0}{ 0}{0}\\
\hline
2& 7A_1& 7& \E{\ast}{0}{ 1}{ 0}{ 0}{ 0}{ 0}{ 0}{\ast}&\E{0}{ 1}{ 0}{ 0}{ 0}{ 0}{ 0}{0}{0}\\
\hline
2& 6A_1& 6& \E{\ast}{\ast}{0}{0}{0}{1}{0}{0}{\ast}&\E{0}{ 1}{ 0}{ 0}{ 0}{ 0}{ 0}{0}{0}\\
\hline
2& 5A_1& 5& \E{1}{1}{1}{\ast}{\ast}{1}{\ast}{1}{\ast}&\E{0}{ 1}{ 0}{ 0}{ 0}{ 0}{ 0}{0}{0}\\
\hline
2& 4A_1'& 4& \E{\ast}{\ast}{1}{1}{\ast}{1}{\ast}{1}{\ast}& 
\E{0}{ 1}{ 0}{ 0}{ 0}{ 0}{ 0}{ 0}{ 0}\\
\hline
2& 4A_1''& 4& \E{\ast}{\ast}{0}{0}{1}{0}{\ast}{\ast}{\ast}& 
\E{0}{ 0}{ 0}{ 0}{ 0}{ 0}{ 0}{ 0}{ 1}\qquad\E{0}{ 1}{ 0}{ 0}{ 0}{ 0}{ 0}{ 0}{ 0}\\
\hline
2& 3A_1& 3& \E{\ast}{\ast}{1}{1}{\ast}{1}{\ast}{\ast}{\ast}& 
\E{0}{ 0}{ 0}{ 0}{ 0}{ 0}{ 0}{ 0}{ 1}\qquad\E{0}{ 1}{ 0}{ 0}{ 0}{ 0}{ 0}{ 0}{ 0}\\
\hline

2& 2A_1& 2& \E{\ast}{\ast}{1}{1}{\ast}{\ast}{\ast}{\ast}{\ast}& 
\E{0}{ 0}{ 0}{ 0}{ 0}{ 0}{ 0}{ 0}{ 1}\qquad\E{0}{ 1}{ 0}{ 0}{ 0}{ 0}{ 0}{ 0}{ 0}\\
\hline
2& A_1& 1& \E{\ast}{\ast}{1}{\ast}{\ast}{\ast}{\ast}{\ast}{\ast}& 
\E{0}{ 0}{ 0}{ 0}{ 0}{ 0}{ 0}{ 0}{ 1}\qquad\E{0}{ 1}{ 0}{ 0}{ 0}{ 0}{ 0}{ 0}{ 0}\\
\hline

\end{array}
$$
\end{center}

\subsection{Tables of positive rank gradings for $E_6, E_7$ and $E_8$}\label{Eposrank}
The previous lists contain the Kac coordinates of all positive rank gradings, usually with multiple occurrences. We now discard those in each $\Kac(w)$ which appear in some $\Kac(w')$ with 
$\rank(w')>\rank(w)$. The remaining elements of $\Kac(w)$ are then the Kac coordinates 
of automorphisms $\theta$ of order $m$ with $\rank(\theta)=\rank(w)$. For each grading $\theta$ there still may be more than one  $w$ with $\rank(\theta)=\rank(w)$. It turns out that every $\theta$ of positive rank is contained in $\Kac(w)_{\un}$ for some $m$-admissible $w$ which is a $\bz$-regular element in the Weyl group Levi of a Levi subgroup $L_\theta$ and $\theta$ is a principal inner automorphism of the Lie algebra of $L_\theta$. This Levi subgroup  corresponds to the subset $J$ of Lemma \ref{Jlowerbound} and is indicated in the right most column of the tables below. 
For example, in $E_7$ the Kac diagrams 
$$8_a:\ \EVII {1}{ 0}{ 0}{ 0}{ 1}{ 0}{ 1}{ 1}\qquad\text{and}\qquad
 8_b:\ \EVII{0}{ 1}{ 1}{ 0}{ 0}{ 1}{ 0}{ 1}
$$
occur in  $\Kac(w)$ for $w$ of types $D_6(a_1)$ and $D_5$. 
Since $D_6(a_1)$ is not regular in any Levi subgroup of $W(E_7)$ and $w=D_5$ is regular in the $D_5$ Levi subgroup, we choose $w=D_5$, discard $w=D_6(a_1)$, and set $L_\theta=D_5$. 

Since $\theta$ is principal on the Lie algebra of $L_\theta$, there is a conjugate $\theta'$ of $\theta$ whose un-normalized Kac diagram has a $1$ on each node of $J$ (cf. Lemma \ref{Jlowerbound}).
There may be more than one such $J$, corresponding to various conjugates $\theta'$, and we just pick one of them. 

In the tables we try to  write $w$ in a  form which exhibits its regularity in the Weyl group 
$W_J$. For example,  in $E_6$ the gradings $4_a, 4_b$ have $w=D_4(a_1)$. 
\footnote{cf. Panyushev, Example 4.5.} 
In case $4_a$, which is stable,
 we give the alternate expression $w=E_6^3$ to make it clear that $w$ is $\bz$-regular in $W(E_6)$. In case $4_b$ there is no $W_{\aff}(R)$-conjugate of $\theta$ with $1$'s on the  $E_6$ subdiagram. However, $w=D_5^2$ is the square of a Coxeter element in $W_{D_5}$, hence is $\bz$-regular in $W_{D_5}$. 

The rows in our tables are ordered by decreasing  $m$. 
The positive rank inner gradings of a given order $m$ are named 
$m_a, m_b, m_c,\dots,$ where $m_a$ is the unique principal grading of order $m$.  
The principal grading $m_a$ has maximal rank and minimal dimension of $\fg_0$ among all gradings of order $m$.
The remaining rows of order $m$ are grouped according to $w$ and $L_\theta$, 
ordered in each group by increasing dimension of $\fg_0$. 

The little Weyl groups $W(\fc,\theta)$ are also given, along with their degrees. 
These are either cyclic or given by their notation in \cite{shephard-todd}.
We explain their computation in section \ref{littleweylE}.

\begin{center}
{\small Table 19: The gradings of positive rank in type $E_6$ (inner case)}
$$
\begin{array}{c c c c  c c c c}
\hline
\text{No.}&\text{Kac diagram}& w& W(\fc,\theta)&\text{degrees}&\theta'&L_\theta\\
\hline
12_a &\EVI{1}{1}{1}{1}{1}{1}{1}  &E_6 &\bfmu_{12}&12&\EVI{1}{1}{1}{1}{1}{1}{1} &E_6\\
9_a& \EVI{1}{1}{1}{1}{0}{1}{1} &E_6(a_1) &\bfmu_{9}&9&\EVI{-2}{1}{1}{1}{1}{1}{1} &E_6\\
8_a&\EVI{1}{ 1}{ 1}{ 0}{ 1}{ 0}{ 1} &D_5 &\bfmu_{8}&8&\EVI{-3}{1}{1}{1}{1}{1}{1}&E_6\\
8_b&\EVI{0}{ 1}{ 1}{ 1}{ 0}{ 1}{ 1}&D_5 &\bfmu_8&8&\EVI{-6}{1}{1}{1}{1}{1}{0}&D_5\\
6_a&\EVI{1}{ 1}{ 0}{ 0}{ 1}{ 0}{ 1}  &E_6(a_2) &G_5&6,12&\EVI{-5}{1}{1}{1}{1}{1}{1}&E_6\\
6_b,{6_b}'&\EVI{0}{ 1}{ 1}{ 1}{ 0}{ 0}{ 1} \qquad\EVI{0}{ 1}{ 1}{ 0}{ 0}{ 1}{ 1}
 &D_4 &\bfmu_6&6&\EVI{-4}{0}{1}{1}{1}{1}{1}&D_4\\
 
6_c&\EVI{1}{ 0}{ 1}{ 0}{ 1}{ 0}{ 0}  &D_4 &\bfmu_{6}&6&\EVI{-3}{0}{1}{1}{1}{1}{0}&D_4\\

6_d&\EVI{0}{ 0}{ 1}{ 1}{ 0}{ 1}{ 0} &A_5     &\bfmu_6&6&\EVI{-3}{1}{0}{1}{1}{1}{1}&A_5\\

5_a&\EVI{0}{ 1}{ 0}{ 0}{ 1}{ 0}{ 1}  &A_4 &\bfmu_5&5&\EVI{-6}{1}{1}{1}{1}{1}{1}&A_5\\
 
5_b&\EVI{1}{ 0}{ 0}{ 1}{ 0}{ 1}{ 0}&A_4 &\bfmu_5&5&\EVI{-8}{1}{2}{1}{1}{1}{1}&A_5\\
 
5_c&\EVI{1}{ 1}{ 1}{ 0}{ 0}{ 0}{ 1}  &A_4 &\bfmu_5&5&\EVI{-10}{1}{3}{1}{1}{1}{1}&A_5\\
 
4_a&\EVI{1}{ 0}{ 0}{ 0}{ 1}{ 0}{ 0}  &D_4(a_1)=E_6^3 &G_8&8,12&\EVI{-7}{1}{1}{1}{1}{1}{1}
&E_6\\
 
4_b&\EVI{0}{ 1}{ 1}{ 0}{ 0}{ 0}{ 1}  &D_4(a_1)=D_5^2 &G(4,1,2)&4,8&\EVI{-6}{0}{1}{1}{1}{1}{1}&D_5\\
 
4_c&\EVI{0}{ 0}{ 0}{ 1}{ 0}{ 1}{ 0} &A_3 &\bfmu_4&4&\EVI{-6}{2}{1}{0}{1}{1}{1}&A_4\\
 
4_d,{4_d}'&\EVI{0}{ 1}{ 0}{ 1}{ 0}{ 0}{ 1}\qquad\EVI{0}{ 1}{ 0}{ 0}{ 0}{ 1}{ 1} &A_3&
\bfmu_{4}&4&\EVI{-4}{0}{0}{1}{1}{1}{1}&A_4\\
 
3_a&\EVI{0}{ 0}{ 0}{ 0}{ 1}{ 0}{ 0} &3A_2&G_{25}&6,9,12&\EVI{-8}{1}{1}{1}{1}{1}{1}&E_6\\
 
3_b&\EVI{1}{ 1}{ 0}{ 0}{ 0}{ 0}{ 1}&2A_2=A_5^2&G(3,1,2)&3,6&\EVI{-6}{1}{0}{1}{1}{1}{1}&A_5\\
 
3_c&\EVI{1}{ 0}{ 1}{ 0}{ 0}{ 0}{ 0} &A_2=D_4^2&\bfmu_6&6&\EVI{-6}{0}{1}{1}{1}{1}{0}&D_4\\
 
3_d,{3_d}' &\EVI{0}{ 1}{ 0}{ 0}{ 0}{ 1}{ 0}\qquad\EVI{0}{ 0}{ 0}{ 1}{ 0}{ 0}{ 1} &A_2=D_4^2&\bfmu_6&6&
\EVI{-7}{0}{1}{1}{1}{1}{1}&D_4\\
 
2_a&\EVI{0}{ 0}{ 1}{ 0}{ 0}{ 0}{ 0} &4A_1=E_6^6&W(F_4)&2,6,8,12&\EVI{-9}{1}{1}{1}{1}{1}{1}&E_6\\
 
2_b&\EVI{0}{ 1}{ 0}{ 0}{ 0}{ 0}{ 1} &2A_1=A_3^2&W(B_2)&2,4&\EVI{-6}{0}{0}{1}{1}{1}{1}&A_3\\
 
1_a&\EVI{1}{ 0}{ 0}{ 0}{ 0}{ 0}{ 0} &e&W(E_6)&2,5,6,8,9,12&\EVI{1}{ 0}{ 0}{ 0}{ 0}{ 0}{ 0}&\varnothing\\
 \hline
\end{array}
$$
\end{center}

\begin{center}
{\small Table 20: The gradings of positive rank in type $E_7$ }
$$
\begin{array}{cccccccc}
\hline
\text{No.}&\text{Kac diagram }&  w& W(\fc,\theta)&\text{degrees}&\theta'&L_\theta\\
\hline
18_a&\EVII{1}{1}{1}{1}{1}{1}{1}{1}  &E_7 &\bfmu_{18}&18&\EVII{1}{1}{1}{1}{1}{1}{1}{1} &E_7\\
14_a&\EVII{1}{1}{1}{1}{0}{1}{1}{1} &E_7(a_1) =-A_6&\bfmu_{14}&14&\EVII{-3}{1}{1}{1}{1}{1}{1}{1} &E_7\\

12_a&\EVII{1}{1}{1}{1}{0}{1}{0}{1}  &E_6 &\bfmu_{12}&12&\EVII{-5}{1}{1}{1}{1}{1}{1}{1} &E_7\\
12_b&\EVII{1}{ 1}{ 1}{ 0}{ 1}{ 0}{ 1}{ 1}  &E_7(a_2)=-E_6 &\bfmu_{12}&12&\EVII{-6}{1}{1}{1}{1}{1}{1}{2} &E_6\\
12_c&\EVII{0}{ 1}{ 0}{ 0}{ 1}{ 1}{ 1}{ 1}  &E_6 &\bfmu_{12}&12&\EVII{-4}{1}{1}{1}{1}{1}{1}{0} &E_6\\

10_a&\EVII{1}{ 0}{ 1}{ 1}{ 0}{ 1}{ 0}{ 1}  &D_6 &\bfmu_{10}&10&\EVII{-7}{1}{1}{1}{1}{1}{1}{1} &
D_6\\
10_b&\EVII {1}{ 1}{ 0}{ 0}{ 1}{ 0}{ 1}{ 1}  &D_6 &\bfmu_{10}&10&\EVII{-9}{2}{1}{1}{1}{1}{1}{1} &D_6\\
10_c&\EVII{0}{ 1}{ 1}{ 0}{ 1}{ 0}{ 1}{ 0} &D_6 &\bfmu_{10}&10&\EVII{-5}{0}{1}{1}{1}{1}{1}{1} &D_6\\

9_a&\EVII{0}{ 1}{ 0}{ 0}{ 1}{ 0}{ 1}{ 1}  &E_6(a_1)=E_7^2 &\bfmu_{18}&18&
\EVII{-8}{1}{1}{1}{1}{1}{1}{1} &E_7\\

9_b&\EVII{1}{ 0}{ 1}{ 1}{ 0}{ 0}{ 1}{ 1}   &E_6(a_1) &\bfmu_{9}&9&
\EVII{-7}{1}{1}{1}{1}{1}{1}{0} & E_6\\
8_a&\EVII {1}{ 0}{ 0}{ 0}{ 1}{ 0}{ 1}{ 1}  &D_5 &\bfmu_{8}&8&
\EVII{-9}{1}{1}{1}{1}{1}{1}{1} &D_5\\
8_b&\EVII {0}{ 1}{ 1}{ 0}{ 0}{ 1}{ 0}{ 1}&D_5 &\bfmu_{8}&8&
\EVII{-11}{2}{1}{1}{1}{1}{1}{1} &D_5\\

8_c&\EVII {0}{ 1}{ 0}{ 0}{ 1}{ 0}{ 1}{ 0}   &D_5 &\bfmu_{8}&8&
\EVII{-12}{0}{1}{1}{1}{1}{1}{6} &D_5\\

8_d&\EVII {1}{ 0}{ 0}{ 1}{ 0}{ 1}{ 0}{ 1}  &D_5 &\bfmu_{8}&8&
\EVII{-10}{1}{1}{1}{1}{1}{1}{2} &D_5\\

8_e&\EVII {1}{ 1}{ 1}{ 0}{ 0}{ 0}{ 1}{ 1}  &D_5 &\bfmu_{8}&8&
\EVII{-8}{0}{1}{1}{1}{1}{1}{2} &D_5\\

8_f&\EVII{1}{ 0}{ 1}{ 0}{ 1}{ 0}{ 0}{ 1}  &D_5 &\bfmu_{8}&8&
\EVII{-12}{1}{1}{1}{1}{1}{1}{4} &D_5\\

8_g&\EVII {0}{ 0}{ 1}{ 1}{ 0}{ 0}{ 1}{ 1} &D_5 &\bfmu_{8}&8&
\EVII{-6}{0}{1}{1}{1}{1}{1}{0} &D_5\\

8_h&\EVII{0}{ 1}{ 0}{ 0}{ 0}{ 1}{ 1}{ 1}   &D_5 &\bfmu_{8}&8&
\EVII{-8}{1}{1}{1}{1}{1}{1}{0} &D_5\\

7_a&\EVII{0}{ 1}{ 0}{ 0}{ 1}{ 0}{ 0}{ 1} &A_6=E_7(a_1)^2 &\bfmu_{14}&14&
\EVII{-10}{1}{1}{1}{1}{1}{1}{1} &E_7\\

6_a& \EVII{1}{ 0}{ 0}{ 0}{ 1}{ 0}{ 0}{ 1}  &E_7(a_4)=E_7^3=-3A_2&G_{26}&6,12,18&
\EVII{-11}{1}{1}{1}{1}{1}{1}{1} &E_7\\

6_b&\EVII{0}{ 1}{ 0}{ 0}{ 0}{ 1}{ 0}{ 1} &E_6(a_2)=E_6^2&G_5&6,12&
\EVII{-10}{1}{1}{1}{1}{1}{1}{0} &E_6\\

6_c& \EVII{0}{ 1}{ 1}{ 0}{ 0}{ 0}{ 1}{ 0} &D_6(a_2)&G(6,2,2)&6,6&
\EVII{-9}{0}{1}{1}{1}{1}{1}{1} &D_6\\

6_d& \EVII{1}{ 0}{ 1}{ 0}{ 0}{ 0}{ 1}{ 1}  &D_4&\bfmu_{6}&6&
\EVII{-12}{0}{1}{1}{1}{1}{1}{2} &D_4\\

6_e&\EVII {0}{ 0}{ 1}{ 1}{ 0}{ 0}{ 0}{ 1} &D_4&\bfmu_{6}&6&
\EVII{-10}{0}{1}{1}{1}{1}{1}{2} &D_4\\

6_f&\EVII{0}{ 0}{ 0}{ 1}{ 0}{ 0}{ 1}{ 1}  &D_4&\bfmu_{6}&6&
\EVII{-13}{0}{1}{1}{1}{1}{1}{5} &D_4\\

6_g&\EVII{0}{ 0}{ 1}{ 0}{ 1}{ 0}{ 0}{ 0} &D_4&\bfmu_{6}&6&
\EVII{-9}{0}{1}{1}{1}{1}{0}{3} &D_4\\

6_h&\EVII{0}{ 0}{ 0}{ 0}{ 0}{ 1}{ 1}{ 1}  &D_4&\bfmu_{6}&6&
\EVII{-6}{0}{1}{1}{1}{1}{0}{0} &D_4\\

6_i&\EVII{0}{ 0}{ 0}{ 1}{ 0}{ 1}{ 0}{ 0} &A_5' &\bfmu_{6}&6&
\EVII{-10}{2}{1}{0}{1}{1}{1}{1} &A_5' \\

6_j&\EVII{0}{ 0}{ 0}{ 0}{ 1}{ 0}{ 1}{ 0} &A_5''&\bfmu_{6}&6&
\EVII{-9}{1}{0}{1}{1}{1}{1}{1} &A_5''\\

6_k&\EVII{1}{ 1}{ 0}{ 0}{ 0}{ 0}{ 1}{ 1} &A_5' &\bfmu_{6}&6&
\EVII{-6}{0}{1}{0}{1}{1}{1}{1} &A_5' \\

\end{array}
$$
\end{center}

\begin{center}
{\small Table 20 continued: The gradings of positive rank in type $E_7$}
$$
\begin{array}{cccccccc}
\hline
\text{No.}&\text{Kac diagram}& w& W(\fc,\theta)
&\text{degrees}&\theta'&L_\theta\\
\hline

5_a&\EVII{0}{ 0}{ 0}{ 0}{ 1}{ 0}{ 0}{ 1} &A_4=D_6^2&\bfmu_{10}&10&
\EVII{-12}{1}{1}{1}{1}{1}{1}{1} &D_6\\

5_b&\EVII{0}{ 0}{ 0}{ 1}{ 0}{ 0}{ 1}{ 0}  &A_4=D_6^2&\bfmu_{10}&10&
\EVII{-14}{2}{1}{1}{1}{1}{1}{1} &D_6\\

5_c&\EVII{0}{ 1}{ 0}{ 0}{ 0}{ 0}{ 1}{ 1} &A_4=D_6^2&\bfmu_{10}&10&
\EVII{-10}{0}{1}{1}{1}{1}{1}{1} &D_6\\

5_d&\EVII{1}{ 0}{ 0}{ 0}{ 0}{ 1}{ 0}{ 1} &A_4&\bfmu_{5}&5&
\EVII{-11}{1}{0}{1}{1}{1}{1}{0} & A_4\\

5_e&\EVII{0}{ 1}{ 1}{ 0}{ 0}{ 0}{ 0}{ 1}  &A_4&\bfmu_{5}&5&
\EVII{-11}{1}{1}{1}{1}{1}{0}{2} &A_4\\

4_a&\EVII{0}{ 0}{ 0}{ 1}{ 0}{ 0}{ 0}{ 1}  &D_4(a_1)=E_6^3&G_8&8,12&
\EVII{-13}{1}{1}{1}{1}{1}{1}{1} &E_6\\

4_b&\EVII{0}{ 0}{ 0}{ 0}{ 1}{ 0}{ 0}{ 0} &D_4(a_1)=E_6^3&G_8&8,12&
\EVII{-14}{1}{1}{1}{1}{1}{1}{2} &E_6\\

4_c&\EVII{0}{ 0}{ 0}{ 0}{ 0}{ 1}{ 0}{ 1} &D_4(a_1)=E_6^3&G_8&8,12&
\EVII{-12}{1}{1}{1}{1}{1}{1}{0} &E_6\\

4_d&\EVII{0}{ 1}{ 0}{ 0}{ 0}{ 0}{ 1}{ 0} &D_4(a_1)&G(4,1,2)&4,8&
\EVII{-12}{0}{1}{1}{1}{1}{1}{2} &D_5\\

4_e&\EVII{1}{ 0}{ 1}{ 0}{ 0}{ 0}{ 0}{ 1} &D_4(a_1)&G(4,1,2)&4,8&
\EVII{-12}{1}{1}{1}{1}{1}{0}{2} &D_5\\

4_f&\EVII{0}{ 0}{ 1}{ 0}{ 0}{ 0}{ 1}{ 0}  &A_3&\bfmu_{4}&4&
\EVII{-9}{1}{1}{1}{1}{0}{0}{2} &A_4\\

4_g&\EVII{1}{ 0}{ 0}{ 0}{ 0}{ 0}{ 1}{ 1}  &A_3&\bfmu_{4}&4&
\EVII{-7}{1}{1}{1}{1}{0}{0}{0} &A_4\\

3_a&\EVII{0}{ 0}{ 0}{ 0}{ 0}{ 1}{ 0}{ 0} &3A_2=E_7^6&G_{26}&6,12,18&
\EVII{-14}{1}{1}{1}{1}{1}{1}{1} &E_7\\

3_b&\EVII{0}{ 1}{ 0}{ 0}{ 0}{ 0}{ 0}{ 1} &2A_2&G(6,2,2)&6,6&
\EVII{-12}{0}{1}{1}{1}{1}{1}{1} &D_6\\

3_c&\EVII{0}{ 0}{ 1}{ 0}{ 0}{ 0}{ 0}{ 1}  &A_2=D_4^2&\bfmu_{6}&6&
\EVII{-13}{0}{1}{1}{1}{1}{1}{2} &D_4\\

3_d&\EVII{0}{ 0}{ 0}{ 0}{ 0}{ 0}{ 1}{ 1}  &A_2=D_4^2&\bfmu_{6}&6&
\EVII{-9}{0}{1}{1}{1}{1}{0}{0} &D_4\\

2_a&\EVII{0}{ 0}{ 1}{ 0}{ 0}{ 0}{ 0}{ 0} &7A_1&W(E_7)&2,6,8,10,12,14,18&
\EVII{-15}{1}{1}{1}{1}{1}{1}{1} &E_7\\

2_b&\EVII{0}{ 0}{ 0}{ 0}{ 0}{ 0}{ 1}{ 0} &4A_1''&W(F_4)&2,6,8,12&
\EVII{-14}{1}{1}{1}{1}{1}{1}{0} &E_6\\

2_c&\EVII{1}{ 0}{ 0}{ 0}{ 0}{ 0}{ 0}{ 1}  &3A_1'&W(B_3)&2,4,6&
\EVII{-10}{0}{0}{1}{1}{1}{1}{1} &A_5'\\

1_a&\EVII{1}{ 0}{ 0}{ 0}{ 0}{ 0}{ 0}{ 0}  &e&W(E_7)&2,6,8,10,12,14,18&
\EVII{1}{0}{0}{0}{0}{0}{0}{0} &E_7\\
\hline
\end{array}
$$
\end{center}

\begin{center}
{\small Table 21: The gradings of positive rank in type $E_8$}
$$
\begin{array}{cccccccc}
\hline
\text{No.}&\text{Kac diagram}&  w& W(\fc,\theta)&\text{degrees}&\theta'&L_\theta\\
\hline
30_a&\E{1}{ 1}{ 1}{ 1}{ 1}{ 1}{ 1}{ 1}{ 1}&E_8 & \bfmu_{30}&30&
\E{1}{ 1}{ 1}{ 1}{ 1}{ 1}{ 1}{ 1}{ 1} &E_8\\

24_a&\E{1}{ 1}{ 1}{ 1}{ 0}{ 1}{ 1}{ 1}{ 1}&E_8(a_1) & \bfmu_{24}&24&
\E{-5}{ 1}{ 1}{ 1}{ 1}{ 1}{ 1}{ 1}{ 1} &E_8\\

20_a&\E{1}{ 1}{ 1}{ 1}{ 0}{ 1}{ 0}{ 1}{ 1}&E_8(a_2) & \bfmu_{20}&20&
\E{-9}{ 1}{ 1}{ 1}{ 1}{ 1}{ 1}{ 1}{1} &E_8\\

18_a&\E{1}{ 1}{ 1}{ 1}{ 0}{ 1}{ 0}{ 1}{ 0}&E_7 & \bfmu_{18}&18&
\E{-11}{ 1}{ 1}{ 1}{ 1}{ 1}{ 1}{ 1}{1}&E_7\\

18_b&\E{1}{ 1}{ 1}{ 0}{ 1}{ 0}{ 1}{ 0}{ 1}&E_7 & \bfmu_{18}&18&
\E{-13}{ 1}{ 1}{ 1}{ 1}{ 1}{ 1}{1}{ 2}&E_7\\

18_c&\E{1}{ 1}{ 0}{ 0}{ 1}{ 0}{ 1}{ 1}{ 1}&E_7 & \bfmu_{18}&18&
\E{-17}{ 1}{ 1}{ 1}{ 1}{ 1}{ 1}{1}{ 4}&E_7\\

18_d&\E{1}{ 0}{ 1}{ 1}{ 0}{ 1}{ 0}{ 1}{ 1}&E_7& \bfmu_{18}&18&
\E{-15}{ 1}{ 1}{ 1}{ 1}{ 1}{ 1}{1}{ 3} &E_7\\

18_e&\E{0}{ 1}{ 0}{ 1}{ 1}{ 0}{ 1}{ 0}{ 1}&E_7 & \bfmu_{18}&18&
\E{-9}{ 1}{ 1}{ 1}{ 1}{ 1}{ 1}{1}{ 0} &E_7\\

15_a&\E{1}{ 1}{ 0}{ 0}{ 1}{ 0}{ 1}{ 0}{ 1}&E_8(a_5) & \bfmu_{30}&30&
\E{-14}{ 1}{ 1}{ 1}{ 1}{ 1}{ 1}{ 1}{ 1} &E_8\\

14_a&\E{1}{ 1}{ 0}{ 0}{ 1}{ 0}{ 0}{ 1}{ 1}&E_7(a_1) & \bfmu_{14}&14&
\E{-15}{ 1}{ 1}{ 1}{ 1}{ 1}{ 1}{ 1}{ 1} &E_7\\

14_b&\E{1}{ 1}{ 1}{ 0}{ 0}{ 1}{ 0}{ 1}{ 0}&E_7(a_1)  & \bfmu_{14}&14&
\E{-19}{ 1}{ 1}{ 1}{ 1}{ 1}{ 1}{1}{ 3} &E_7\\

14_c&\E{1}{ 0}{ 1}{ 1}{ 0}{ 0}{ 1}{ 0}{ 1}&E_7(a_1)  & \bfmu_{14}&14&
\E{-17}{ 1}{ 1}{ 1}{ 1}{ 1}{ 1}{1}{ 2} &E_7\\

14_d&\E{0}{ 1}{ 0}{ 0}{ 1}{ 0}{ 1}{ 0}{ 1}&E_7(a_1)
& \bfmu_{14}&14&
\E{-21}{ 1}{ 1}{ 1}{ 1}{ 1}{ 1}{1}{ 4} &E_7\\

12_a&\E{1}{ 1}{ 0}{ 0}{ 1}{ 0}{ 0}{ 1}{ 0}&E_8(a_3)  & G_{10}&12,24&
\E{-17}{ 1}{ 1}{ 1}{ 1}{ 1}{ 1}{ 1}{ 1} &E_8\\

12_b&\E{1}{ 1}{ 1}{ 0}{ 0}{ 0}{ 1}{ 0}{ 1}&E_6  & \bfmu_{12}&12&
\E{-15}{ 1}{ 1}{ 1}{ 1}{ 1}{ 1}{1}{ 0} &E_6\\

12_c&\E{1}{ 0}{ 1}{ 1}{ 0}{ 0}{ 1}{ 0}{ 0}&E_6  & \bfmu_{12}&12&
\E{-21}{ 1}{ 1}{ 1}{ 1}{ 1}{ 1}{ 1}{ 3} &E_6\\

12_d&\E{1}{ 0}{ 0}{ 1}{ 0}{ 1}{ 0}{ 0}{ 1}&E_6  &\bfmu_{12}&12&
\E{-16}{ 1}{ 1}{ 1}{ 1}{ 1}{ 1}{ 0}{ 2} &E_6\\

12_e&\E{1}{ 0}{ 0}{ 0}{ 1}{ 0}{ 0}{ 1}{ 1}&E_6  & \bfmu_{12}&12&
\E{-14}{ 1}{ 1}{ 1}{ 1}{ 1}{ 1}{ 0}{ 1} &E_6\\

12_f&\E{0}{ 1}{ 1}{ 0}{ 0}{ 1}{ 0}{ 0}{ 1}&E_6  & \bfmu_{12}&12&
\E{-19}{ 1}{ 1}{ 1}{ 1}{ 1}{ 1}{ 1}{ 2} &E_6\\

12_g&\E{0}{ 1}{ 0}{ 1}{ 0}{ 0}{ 1}{ 0}{ 1}&E_6  & \bfmu_{12}&12&
\E{-24}{ 1}{ 1}{ 1}{ 1}{ 1}{ 1}{ 0}{ 6} &E_6\\

12_h&\E{0}{ 0}{ 1}{ 0}{ 1}{ 0}{ 0}{ 1}{ 0}&E_6  & \bfmu_{12}&12&
\E{-20}{ 1}{ 1}{ 1}{ 1}{ 1}{ 1}{ 0}{ 4} &E_6\\

12_i&\E{1}{ 1}{ 0}{ 0}{ 0}{ 0}{ 1}{ 1}{ 1}&E_6  & \bfmu_{12}&12&
\E{-12}{ 1}{ 1}{ 1}{ 1}{ 1}{ 1}{ 0}{ 0} &E_6\\

12_j&\E{0}{ 0}{ 0}{ 0}{ 1}{ 0}{ 1}{ 0}{ 1}&D_7  & \bfmu_{12}&12&
\E{-15}{ 0}{ 1}{ 1}{ 1}{ 1}{ 1}{ 1}{ 1} &D_7\\

10_a&\E{1}{ 0}{ 0}{ 0}{ 1}{ 0}{ 0}{ 1}{ 0}&E_8(a_6)=-2A_4  & G_{16}&20,30&
\E{-19}{ 1}{ 1}{ 1}{ 1}{ 1}{ 1}{ 1}{ 1} &E_8\\

10_b&\E{1}{ 1}{ 1}{ 0}{ 0}{ 0}{ 1}{ 0}{ 0}&D_6  & \bfmu_{10}&10&
\E{-17}{ 1}{ 1}{ 1}{ 1}{ 1}{ 1}{ 1}{ 0} &D_6\\

10_c&\E{1}{ 1}{ 0}{ 0}{ 0}{ 1}{ 0}{ 0}{ 1}&D_6  & \bfmu_{10}&10&
\E{-21}{ 2}{ 1}{ 1}{ 1}{ 1}{ 1}{ 1}{ 1} &D_6\\

10_d&\E{0}{ 1}{ 0}{ 0}{ 1}{ 0}{ 0}{ 0}{ 1}&D_6  & \bfmu_{10}&10&
\E{-17}{ 0}{ 1}{ 1}{ 1}{ 1}{ 1}{ 1}{ 1} &D_6\\

10_e&\E{0}{ 0}{ 0}{ 1}{ 0}{ 0}{ 1}{ 0}{ 1}&D_6  & \bfmu_{10}&10&
\E{-19}{ 0}{ 1}{ 1}{ 1}{ 1}{ 1}{ 1}{2} &D_6\\

10_f&\E{0}{ 1}{ 0}{ 1}{ 0}{ 0}{ 1}{ 0}{ 0}&D_6  & \bfmu_{10}&10&
\E{-15}{ 0}{ 1}{ 1}{ 1}{ 1}{ 1}{ 1}{ 0} &D_6\\
\end{array}
$$
\end{center}

\begin{center}
{\small Table 21 continued: The gradings of positive rank in type $E_8$ }
$$
\begin{array}{cccccccc}
\hline
\text{No.}&\text{Kac diagram}&  w& W(\fc,\theta)&\text{degrees}&\theta'&L_\theta\\
\hline
9_a&\E{1}{ 0}{ 0}{ 0}{ 1}{ 0}{ 0}{ 0}{ 1}&E_6(a_1)=E_7^2  &\bfmu_{18}&18&
\E{-20}{ 1}{ 1}{ 1}{ 1}{ 1}{ 1}{ 1}{ 1} &E_7\\

9_b&\E{1}{ 0}{ 0}{ 1}{ 0}{ 0}{ 1}{ 0}{ 0}&E_6(a_1)=E_7^2  &\bfmu_{18}&18&
\E{-22}{ 1}{ 1}{ 1}{ 1}{ 1}{ 1}{ 1}{ 2} &E_7\\

9_c&\E{0}{ 1}{ 0}{ 0}{ 0}{ 1}{ 0}{ 0}{ 1}&E_6(a_1)=E_7^2  &\bfmu_{18}&18&
\E{-26}{ 1}{ 1}{ 1}{ 1}{ 1}{ 1}{ 1}{ 4} &E_7\\

9_d&\E{0}{ 0}{ 0}{ 0}{ 1}{ 0}{ 0}{ 1}{ 0}&E_6(a_1)=E_7^2  &\bfmu_{18}&18&
\E{-24}{ 1}{ 1}{ 1}{ 1}{ 1}{ 1}{ 1}{ 3} &E_7\\

9_e&\E{1}{ 1}{ 0}{ 0}{ 0}{ 0}{ 1}{ 0}{ 1}&E_6(a_1) =E_7^2 &\bfmu_{18}&18&
\E{-18}{ 1}{ 1}{ 1}{ 1}{ 1}{ 1}{ 1}{ 0} &E_7\\

9_f&\E{1}{ 1}{ 1}{ 0}{ 0}{ 0}{ 0}{ 1}{ 0}&E_6(a_1)  &\bfmu_{9}&9&
\E{-15}{ 1}{ 1}{ 1}{ 1}{ 1}{ 1}{ 0}{ 0} &E_6\\

8_a&\E{0}{ 0}{ 0}{ 0}{ 1}{ 0}{ 0}{ 0}{ 1}&D_8(a_3)  &G_9&8,24&
\E{-21}{ 1}{ 1}{ 1}{ 1}{ 1}{ 1}{ 1}{ 1} &E_8\\

8_b&\E{1}{ 1}{ 0}{ 0}{ 0}{ 1}{ 0}{ 0}{ 0}&D_5  &\bfmu_{8}&8&
\E{-23}{ 1}{ 1}{ 1}{ 1}{ 1}{ 1}{ 1}{ 2} &D_5\\

8_c&\E{1}{ 0}{ 0}{ 1}{ 0}{ 0}{ 0}{ 1}{ 0}&D_5    &\bfmu_{8}&8&
\E{-19}{ 1}{ 1}{ 1}{ 1}{ 1}{ 1}{ 1}{ 0} &D_5\\

8_d&\E{0}{ 0}{ 0}{ 1}{ 0}{ 0}{ 1}{ 0}{ 0}&D_5  &\bfmu_{8}&8&
\E{-22}{ 1}{ 1}{ 1}{ 1}{ 1}{ 0}{ 2}{ 2} &D_5\\

8_e&\E{1}{0}{ 1}{ 0}{ 0}{ 0}{ 1}{ 0}{ 0}&D_5  &\bfmu_{8}&8&
\E{-20}{ 1}{ 1}{ 1}{ 1}{ 1}{ 1}{ 0}{ 2} &D_5\\

8_f&\E{1}{ 0}{ 0}{ 0}{ 0}{ 1}{ 0}{ 0}{ 1}&D_5    &\bfmu_{8}&8&
\E{-31}{ 1}{ 1}{ 1}{ 1}{ 1}{ 1}{ 1}{ 6} &D_5\\

8_g&\E{0}{ 1}{ 1}{ 0}{ 0}{ 0}{ 0}{ 1}{ 0}&D_5    &\bfmu_{8}&8&
\E{-22}{ 1}{ 1}{ 1}{ 1}{ 1}{ 1}{ 0}{ 3} &D_5\\

8_h&\E{0}{ 1}{ 0}{0}{ 0}{ 0}{ 1}{ 0}{ 1}&D_5 &\bfmu_{8}&8&
\E{-12}{ 1}{ 1}{ 1}{ 1}{ 1}{ 0}{ 0}{ 0} &D_5\\

8_i&\E{1}{ 1}{ 1}{ 0}{ 0}{ 0}{ 0}{ 0}{ 1}&D_5  &\bfmu_{8}&8&
\E{-14}{ 0}{ 1}{ 1}{ 1}{ 1}{ 1}{ 0}{ 0} &D_5\\

8_j&\E{0}{ 1}{ 0}{ 1}{ 0}{ 0}{ 0}{ 0}{ 1}&D_5  &\bfmu_{8}&8&
\E{-24}{ 1}{ 1}{ 1}{ 1}{ 1}{ 1}{ 0}{ 4} &D_5\\

8_k&\E{1}{ 1}{ 0}{ 0}{ 0}{ 0}{ 0}{ 1}{ 1}&D_5  &\bfmu_{8}&8&
\E{-16}{ 1}{ 1}{ 1}{ 1}{ 1}{ 1}{ 0}{ 0} &D_5\\

7_a&\E{0}{ 0}{ 0}{ 0}{ 0}{ 1}{ 0}{ 0}{ 1}&A_6=E_7(a_1)^2  &\bfmu_{14}&14&
\E{-22}{ 1}{ 1}{ 1}{ 1}{ 1}{ 1}{ 1}{ 1} &E_7\\

7_b&\E{1}{ 0}{ 0}{ 0}{ 1}{ 0}{ 0}{ 0}{ 0}&A_6=E_7(a_1)^2  &\bfmu_{14}&14&
\E{-26}{ 1}{ 1}{ 1}{ 1}{ 1}{ 1}{ 1}{ 3} &E_7\\

7_c&\E{0}{ 0}{ 0}{ 1}{ 0}{ 0}{ 0}{ 1}{ 0}&A_6=E_7(a_1)^2  &\bfmu_{14}&14&
\E{-24}{ 1}{ 1}{ 1}{ 1}{ 1}{ 1}{ 1}{ 2} &E_7\\

7_d&\E{1}{ 1}{ 0}{ 0}{ 0}{ 0}{ 1}{ 0}{ 0}&A_6=E_7(a_1)^2  &\bfmu_{14}&14&
\E{-28}{ 1}{ 1}{ 1}{ 1}{ 1}{ 1}{ 1}{ 4} &E_7\\

6_a&\E{1}{ 0}{ 0}{ 0}{ 0}{ 1}{ 0}{ 0}{ 0}&E_8(a_8)=-4A_2 &G_{32}&12,18,24,30&
\E{-23}{ 1}{ 1}{ 1}{ 1}{ 1}{ 1}{ 1}{ 1} &E_8\\

6_b&\E{0}{ 0}{ 0}{ 1}{ 0}{ 0}{ 0}{ 0}{ 1}&E_7(a_4)=E_7^3  &G_{26}&6,12,18&
\E{-21}{ 1}{ 1}{ 1}{ 1}{ 1}{ 1}{ 1}{ 0} &E_7\\

6_c&\E{0}{ 1}{ 0}{ 0}{ 0}{ 0}{ 1}{ 0}{ 0}&D_6(a_2)  &G(6,1,2)&6,12&
\E{-21}{ 0}{ 1}{ 1}{ 1}{ 1}{ 1}{ 1}{ 1} &D_7\\

6_d&\E{0}{ 0}{ 1}{ 0}{ 0}{ 0}{ 0}{ 1}{ 0}&E_6(a_2)  &G_5&6,12&
\E{-22}{ 1}{ 1}{ 1}{ 1}{ 1}{ 1}{ 0}{ 2} &E_6\\

6_e&\E{1}{ 1}{ 0}{ 0}{ 0}{ 0}{ 0}{ 1}{ 0}&E_6(a_2)  &G_5&6,12&
\E{-18}{ 1}{ 1}{ 1}{ 1}{ 1}{ 1}{ 0}{ 0} &E_6\\

6_f&\E{0}{ 0}{ 0}{ 0}{ 1}{ 0}{ 0}{ 0}{ 0}&A_5  &\bfmu_{6}&6&
\E{-22}{ 1}{ 0}{ 1}{ 1}{ 1}{ 1}{ 1}{ 2} &A_5\\

6_g&\E{0}{ 0}{ 0}{ 0}{ 0}{ 0}{ 1}{ 0}{ 1}&A_5  &\bfmu_{6}&6&
\E{-18}{ 1}{ 0}{ 1}{ 1}{ 1}{ 1}{ 1}{ 0} &A_5\\

6_h&\E{1}{ 1}{ 1}{ 0}{ 0}{ 0}{ 0}{ 0}{ 0}&D_4  &\bfmu_{6}&6&
\E{-18}{ 1}{ 1}{ 1}{ 1}{ 1}{ 0}{ 0}{ 2} &D_4\\

6_i&\E{1}{ 0}{ 1}{ 0}{ 0}{ 0}{ 0}{ 0}{ 1}&D_4  &\bfmu_{6}&6&
\E{-25}{ 1}{ 1}{ 1}{ 1}{ 1}{ 2}{ 1}{ 0} &D_4\\

6_j&\E{1}{ 0}{ 0}{ 0}{ 0}{ 0}{ 0}{ 1}{ 1}&D_4  &\bfmu_{6}&6&
\E{-12}{ 0}{ 1}{ 1}{ 1}{ 1}{ 0}{ 0}{ 0} &D_4\\

6_k&\E{0}{ 1}{ 0}{ 1}{ 0}{ 0}{ 0}{ 0}{ 0}&D_4  &\bfmu_{6}&6&
\E{-18}{ 0}{ 1}{ 1}{ 1}{ 1}{ 0}{ 0}{ 3} &D_4\\

\end{array}
$$
\end{center}

\begin{center}
{\small Table 21 continued: The gradings of positive rank in type $E_8$}
$$
\begin{array}{cccccccc}
\hline
\text{No.}&\text{Kac diagram}& w& W(\fc,\theta)&\text{degrees}&\theta'&L_\theta\\
\hline
5_a&\E{0}{ 0}{ 0}{ 0}{ 0}{ 1}{ 0}{ 0}{ 0}&2A_4=E_8^6  &G_{16}&20,30&
\E{-24}{ 1}{ 1}{ 1}{ 1}{ 1}{ 1}{ 1}{ 1} &E_8\\

5_b&\E{1}{ 0}{ 0}{ 1}{ 0}{ 0}{ 0}{ 0}{ 0}&A_4  &\bfmu_{10}&10&
\E{-22}{ 1}{ 1}{ 1}{ 1}{ 1}{ 1}{ 1}{ 0} &D_6\\

5_c&\E{0}{ 0}{ 1}{ 0}{ 0}{ 0}{ 0}{ 0}{ 1}&A_4  &\bfmu_{10}&10&
\E{-26}{ 2}{ 1}{ 1}{ 1}{ 1}{ 1}{ 1}{ 1} &D_6\\

5_d&\E{1}{ 0}{ 0}{ 0}{ 0}{ 0}{ 1}{ 0}{ 0}&A_4  &\bfmu_{10}&10&
\E{-30}{ 1}{ 1}{ 1}{ 1}{ 1}{ 1}{ 1}{ 4} &D_6\\

5_e&\E{0}{ 1}{ 0}{ 0}{ 0}{ 0}{ 0}{ 1}{ 0}&A_4  &\bfmu_{10}&10&
\E{-24}{ 0}{ 1}{ 1}{ 1}{ 1}{ 1}{ 1}{ 2} &D_6\\

5_f&\E{1}{ 1}{ 0}{ 0}{ 0}{ 0}{ 0}{ 0}{ 1}&A_4  &\bfmu_{10}&10&
\E{-20}{ 0}{ 1}{ 1}{ 1}{ 1}{ 1}{ 1}{ 0} &D_6\\

4_a&\E{0}{ 0}{ 0}{ 0}{ 0}{ 0}{ 1}{ 0}{ 0}&2D_4(a_1)  &G_{31}&8,12,20,24&
\E{-25}{ 1}{ 1}{ 1}{ 1}{ 1}{ 1}{ 1}{ 1} &E_8\\

4_b&\E{1}{ 0}{ 1}{ 0}{ 0}{ 0}{ 0}{ 0}{ 0}&D_4(a_1)=E_6^3  &G_8&8,12&
\E{-27}{ 1}{ 1}{ 1}{ 1}{ 1}{ 1}{ 1}{ 2} &E_6\\

4_c&\E{0}{ 0}{ 0}{ 1}{ 0}{ 0}{ 0}{ 0}{ 0}&D_4(a_1)=E_6^3 &G_8&8,12&
\E{-24}{ 1}{ 1}{ 1}{ 1}{ 1}{ 1}{ 0}{ 2} &E_6\\

4_d&\E{1}{ 0}{ 0}{ 0}{ 0}{ 0}{ 0}{ 1}{ 0}&D_4(a_1)=E_6^3  &G_8&8,12&
\E{-20}{ 1}{ 1}{ 1}{ 1}{ 1}{ 1}{ 0}{ 0} &E_6\\

4_e&\E{0}{ 1}{ 0}{ 0}{ 0}{ 0}{ 0}{ 0}{ 1}&D_4(a_1)=D_5^2  &G(4,1,2)&4,8&
\E{-16}{ 1}{ 1}{ 1}{ 1}{ 1}{ 0}{ 0}{ 0} &D_5\\

3_a&\E{0}{ 0}{ 1}{ 0}{ 0}{ 0}{ 0}{ 0}{ 0}&4A_2=E_8^{10}  &G_{32}&12,18,24,30&
\E{-26}{ 1}{ 1}{ 1}{ 1}{ 1}{ 1}{ 1}{ 1} &E_8\\

3_b&\E{0}{ 0}{ 0}{ 0}{ 0}{ 0}{ 0}{ 1}{ 0}&3A_2=E_7^6  &G_{26}&6,12,18&
\E{-24}{ 1}{ 1}{ 1}{ 1}{ 1}{ 1}{ 1}{ 0} &E_7\\

3_c&\E{1}{ 1}{ 0}{ 0}{ 0}{ 0}{ 0}{ 0}{ 0}&2A_2 =D_7^4 &G(6,1,2)&6,12&
\E{-24}{ 0}{ 1}{ 1}{ 1}{ 1}{ 1}{ 1}{ 1} &D_7\\

3_d&\E{1}{ 0}{ 0}{ 0}{ 0}{ 0}{ 0}{ 0}{ 1}&A_2=D_4^2  &\bfmu_{6}&6&
\E{-15}{ 0}{ 1}{ 1}{ 1}{ 1}{ 0}{ 0}{ 0} &D_4\\

2_a&\E{0}{ 1}{ 0}{ 0}{ 0}{ 0}{ 0}{ 0}{0}&8A_1=-1  &W(E_8)&2,8,12,14,18,20,24,30&
\E{-27}{ 1}{ 1}{ 1}{ 1}{ 1}{ 1}{ 1}{ 1} &E_8\\

2_b&\E{0}{ 0}{ 0}{ 0}{ 0}{ 0}{ 0}{ 0}{ 1}&4A_1'' =E_6^6 &W(F_4)&2,6,8,12&
\E{-22}{ 1}{ 1}{ 1}{ 1}{ 1}{ 1}{ 0}{ 0} &E_6\\

1_a&\E{1}{0}{ 0}{ 0}{ 0}{ 0}{ 0}{ 0}{ 0} &1 &W(E_8)&2,8,12,14,18,20,24,30&
\E{1}{ 0}{ 0}{ 0}{ 0}{ 0}{ 0}{ 0}{ 0} &E_8\\
\hline

\end{array}
$$
\end{center}

\section{Little Weyl groups for inner type $E$ and Kostant sections}\label{littleweylE}
In this section we compute the little Weyl groups $W(\fc,\theta)$ and their degrees when 
$\theta$ is inner of positive rank in type $E$. As a byproduct we show that every positive rank inner automorphism is principal in a Levi subgroup. This leads to a verification of Popov's conjecture on the existence of Kostant sections, and gives a characterization of the orders of positive-rank automorphisms. 

\subsection{The Levi subgroup $L_\theta$}
In  tables 19-21 above we have indicated a Levi subgroup 
$L_\theta$ whose corresponding subset $J\subset\{1,\dots,\ell\}$ satisfies the conditions of Lemma \ref{Jlowerbound}, giving an embedding 
\begin{equation}\label{embedding}
C_{W_J}(w)\hra W(\fc,\theta).
\end{equation}
In each case, the embedding \eqref{embedding} turns out to be an isomorphism. 
It follows that the degrees of $W(\fc,\theta)$ are
those degrees of $W_J$ which are divisible by $m$. 

We verify that \eqref{embedding} is an isomorphism as follows.  
Let  $U_J\subset W$ be the subgroup acting trivially on the span of the roots $\al_j$ for $j\in J$ and set $c_J(w)=|C_W(w)|/|U_J|$. 
Lemma \ref{simplebound} 
shows that $|W(\fc,\theta)|$ divides $ c_J(w)$.
The subgroup $U_J$ can be found in the tables of 
\cite{carter:weyl} (it is denoted there by $W_2$). 
In all but eight cases we find that 
$$|C_{W_J}(w)|=c_J(w),$$
showing that $C_{W_J}(w)=W(\fc,\theta)$. 

We list the exceptional cases for which $|C_{W_J}(w)|<c_J(w)$.   
We write $|C_{W_J}(w)|$   as the product of degrees divisible by $m$. 
\begin{center}
$$
{\renewcommand{\arraystretch}{1.3}
\begin{array}{ccccccc}
\hline
G& \text{no.} & w& J & |C_{W_J}(w)| & c_J(w)\\
\hline
E_6& 4_b & D_4(a_1) & D_5 & 4\cdot 8 & 8\cdot 12\\
E_7 & 9_b & E_6(a_1)& E_6 & 9 & 18\\
E_7 & 5_d,5_e & A_4 & A_4& 5 & 10\\
E_7 & 4_d,4_e & D_4(a_1) & D_5 & 4\cdot 8 & 8\cdot 12\\
E_8 & 9_f& E_6(a_1) & E_6 & 9 & 18\\
E_8 & 4_e& D_4(a_1) & D_5 & 4\cdot 8 & 8\cdot 12\\
\hline
\end{array}}
$$
\end{center}
To show that $W(\fc,\theta)=C_{W_J}(w)$ in all of these cases as well, 
it suffices to show that $G_0$ has an invariant polynomial of degree $d=4,9,5,4,9,4$ 
for the respective rows. If $k$ has characteristic zero this can be done using the computer algebra system LiE to find the dimension of the $G_0$-invariants in $\Sym^d(\fg_1^\ast)$. In fact we did this for all of the positive rank cases in exceptional groups, as a confirmation of our tables. If $k$ has positive characteristic $p$ (not dividing $m$) the desired invariant is provided by the following result which is apparently standard, but we could not find a reference. 

\begin{lemma}\label{invarianttheory} Let $\rho:H\to \GL(V)$ be a rational representation of a reductive algebraic group $H$ over the ring $\bz[\zeta]$, where $\zeta\in\overline\bq$ is a primitive $m^{th}$-root of unity. Assume that $H(\overline\bq)$ has a nonzero invariant vector in $V(\overline\bq)$ with multiplicity one.  Then $H(k)$ has a nonzero invariant in $V(k)$ for any algebraically-closed field $k$ of characteristic $p$ not dividing $m$.  
\end{lemma}
\proof
Let $W(k)$ be the ring of Witt vectors of $k$, let $K$ be the quotient field of $W(k)$ and let $L$ be an algebraic closure of $K$. Our assumption implies, via complete reducibility,  that $\dim_LV(L)^{H(L)}=\dim_{\overline\bq}V(\overline\bq)^{H(\overline\bq)}=1$. Let  $f\in V(L)$ be a generator of $V(L)^{H(L)}$. 
The line $L\cdot f$ is preserved by $\Gal(L/K)$, 
so Hilbert's theorem 90 implies that $L\cdot f\cap V(K)$ is nonzero.
We may therefore assume that $f\in V(K)$. 
Clearing denominators, we may further assume that $f\in V(W(k))$ and 
is nonzero modulo the maximal ideal $M$ of $W(k)$. The reduction of $f$ modulo $M$ gives a nonzero invariant of $H(k)$ in $V(k)$. 
\qed

As illustrated in the following examples, we can often compute the desired invariant by hand.

\subsubsection{Example: $E_6$ no. $4_b$}
We label the affine diagram  of $E_6$ and write the Kac diagram of $\theta$ respectively as as shown:
$$\EVI{0}{1}{6}{2}{3}{4}{5}\qquad\qquad  \EVI{1}{0}{0}{1}{0}{0}{1}.$$
We view $\fg_1$ as a representation of $\SL_2\times\SL_4\times T_2$, 
where $T_2$ is the two dimensional torus whose cocharacter group has basis 
$\{\check\om_2,\check\om_5\}$, where $\check\om_i$ are the fundamental co-weights of $E_6$. Each node $i$ labelled $1$ in the Kac diagram gives a summand $V_i$ of $\fg_1^\ast$
whose highest weight is the fundamental weight on each node adjacent to $i$ and 
with central character $\al_i$ restricted to $T_2$. Thus, we have 
$$
\begin{array}{cccccc}
\fg_1^\ast\simeq 
&(\mathbf{2}\boxtimes\mathbf{6})
&\oplus&(\mathbf{1}\boxtimes\check{\mathbf{4}})
&\oplus&(\mathbf{1}\boxtimes\mathbf{4})\\
\check\om_2=
&1&&-2&&0\\
\check\om_5=
&0&&-1&&1
\end{array}
$$
Here $\mathbf{2}$ and $\mathbf{4}$ are the standard representations of $\SL_2$ and $\SL_4$, 
$\check {\mathbf{4}}$ is the dual of $\mathbf{4}$ and 
$\mathbf{6}=\Lambda^2\mathbf{4}$. 
It follows that the symmetric algebra of $\fg_1^\ast$ can have 
$G_0$-invariants only in tri-degrees $(2k,k,k)$. To find the expected invariant of degree four, 
we must find an $\SL_2\times\SL_4$-invariant in the summand for $k=1$: 
$$\Sym^2(\mathbf{2}\boxtimes\mathbf{6})
\otimes (\mathbf{1}\boxtimes\check{\mathbf{4}})
\otimes(\mathbf{1}\boxtimes\mathbf{4})=
\Sym^2(\mathbf{2}\boxtimes\mathbf{6})
\otimes (\mathbf{1}\boxtimes\End(\mathbf {4})).
$$
Since $m=4$ we have $p\neq 2$, so $\End(\mathbf {4})=\mathbf{1}\oplus\fsl_4$.
Since $\mathbf{2}$ and $\mathbf{6}$ are self-dual, 
affording alternating and symmetric forms, respectively, 
our invariant must be given by an $\SL_2\times\SL_4$-equivariant mapping 
$\Sym^2(\mathbf{2}\boxtimes\mathbf{6})\to \mathbf{1}\otimes\fsl_4.$
Indeed, wedging in both factors gives a map
$$\Sym^2(\mathbf{2}\boxtimes\mathbf{6})\lra 
\Lambda^2\mathbf{2}\boxtimes\Lambda^2\mathbf{6}
=\mathbf{1}\boxtimes\fso_6\simeq\mathbf{1}\boxtimes\fsl_4,
$$
exhibiting the desired invariant of degree four. 

\subsubsection{Example: $E_7$ no. $5_d$}
The Kac diagram is 
$$\EVII{0}{ 1}{ 1}{ 0}{ 0}{ 0}{ 0}{ 1}$$
with $G_0^{sc}=\SL_2\times\SL_5$ and 
$$\fg_1^\ast=
(\mathbf{2}\boxtimes\mathbf{5})\oplus
(\mathbf{1}\boxtimes\check{\mathbf{5}})\oplus
(\mathbf{1}\boxtimes\Lambda^2\mathbf{5}).
$$
The center of $G_0$ has invariants in tri-degrees $(2k,k,2k)$, leading us to seek an 
$\SL_5$-equivariant mapping 
$$\mathbf{5}\otimes\Sym^{2}(\Lambda^2\check{\mathbf{5}})\lra 
\Sym^2(\mathbf{2}\boxtimes\mathbf{5})^{\SL_2}.
$$
Let $U$ and $V$ be $k$-vector spaces of dimensions $2$ and arbitrary $n<\infty$, respectively. 
Let $P_2(\Hom(V,U))$ be the space of degree two-polynomials on $\Hom(V,U)$, with the natural $SL(V)\times\SL(U)$-action. Then we have a nonzero (hence injective) mapping 
$$
\vp:\Lambda^2(V)\lra P_2(\Hom(V,U))^{\SL(U)}, \quad \om\mapsto \vp_\om,
$$
given by $\vp_\om(f)=f_\ast(\om)$, where $f_\ast:\Lam^2(V)\to \Lam^2(U)\simeq k$ 
is the map induced by $f$. One checks that 
$\dim P_2(\Hom(V,U))^{\SL(U)}=\binom{n}{2}$, so that $\vp$ is an isomorphism of $\SL(V)$-modules 
\begin{equation}\label{SLV}
\Lambda^2(V)\simeq P_2(\Hom(V,U))^{\SL(U)}.
\end{equation}
Returning to our task, we now must find an $\SL_5$-equivariant mapping 
$$\mathbf{5}\otimes\Sym^{2}(\Lambda^2\check{\mathbf{5}})\lra\Lam^2\mathbf{5}.$$
The contraction mapping  
$$\mathbf{5}\otimes\Lambda^2\check{\mathbf{5}}\lra \check{\mathbf{5}},\quad v\otimes\om\mapsto c_v(\om),
$$
where $c_v(\lam\wedge\mu)=\la \lam,v\ra\mu-\la\mu,v)\lam$, extends to a mapping
$$
\mathbf{5}\otimes\Sym^{2}(\Lambda^2\check{\mathbf{5}})
\lra\Lam^3\check{\mathbf{5}}, \qquad
v\otimes(\om\cdot\eta)\mapsto c_v(\om)\wedge\eta+c_v(\eta)\wedge\om.
$$
Since $\Lam^3\check{\mathbf{5}}\simeq\Lam^2\mathbf{5}$ as $\SL_5$-modules, 
we have the desired invariant.

\subsubsection{Example: $E_7$ no. $4_d$} The Kac diagram is 
$$\EVII{1}{0}{1}{0}{0}{0}{0}{1}$$
and $G_0^{sc}=\SL_6$ with 
$\fg_1=\mathbf{6}\oplus\check{\mathbf{6}}\oplus\Lam^3\mathbf 6.$
The action of the center leads us to seek an $\SL_6$-invariant in 
$$\mathbf{6}\otimes\check{\mathbf{6}}\otimes\Sym^2(\Lam^3\mathbf 6).$$
If $V$ is a $k$-vector space of even dimension $2m$, we have a nonzero $\SL(V)$-equivariant mapping 
$$\vp:\End(V)\lra P_2(\Lam^m V),\qquad A\mapsto \vp_A,$$
given by $\vp_A(\om)=\om\wedge A_\ast\om$. Since the $\SL(V)$-module 
$\Lam^mV$ is self-dual this may be viewed as a nonzero mapping $\End(V)\to \Sym^2(\Lam^mV)$. 
Taking $m=3$ gives the desired invariant.

\subsubsection{Example: $E_7$ no. $4_e$} The Kac diagram is 
$$\EVII{0}{1}{0}{0}{0}{0}{1}{0}$$
with $G_0^{sc}= H_1\times \Spin_8\times H_2$, where $H_1\simeq H_2\simeq\SL_2$, and 
$\fg_1^\ast=
(\mathbf{2}\boxtimes \mathbf{8}\boxtimes\mathbf{1})\oplus
(\mathbf{1}\boxtimes \mathbf{8'}\boxtimes\mathbf{2}),
$
where $\mathbf{8}$ and $\mathbf{8'}$ are non-isomorphic eight dimensional irreducible representations of $\Spin_8$. The action of the center leads us to seek an invariant in 
$$\Sym^2(\mathbf{2}\boxtimes \mathbf{8}\boxtimes\mathbf{1})\otimes
\Sym^2(\mathbf{1}\boxtimes \mathbf{8'}\boxtimes\mathbf{2}).
$$
Since every representation in sight is self-dual, we require a $\Spin_8$-equivariant mapping 
from the $H_1$-coinvariants to the $H_2$-invariants:
$$\Sym^2(\mathbf{2}\boxtimes \mathbf{8}\boxtimes 1)_{H_1}\lra
\Sym^2(\mathbf{1}\boxtimes \mathbf{8'}\boxtimes\mathbf{2})^{H_2}.
$$

Since $m=4$, the characteristic of $k$ is not two, so for a $k$-vector space $V$ of arbitrary finite dimension $n$ the symmetric square 
$$\Sym^2(V\otimes\mathbf{2})=
\Sym^2(\mathbf{2}^{\oplus n})=
n\cdot\Sym^2(\mathbf{2})\oplus\binom{n}{2}(\mathbf{2}\otimes\mathbf{2})
$$ is completely reducible as an $\SL_2$-module. Hence the canonical map 

\begin{equation}\label{coinvar}
\Sym^2(V\otimes\mathbf{2})^{\SL_2}\lra 
\Sym^2(V\otimes\mathbf{2})_{\SL_2},
\end{equation}
from the invariants to the coinvariants, is an isomorphism of $\SL(V)$-modules. 
From \eqref{SLV}, both modules are isomorphic to $\Lam^2V$. 

Returning to our task, we now require a $\Spin_8$-equivariant mapping
$$\Lam^2 \mathbf{8}\to \Lam^2\mathbf{8'}.$$
But both of these exterior squares are isomorphic to the adjoint representation of $\Spin_8$, 
whence the desired invariant.

\subsubsection{Example:  $E_8$ no. $4_e$} The Kac diagram is 
$$\E{0}{ 1}{ 0}{ 0}{ 0}{ 0}{ 0}{ 0}{ 1}$$
with $G_0^{sc}=\Spin(12)\times \SL_2$ and
$$\fg_1^\ast=(\mathbf{32}\boxtimes\mathbf{1})\oplus (\mathbf{12}\boxtimes\mathbf{2}),$$
where $\mathbf{32}$ is one of the half-spin representations of $\Spin_{12}$, which is symplectic since $6\equiv 2\mod 4$.  The action of the center of $G_0$ leads us to seek an invariant in 
$$\Sym^2(\mathbf{32}\boxtimes\mathbf{1})\otimes\Sym^2(\mathbf{12}\boxtimes\mathbf{2}).
$$
We must therefore find a $\Spin_{12}\times\SL_2$-equivariant  mapping 
$$\Sym^2(\mathbf{12}\boxtimes\mathbf{2})\lra \Sym^2(\mathbf{32}\boxtimes\mathbf{1}).$$

From \eqref{SLV} and \eqref{coinvar} this is equivalent to a $\Spin_{12}$-equivariant mapping 
$$\Lam^2\mathbf{12}\lra \Sym^2(\mathbf{32}).$$
But $\Lam^2\mathbf{12}\simeq \fso_{12}$ and $\Sym^2(\mathbf{32})\simeq \fsp_{32}$. 
The desired mapping $\fso_{12}\to\fsp_{32}$ is simply the representation of 
$\fso_{12}$ on the symplectic half-spin representation $\mathbf{32}$. 

\subsection{A remark on saturation} Let 
$$W^*({\mathfrak c},\theta):=N_{G^\theta}({\mathfrak c})/C_{G^\theta}({\mathfrak c}).$$
Clearly $W({\mathfrak c},\theta)\subset W^*({\mathfrak c},\theta)\subset W_1^\theta$ 
(see \eqref{W1}). 
We say that $\theta$ is {\bf saturated}  if $W({\mathfrak c},\theta)=W^*({\mathfrak c},\theta)$. (For the adjoint group $G$ this is equivalent to the definition given in 
section 5 of \cite{vinberg:graded}.) Clearly $\theta$ is saturated if $G^\vt=G_0$. 
As remarked in section \ref{stableisotropy} this holds whenever the group $\Om_\vt(x)$ is trivial. 
In particular, saturation holds in types $G_2, {^3D_4}, F_4, E_8, {^2E_6}$, where $\Om_\vt$ itself is trivial. 
It is known (\cite{vinberg:graded}, \cite{levy:thetap}) that all gradings on classical Lie algebras are saturated except for certain outer automorphisms of order divisible by 4 in type $D_n$. 
It remains to consider only those inner automorphisms of $E_6$ and $E_7$ where the Kac diagram is invariant under the symmetries of the affine Dynkin diagram and we have $W(\theta,\fc)\neq W^\theta_1$. The latter implies that $|C_{W_J}(w)|<c_J(w)$. The only cases not thus eliminated are $4_d$ and $4_e$ in type $E_7$. But in these two cases we have $c_J(w)/|W(\theta,\fc)|=3$, while $[G^\theta:G_0]=2$, so saturation holds in these cases as well. We conclude that  all gradings on exceptional Lie algebras are saturated.

\subsection{Kostant sections and the Levi subgroup $L_\theta$}
 
 A {\bf Kostant section} 
\footnote{In the literature, this is also called a "Kostant-Weierstrass" or "KW" section because in the case of  the non-pinned outer triality automorphism of $\fso_8$ such a section is equivalent to the Weierstrass-normal form of a nonsingular homogeneous cubic polynomial in three variables.} 
for the grading $\fg=\oplus_{i\in\bz/m}\ \fg_i$  is an affine subspace $\fv\subset \fg_1$ such that the embedding $\fv\hra\fg_1$ induces an isomorphism of affine varieties 
 $\fv\overset\sim\lra \fg_1//G_0$, or equivalently, if the restriction map 
 $k[\fg_1]^{G_0}\lra k[\fv]$ is bijective. 
 
 Recall that we have fixed a pinning $(X,R,\check X,\check R,\{E_i\})$ in $G$, which determines the co-character $\check\rho\in X_\ast(T)$ and principal nilpotent element $E=\sum E_i$, such that 
 $\check\rho(t)\cdot E=tE$. 
From \cite[Thm.3.5]{panyushev:theta} and 
 \cite[Prop.5.2]{levy:thetap} we have the following existence result for Kostant sections. 
 \begin{thm}\label{principalkostant} \label{kostant}  Assume the characteristic of $k$ is not a torsion prime for $G$,  that $m$ nonzero  in $k$. Then the grading $\fg=\oplus_{i\in \bz/m}\ \fg_i$ associated to the principal automorphism 
$\theta_m=\check\rho(\zeta)\vt$ has a Kostant section $E+\fu$, where $\fu$ is any vector space complement to $[\fg_0,E]$ in $\fg_1$ such that $\fu$ is stable under 
$\check\rho(k^\times)$.
\end{thm}
 
We have seen that for each positive-rank torsion inner automorphism in type $E_{6,7,8}$ there exists a subset 
 $J\subseteq\{1,2,\dots,\ell\}$ such that $W(\fc,\theta)=C_{W_J}(w)$. 
This can also be checked for the classical groups and types $F_4, G_2$. Thus, we have a case-by-case proof of the following theorem.

\begin{thm} Let $\theta$ be an inner automorphism of $\fg$ whose order $m$ is nonzero in $k$ and let $\fc$ be a Cartan subspace of $\fg_1$. Then there exists a $\theta$-stable Levi subgroup $L=L_\theta$ whose  Lie algebra $\fl$ contains $\fc$ in its derived subalgebra, such that the following hold:
\begin{enumerate}
\item $\theta\vert_\fl=\Ad(\check\rho_L(\zeta))$.
\item The inclusion of little  Weyl groups $W_L(\fc,\theta)\hra W(\fc,\theta)$ is a bijection. 
In particular, the degrees of $W(\fc,\theta)$ are precisely the degrees of 
$W_L$ which are divisible by $m$. 
\item The restriction map 
$k[\fg_1]^{G_0}\lra k[\fl_1]^{L_0}$
is a bijection.
\end{enumerate}
\end{thm}

In view of Thm. \ref{kostant}, we conclude:
\begin{cor}\label{cor:kostantE}
Every positive-rank torsion inner automorphism in type $E_{6,7,8}$ has a Kostant section contained in the Levi subalgebra $\fl$ of the previous theorem. 
\end{cor}

We also observe:
\begin{cor}\label{cor:m}
A positive integer $m$ is the order of a torsion inner automorphism of positive rank precisely if $m$ is the order of a $\bz$-regular element in the Weyl group of a Levi subgroup of $G$.
\end{cor}

\section{Outer gradings of positive rank in type $E_6$}\label{2E6}
We realize the outer pinned automorphism of $E_6$ as the restriction of an affine pinned automorphism of $E_7$, as in section \ref{affine-pinned}. 

\subsection{Root systems of type $E_7$ and ${^2E_6}$}\label{E7to2E6}
Let $(Y,R,\check  Y, \check  R)$ be a root datum of adjoint type $E_7$
and fix a base $\De=\{\al_1,\dots,\al_7\}\subset R$ with lowest root $\al_0$, 
according to the numbering 
\begin{equation}\label{E7numbering}
\EVII{0}{1}{4}{2}{3}{5}{6}{7}.
\end{equation}
The set $\Pi:=\{\al_0\}\cup\De$ has stabilizer $W_\Pi=\{1,\vt\}$ of order two, where 
$\vt=r_1r_2r_3$ is a product of reflections about mutually orthogonal roots 
$\ga_1,\ga_2,\ga_3$ in which the coefficients of simple roots $\{\al_1,\dots,\al_7\}$ are given by 
$$
\ga_1=\EVII{}{0}{1}{1}{2}{2}{2}{1},\qquad 
\ga_2=\EVII{}{1}{1}{1}{2}{2}{1}{1},\qquad
\ga_3=\EVII{}{1}{1}{2}{2}{1}{1}{1}.
$$
The sum
$$\check \ga_1+\check \ga_2+\check \ga_3=2\check \mu,$$
where $\check \mu=\check \om_7$ is the nontrivial minuscule co-weight. 

Regard the vector space $V=\br\otimes\check  Y$ as an affine space with $0$ as basepoint. Each linear functional $\lam:V\to \br$ 
is then regarded as an affine function on $V$ vanishing at $0$, we have the affine root system 
$$\Phi=\{\al+n:\ \al\in R,\ n\in\bz\}$$
with basis $\{\phi_0,\phi_1,\dots,\phi_7\}$ where
$\phi_0=1+ \al_0,\phi_1= \al_1,\dots,\phi_7= \al_7$ satisfy the relation 
\begin{equation}\label{E7relation}
\phi_0+2\phi_1+3\phi_2+4\phi_3+2\phi_4+3\phi_5+2\phi_6+\phi_7=1.
\end{equation}
A point $x\in V_\bq$ of order $m$ has Kac diagram 
\begin{equation}\label{E7coord}
\eVII{s_0}{s_1}{s_2}{s_3}{s_5}{s_6}{s_7}{s_4},
\end{equation}
where $s_i/m=\phi_i(x)$. 

The affine transformation $\widetilde\vt:V\to V$ given by 
$$\widetilde\vt(x)= \check \mu+\vt\cdot x$$
permutes the simple affine roots $\{\phi_0,\dots,\phi_7\}$ according to the nontrivial symmetry of the affine diagram of $E_7$.  
The fixed-point space of $\widetilde\vt$ in $V$ is given by
$$\sca^{\vt}:=V^\vt+\tfrac{1}{2}\check \mu,$$
which is an affine space under the vector space $V^\vt=\br\otimes \check  Y_\vt$, 
with basepoint $\tfrac{1}{2}\check\mu$. 
The rational points in $\sca^{\vt}$ are precisely those points $x\in V_\bq$  whose Kac diagram  has the symmetric form
\begin{equation}\label{si}
\eVII{s_0}{s_1}{s_2}{s_3}{s_2}{s_1}{s_0}{s_4},
\end{equation}
in which case equation \eqref{E7relation} implies that 
\begin{equation}\label{s0}
s_0+2s_1+3s_2+2s_3+s_4=m/2,
\end{equation}
where $m$ is the order of $x$. 

The automorphism $\vt$ permutes the roots $\al_1,\dots,\al_6$ which generate a root subsystem $R'$  of type $E_6$. The co-weight lattice $\check  X=\Hom(\bz R',\bz)$ has dual basis $\{\check \om_1,\dots,\check \om_6\}$ and we have 
$$\check X^\vt=\check  Y^\vt.$$
Hence $\sca^{\vt}$ is also an affine space under $\br\otimes \check X^\vt$ 
and we may construct the affine root system $\Psi(R',\vt)$ as in section \ref{affine}, 
using the point $x_0=\tfrac{1}{2}\check \mu$. We have $\ell_\vt=4$ and 
$\Psi(R',\vt)$ has basis $\psi_0,\dots,\psi_4$, where
$\psi_i=\al_i\vert_{\sca^{\vt}}$ for $1\leq i\leq 4$ and
\begin{equation}\label{2E6relation}
\psi_0+2\psi_1+3\psi_2+2\psi_3+\psi_4=1/2.
\end{equation}

A rational point $x\in \sca^{\vt}_\bq$ with $E_7$ Kac-diagram \eqref{si}
has ${^2E_6}$ Kac-diagram 
$$s_0\ s_1\ s_2\Leftarrow s_3\  s_4.$$
This is clear for $s_1,\dots, s_4$ since $\psi_i$ is the restriction of $\phi_i$,  and follows for $s_0$ by comparing the relations \eqref{E7relation} and \eqref{2E6relation}.

\subsection{Lie algebras of type $E_7$ and ${^2E_6}$}

Let $k$ be an algebraically closed field of characteristic $\neq 2,3$
and let $\fg$ be a simple Lie algebra over $k$ of type $E_7$ with automorphism group 
$G=\Aut(\fg)$. 
We fix a maximal torus $T\subset G$ with Lie algebra $\ft$ and we choose an affine pinning $\widetilde\Pi=\{E_0,\dots,E_7\}$ for $T$ in $\fg$, numbered as in 
\eqref{E7numbering}. 
As above we let $\vt=r_1r_2r_3\in W_\Pi$ be the unique involution acting on  $\Pi$ via the permutation $(0 7)(1 6)(2 5)$. Recall from section \ref{affine-pinned} that $\vt$ has a lift  $n\in N$ of order two defined via the homomorphism $\vp:\SL_2\to G$ as in  equation \eqref{npinned}. 

Let $S=(T^\vt)^\circ$ be the identity component of the group of fixed-points of $\vt$ in $T$. The co-weight group of $S$ is $\check X^\vt$ and we have 
$$T^\vt=S\times\la \check \mu(-1)\ra,$$
where $\check \mu=\check \om_7$ is the nontrivial minuscule co-weight. 
The automorphism $\vep:=\Ad(\check \mu(-1))$ has order two; its fixed-point subalgebra $\fg^{\vep}$ decomposes as 
$$\fg^{\vep}=\fh\oplus \fz$$
where $\fz=d\check \mu(k)$ and $\fh=[\fg^\vep,\fg^\vep]$, the  derived subalgebra of $\fg^{\vep}$, 
has type $E_6$ and is generated by the root spaces $\fg_\al$ for  
$\al\in\pm\{\al_1,\dots,\al_6\}$.  Note that $\vep$ and $n$ both lie in the subgroup $\vp(\SL_2)$ and are conjugate therein.

The centralizer $C_G(\vep)$ is the normalizer in $G$ of $\fh$, surjecting onto 
$\Aut(\fh)$, and is also the normalizer of in $G$ of $\fz$. 
The centralizer $C_G(\fz)$ of $\fz$ is  the identity component of $C_G(\vep)$, 
and the image of $C_G(\fz)$ in $\Aut(\fh)$ is the group 
$$H:=\Aut(\fh)^\circ$$ 
of inner automorphisms of $\fh$. 
It follows that 
we have an exact sequence
\begin{equation}\label{Hexact}
1\lra \mu(k^\times)\lra C_G(\fz)\lra H\lra 1.
\end{equation}

\begin{prop}\label{symmetric} Let $\theta\in\Aut(\fg)$ be a torsion automorphism whose order $m$ is nonzero in $k$. Then the centralizer $G^\theta$ has at most two components, and the following are equivalent.
\begin{enumerate}
\item The normalized Kac diagram of $\theta$ has the symmetric form 
$\ \EVII{a}{b}{e}{c}{d}{c}{b}{a}.$
\item The $G$-conjugacy class of $\theta$ meets $Sn$. 
\item The centralizer $G^\theta$ has two components and $n$ lies in the non-identity component.
\end{enumerate}
\end{prop}

\proof After conjugating by $G$, we may assume $\theta=\Ad(t)$, where $t=\check \lam(\zeta)$, for some $\check \lam\in \check  X$ and $\zeta\in k^\times$ of order $m$. 
We set $x=\frac{1}{m}\check \lam$.

Over $\bc$, the equivalence $1\Leftrightarrow 3$ follows from \cite[Prop. 2.1]{reeder:torsion}, whose proof, once we replace $\exp(x)$ by $\check\lam(\zeta)$, is also valid over $k$.  

We prove $1\Leftrightarrow 2$. 
From the previous section the Kac coordinates of $\theta$ are symmetric precisely if 
$$x=\check \mu+\vt\cdot x.$$
This is equivalent to having $\check \lam-\frac{m}{2}\check \mu\in \check X^\vt$. 
Evaluating at $\zeta$ this is in turn equivalent to having $t\vep\in S$,  or $t\in S\vep$. 
Since $n$ and $\vep$ are conjugate in $\vp(\SL_2)$ which centralizes $S$ 
(see  Lemma \ref{nfixed}) we can replace $\vep$ by $n$.  
\qed

\begin{prop}\label{E7vsE6} Let $s\in S$ and suppose $sn$ has order $m$ invertible in $k$ and let 
$\theta=\Ad(sn)$ have symmetric normalized Kac diagram 
$\EVII{a}{b}{e}{c}{d}{c}{b}{a}.$
Then 
\begin{enumerate}
\item
$\theta$ normalizes $\fh$ and $\theta\vert_\fh$ is an outer automorphism of 
$\fh$ with Kac diagram 
$$a\ b\ c\Leftarrow d\ e.$$
\item
Every torsion outer automorphism of $\fh$ is conjugate to $\theta\vert_\fh$, 
where $\theta=\Ad(sn)$ for some $s\in S$. 
\item We have $\rank(\theta\vert_\fh)\leq\rank(\theta)$. 
\end{enumerate}
\end{prop}
\proof
Since $\Ad(n)=\vt$ normalizes $\fh$, acting there via a pinned automorphism, 
and $s\in S\subset H$, we have that $\theta\vert_\fh$ is an outer automorphism of $\fh$. The relation between the Kac diagrams of $\theta$ and $\theta\vert_\fh$ follows from the discussion in section \ref{E7to2E6}. 

Assertion 2 is now clear, since every Kac diagram $s_0\ s_1\ s_2\Leftarrow s_3\ s_4$ corresponds to $\Ad(sn)\vert_\fh$ for some $s\in S$. We can also prove assertion 2 directly, as follows: Since $\vt$ preserves the maximal torus $T\cap H$ of $H$, and permutes the simple roots $\{\al_1,\dots,\al_6\}$, every torsion outer automorphism of $\fh$ is $H$-conjugate to one of the form $\Ad(s)\vt$ for some $s\in (T\cap H)^\vt$ (see \cite[Lemma 3.2]{reeder:torsion}, whose proof is valid for $k$). We must therefore show that $(T\cap H)^\vt=S$. Since the Lie algebra of $S$ is $\ft^\vt$ which is contained in $(\ft\cap\fh)^\vt$, it suffices to show that $\ft^\vt\subset\fh$. But $\ft^\vt$ has dimension four and is spanned by $d\check \al_i(1)+d\check \al_{7-i}(1)$ for $1\leq i\leq 4$, and these vectors lie in $\fh$. 

Finally, a Cartan subspace for $\theta\vert_\fh$ is contained in a Cartan subspace for $\theta$, so assertion 3 is obvious. 
\qed

Prop. \ref{E7vsE6}  implies that the Kac diagram of any outer positive rank automorphism of $\fh$ must have the form $a\ b\ c\ \Leftarrow d\ e$, where 
$\EVII{a}{b}{e}{c}{d}{c}{b}{a}$ is a positive rank diagram for $E_7$ appearing in section \ref{Eposrank}. 


For example, there are two outer automorphisms of $\fh$ having order $m=2$, 
namely the restrictions to $\fh$ of $\vt=\Ad(n)$ and $\vt_0=\Ad(n_0)$ where $n_0$ is a lift of $-1\in W(E_7)$. These are the involutions in $E_7$ numbered $2_c$ and $2_a$ respectively Table 20 of section \ref{Eposrank}. The  Kac diagrams in $E_7$ and ${^2E_6}$ are shown:
$$ \begin{array}{c c c c c}
\vt:&\eVII{1}{0}{0}{0}{0}{0}{1}{0}&\qquad&\vt_0:&\eVII{0}{0}{0}{0}{0}{0}{0}{1}\\
&&&&\\
\vt\vert_{\fh}:&1\ 0\ 0\Leftarrow 0\ 0&\qquad&
\vt_0\vert_{\fh}:&0\ 0\ 0\Leftarrow 0\ 1.
\end{array}
$$
Both $\vt$ and $\vt_0$ act by $-1$ on $\fz$. It follows that their ranks in $E_6$ are one less than their ranks in $E_7$, namely
$$\rank(\vt\vert_\fh)=2,\qquad \rank(\vt_0\vert_\fh)=6.$$

\subsection{Positive rank gradings on $E_6$ (outer case)}
From Props.  \ref{symmetric} and \ref{E7vsE6} we know that the Kac diagrams for positive rank gradings in outer type $E_6$ are obtained from symmetric positive-rank diagrams for $E_7$. 
We now adapt our methods for the inner case to complete the classification of positive rank outer gradings of $E_6$. 

 We regard $W(E_6)$ as the subgroup of $W(E_7)$ generated by the reflections for the roots $\al_1,\dots,\al_6$. Equivalently,  $W(E_6)$ is the centralizer of $\fz$ in $W(E_7)$. 
The coset $-W(E_6)=\{w\vt_0:\ w\in W(E_6)\}$ consists of the elements in $W(E_7)$ acting by $-1$ on $\fz$ and contains both $\vt$ and $\vt_0$. 

 \begin{lemma}\label{nw}
  Let $n_w\in N_G(\ft)$ be a lift of an element $w\in -W(E_6)$. 
 Then $\Ad(n_w)$ normalizes $\fh$ and acts on $\fh$ as an outer automorphism. 
 \end{lemma}
 \proof Since $w$ permutes the root spaces in $\fh$ 
 it follows that $n_w$ normalizes $\fh$. 
 Let $n\in N_G(\ft)$ be the lift of $\vt$ constructed above. 
 Both $n$ and  $ n_w$ act by $-1$ on $\fz$, so $n\cdot n_w$ lies in the connected subgroup $C_G(\fz)$  and the image of $n\cdot n_w$ in $\Aut(\fh)$ lies in the subgroup $\Aut(\fh)^\circ$ of inner automorphisms. Since $\Ad(n)=\vt$ is outer on $\fh$, it follows that $\Ad(n_w)$ is outer on $\fh$ as well. 
 \qed

Let  $w\in W(E_7)$ be any element whose order $m$ is invertible in $k$ and such that $w$ has an eigenvalue $\zeta$ of order $m$ on $\ft$. 
Recall that $\Kac(w)$ is the set of normalized Kac diagrams of torsion automorphisms 
$\theta\in \Aut(\fg)$ of order $m$ such that $\theta$ normalizes $\ft$ and acts on $\ft$ 
via $w$. 

Let $\tau\in\Aut(\fh)$ be a torsion outer automorphism with Kac coordinates 
$a\ b\ c\Leftarrow d\ e$.  We write 
$$\tau\leadsto w$$ to mean that the symmetric Kac diagram 
 $\eVII{a}{b}{c}{d}{c}{b}{a}{e}$ appears in $\Kac(w)$. Let $\Kac(w)_{\sym}$ denote the set of symmetric diagrams in $\Kac(w)$.

\begin{prop}\label{prop:2E6} Let $\tau\in\Aut(\fh)$ be a torsion outer 
automorphism whose order $m>2$ is invertible in $k$. Assume that $\rank(\tau)>0$. 
Then there exists $w\in -W(E_6)$ such $\tau\leadsto w$. Moreover, we have 
$$\rank(\tau)=\max\{\rank(w):\ w\in -W(E_6),\ \tau\leadsto w\}.$$
\end{prop}
\proof Let 
$\fc\subset\fh(\tau,\zeta)$ be a Cartan subspace. 
Then $\fc$ is contained in a $\tau$-stable Cartan subalgebra 
$\ft'$ of $\fh$ so that $\fc=\ft'(\tau,\zeta)$. Conjugating by $H$, 
we may assume that $\ft'\subset\ft$ and therefore $\ft=\ft'\oplus\fz$. 

We have $\tau=\theta\vert_\fh$ for some $\theta\in \Aut(\fg)$ normalizing $\fh$. 
Then $\theta$ also normalizes the centralizer $\fz$ of $\fh$. 
Since $\theta\vert_\fh$ is outer but $\theta^2\vert_\fh$ is inner, it follows that $\theta$ acts by $-1$ on $\fz$. 

Since $\theta$ normalizes $\ft$, it projects to an element $w\in W(E_7)$.
The subgroup of $W(E_7)$ normalizing $\fz$ is $\{\pm 1\}\times W(E_6)$ 
and $W(E_6)$ is the subgroup centralizing $\fz$. It follows that $w\in -W(E_6)$. 

Since the normalized Kac diagram of $\theta$ belongs to $\Kac(w)$ 
and $\tau=\theta\vert_\fh$, we have $\tau\leadsto w$. We also have 
$$
\rank(w)=\ft(w,\zeta)=\ft'(w,\zeta).
$$

Suppose now that $w\in -W(E_6)$ is any element for which $\tau\leadsto w$. 
Let $a\ b\ c\Leftarrow d\ e$ be the normalized Kac coordinates for $\tau$. 
Since $\tau\leadsto w$ there is a lift $n_w\in N_G(\ft)$ such that $\Ad(n_w)$ has normalized Kac diagram   $\eVII{a}{b}{c}{d}{c}{b}{a}{e}$.

By Lemma \eqref{nw}, we have that $\Ad(n_w)$ is an outer automorphism of $\fh$. 
Hence there is $s\in S$ such that 
$\Ad(n_w)\vert_\fh$ is $H$-conjugate to $\Ad(sn)\vert_\fh$. 
From the exact sequence \eqref{Hexact}
there are $g\in C_G(\fz)$ and $z\in Z$ such that 
$$gn_wzg^{-1}=sn.$$
But $n_w$ is $Z$-conjugate to $n_wz$, since $w=-1$ on $\fz$. 
Therefore $n_w$ and $sn$ are conjugate under $C_G(\fz)$,
so $\Ad(sn)$ also has normalized Kac diagram $\eVII{a}{b}{c}{d}{c}{b}{a}{e}$. 

By Prop. \ref{E7vsE6}, $\Ad(sn)\vert_\fh$ has Kac diagram $a\ b\ c\Leftarrow d\ e$, 
and therefore $\Ad(sn)\vert_\fh$ is $H$-conjugate to $\tau$. 
But 
$$\Ad(sn)\vert_\fh=\Ad(gn_wzg^{-1})\vert_\fh=\Ad(gn_wg^{-1})\vert_\fh$$
is conjugate to $\Ad(n_w)\vert_\fh$, via the element $h=\Ad(g)\vert_\fh\in H$. 
Thus,  $\tau$ and $\Ad(n_w)\vert_\fh$ are $H$-conjugate. Since $\ft(w,\zeta)\subset \fh$, an $H$-conjugate of $\ft(w,\zeta)$ is contained in a Cartan subspace of $\tau$, so $\rank(w)\leq \rank(\tau)$. This completes the proof.
\qed

The Kac diagrams of positive rank for ${^2E_6}$ are obtained from symmetric positive rank diagrams for $E_7$, of which there are $20$ (see Table 20). 

Three of these ($14_a, 8_d, 8_e$) have rank zero for ${^2E_6}$ as will be explained. Two more have order $m=2$ and are easily handled by known results. The ranks for the remaining $15$
are found as follows. Using Prop. \ref{prop:2E6}, it is enough to extract the symmetric diagrams 
from the preliminary table for $E_7$ in section \ref{preliminary}. The results are shown below, where $r$ is the rank of $\tau$ in ${^2E_6}$.  

\begin{center}
{\small  Table 22:  $\Kac(w)_{\sym}$ for  certain $w$ in $-W(E_6)$ }
$$
\begin{array}{|c|c|c|c|l|l|}
\hline
m& w\in-W(E_6)&w\in W(E_7) & r&\Kac(w)_{\un}& \Kac(w)_{\sym}\\
\hline\hline
18& -E_6(a_1)& E_7 & 1& \EVII{1}{1}{1}{1}{1}{1}{1}{1}& \EVII{1}{1}{1}{1}{1}{1}{1}{1}\\
\hline
12& -E_6&E_7(a_2) & 1& \EVII{1}{ 1}{ 1}{ 0}{ 1}{ 0}{ 1}{ 1}& 
\EVII{1}{ 1}{ 1}{ 0}{ 1}{ 0}{ 1}{ 1}\\
\hline
10& -(A_4+A_1)&D_6& 1& \EVII{\ast}{\ast}{1}{1}{1}{1}{1}{1}
&\EVII{0}{ 1}{ 1}{ 0}{ 1}{ 0}{ 1}{ 0}\qquad\EVII {1}{ 0}{ 1}{ 1}{ 0}{ 1}{ 0}{ 1}
\qquad\EVII {1}{ 1}{ 0}{ 0}{ 1}{ 0}{ 1}{ 1}\\
\hline
8& -D_5&D_5+A_1& 1& \EVII{1}{\ast}{1}{1}{1}{1}{1}{\ast}
&\EVII{0}{ 1}{ 0}{ 0}{ 1}{ 0}{ 1}{ 0}\qquad\EVII{1}{ 0}{ 1}{ 0}{ 1}{ 0}{ 0}{ 1}\\
\hline
6& -(3A_2)&E_7(a_4) & 3& \EVII{1}{0}{0}{0}{1}{0}{0}{1}& 
\EVII{1}{0}{0}{0}{1}{0}{0}{1}\\
\hline
6& -(2A_2)&A_1+D_6(a_2) & 2&  \EVII{1}{\ast}{1}{1}{0}{1}{0}{1}& 
\EVII{0}{ 1}{ 1}{ 0}{ 0}{ 0}{ 1}{ 0}\qquad\EVII{1}{ 0}{ 0}{ 0}{ 1}{ 0}{ 0}{ 1}\\
\hline
6& -A_2&3A_1+D_4& 1&\EVII{1}{\ast}{1}{1}{1}{1}{\ast}{1}& 
\EVII{0}{ 0}{ 1}{ 0}{ 1}{ 0}{ 0}{ 0}\\
\hline
6& -(A_1+A_5'')&A_5'& 1&\EVII{\ast}{1}{\ast}{1}{1}{1}{1}{\ast}& 
\EVII{0}{ 0}{ 0}{ 1}{ 0}{ 1}{ 0}{ 0}\qquad\EVII {1}{ 1}{ 0}{ 0}{ 0}{ 0}{ 1}{ 1}
\qquad\EVII{0}{ 1}{ 1}{ 0}{ 0}{ 0}{ 1}{ 0}\\
&&&&&\EVII{1}{ 0}{ 0}{ 0}{ 1}{ 0}{ 0}{ 1}\\
\hline
4& -D_4(a_1)&A_1+2A_3& 2&\EVII{1}{1}{1}{1}{\ast}{1}{1}{1}& 
\EVII{0}{ 0}{ 0}{ 0}{ 1}{ 0}{ 0}{ 0}\\
\hline
4& -A_3+2A_1&(A_1+A_3)''& 1&\EVII{1}{1}{1}{1}{\ast}{1}{1}{1}& 
\EVII {0}{ 0}{ 0}{ 0}{ 1}{ 0}{ 0}{ 0}\qquad \EVII{0}{ 1}{ 0}{ 0}{ 0}{ 0}{ 1}{ 0}\qquad
\EVII{1}{ 0}{ 1}{ 0}{ 0}{ 0}{ 0}{ 1}\\
\hline

\end{array}
$$
\end{center}

Case $14_a$ has rank zero since there are no elements of order $7$ or $14$ in $W(E_6)$. 
Cases $8_d$ and $8_e$ have rank zero since $D_5$ is the only element of order $8$ in $W(E_6)$ and the Kac diagrams for $8_{d,e}$ do not appear in the row for $w=-D_5$ in Table 22 above.

\subsection{Little Weyl groups for ${^2E_6}$}
The little Weyl groups $W_H(\fc,\tau)$ and their degrees are determined as follows. 

{\bf Cases $18_a,\ 12_b,\ 6_a,\ 4_b,\ 2_a$:\ } These cases are stable, hence by Cor. \ref{vinberg:stable}
we have $W_H(\fc,\tau)=W(\ft')^\theta$, where 
$\ft'$ is the unique Cartan subalgebra of $\fh$ containing $\fc$. Then $W(\ft')^\theta$ and its degrees are determined from \cite{springer:regular}. 

\begin{lemma}\label{m/2}
If $\dim \fc=1$ then $W_H(\fc,\tau)\simeq\bfmu_d$ for some integer $d$ divisible by $m/2$.
\end{lemma}
\proof Since $\dim\fc=1$ we have $W_H(\fc,\tau)\simeq\bfmu_d$ for some integer $d$. 
We may assume $\tau=\Ad(n_w)\vert_\fh$, where $n_w\in N_G(\ft)$ has image 
$w\in-W(E_6)$. Then $n_w^2\in H_0$ has eigenvalue $\zeta^2$ on $\fc$, where $\zeta\in k^\times$ has order $m$ equal to the order of $\tau$. 
It follows that so $m/2$ divides $d$. 
\qed

{\bf Cases $10_a, 10_b, 10_c$:\ } In these cases we have $m=10$ and $\dim\fc=1$ so $\bfmu_5\leq W_H(\fc,\tau)$, by Lemma \ref{m/2}. And $W_H(\fc,\tau)\leq W_H(\fc,\tau^2)$. Now $w^2$ has type $A_4$ in $E_6$, and all lifts of this type have little Weyl group $\bfmu_5$, from Table 19.
\footnote {In fact, using Kac diagrams one can check that classes $10_{a,b,c}$ in Table 20 square to classes $5_{a,b,c}$, respectively, in Table 19.}
So And $W_H(\fc,\tau)\leq W_H(\fc,\tau^2)\simeq\bfmu_5$.

{\bf Cases $8_c,\ 8_f$:\ } In these cases we have $m=8$ and $\dim\fc=1$ so $\bfmu_4\leq W_H(\fc,\tau)\leq \bfmu_8$, by Lemma \ref{m/2}.

In case $8_f$ the diagram for $\theta'$ in Table 20 shows that $\tau$ is principle in $\Aut(\fh)$.  Hence
$W_H(\fc,\tau)=N_{W_H}(\fc)/Z_{W_H}(\fc)$, by Prop. \ref{pan}.  The element $w$ has type $-D_5$ and $\fc$ may be chosen to be the $-\zeta$-eigenspace for $y=-w$ in $\ft$. Since $\la y\ra$ acts faithfully on $\bc$, there is a copy of $\bfmu_8$ in $W_H(\fc,\tau)$.

In case $8_c$ we rule out $\bfmu_4$ using invariant theory, as in section \ref{littleweylE}.
A  degree-four invariant in $\fh_1$ would correspond to an element of 
\begin{equation}\label{hom}
\Hom_{M}\left(
\Sym^2(\mathbf{2}\boxtimes \mathbf{2})^L,
\Sym^2(\mathbf{3}\boxtimes\mathbf{2})^R\right), 
\end{equation}
arising from the action of  $L\times M\times R=\SL_2\times\SL_2\times\SL_2$ on 
$$\fh_1\quad\simeq \quad
\mathbf{2}\boxtimes\mathbf{2}\boxtimes\mathbf{1}\quad\oplus\quad
\mathbf{1}\boxtimes\mathbf{3}\boxtimes\mathbf{2}.
$$
 But 
$\Sym^2(\mathbf{2}\boxtimes \mathbf{2})^L$ is the trivial representation of $M$ 
and $\Sym^2(\mathbf{3}\boxtimes\mathbf{2})^R$ is the adjoint representation of $M$, 
which is irreducible since $p>2$. Hence the vector space \eqref{hom} is zero. 
 
 {\bf Case $6_c$:\ } Here the centralizer of $w=-2A_2$ in $W(E_6)$ has order $108$ and contains a subgroup $W(A_2)$ acting trivially on the root subsystem spanned by the $2A_2$. It 
 follows that $|W_H(\tau,\fc)|\leq 18$. Results in the next section show that $W_H(\tau,\fc)$ contains the centralizer of a $[33]$-cycle in the symmetric group $S_6$, which has order $18$. 
 
{\bf Case $6_g$:\ } Here $\dim\fc=1$ and $w^2$ has type $A_2$, of which  all lifts  in $H$ have little Weyl group $\mu_6$. Hence $\mu_3\leq W_H(\fc,\tau)\leq \mu_6$. 
One checks that an $H_0$-invariant in degree $3$ in $\fh_1$ is a quadratic form on $S^2\mathbf{\check 4}$, which must be trivial. Hence $W_H(\fc,\tau)\simeq \mu_6$. 

{\bf Cases $6_i$, $6_k$:\ } These cases have $m=6$ and $\dim\fc=1$ so $\mu_3\leq W_H(\fc,\tau)$, by Lemma \ref{m/2}. We show this is equality by finding an $H_0$-invariant of degree $3$ on $\fh_1$. 

In case $6_i$, $\fh_1$ is the respresentation $\mathbf{3}\boxtimes\mathbf{\check 3}=\End(\mathbf{3})$ of $\SL_3\times \SL_3$, and the determinant is a cubic invariant. 

In case $6_k$, $\fh_1$ is the respresentation $\mathbf{1}\oplus\mathbf{8}$ of $\Spin_7$, 
where $\mathbf{8}$ is the Spin representation, which affords an invariant quadratic form $q$. The map $(x,v)\mapsto x\cdot q(v)$ is a cubic invariant.

{\bf Cases $4_e, 4_d$:\ } These cases have $m=4$ and $\dim\fc=1$ so $\mu_2\leq W_H(\fc,\tau)$, by Lemma \ref{m/2}. We show that in both cases there is a quartic invariant but no quadratic invariant. 

In case $4_e$, $\fh_1$ is the representation $\Lambda^3(\mathbf{6})=\mathbf{6}\oplus\mathbf{14}$ of $\Sp_6\times T_1$, where $t\in T_1$ acts by $t,t^{-1}$ on the respective summands. Since $p>2$ both summands are irreducible so there is no invariant in bidegree $(1,1)$. In characteristic zero one computes that $\Sym^2(\mathbf{6})$ appears in $\Sym^2(\mathbf{14})$, giving a nonzero $H_0$ quartic invariant, which persists in positive characteristic by Lemma \ref{invarianttheory}. 

In case $4_d$, $\fh_1$ is the representation $\mathbf{2}\boxtimes\mathbf{8}$ of $\SL_2\times\Spin_7$. Since this representation is irreducible and symplectic there is no quadratic invariant. To find a quartic invariant we may assume the characteristic of $k$ is zero. Write 
$$\fh_1=\mathbf{8}_+\oplus \mathbf{8}_-,$$
according to the characters $t\mapsto t^{\pm 1}$ of the maximal torus of $\SL_2$. 
One checks that 
$$\dim\left[\Sym^{4-i}(\mathbf{8_+})\otimes \Sym^i(\mathbf{8_-})\right]^{\Spin_7}=
\begin{cases}
1&\quad\text{for\ $i\neq 2$}\\
2&\quad\text{for\ $i=2$}.
\end{cases}
$$
Since this summand affords the character $t^{4-2i}$ of the maximal torus of $SL_2$, 
it follows that there is a one-dimensional space of quartic invariants in $\fh_1$ for $\SL_2\times\Spin_7$.

\subsection{Standard subalgebras and Kostant sections}\label{E6outerKostant}

Fix a torsion automorphism $\theta=\Ad(s)\vt$ of $\fh=\fe_6$, 
with $s\in S=(T^\vt)^\circ$, 
and let $\tau\in\Aut(\fh)$ be another torsion automorphism of the form 
$\tau=\Ad(t)$  (inner case)  or $\tau=\Ad(t)\vt$ (outer case), for some $t\in S$. 
We call the fixed-point subalgebra $\fh^\tau$ a {\bf standard subalgebra}. 
The standard subalgebras $\fh^\tau$ for inner automorphisms $\tau=\Ad(t)$
are in bijection with proper subdiagrams of the affine diagram of type $E_6$; 
these subalgebras all contain $\ft$ as a Cartan subalgebra. The standard subalgebras $\fh^\tau$ for outer automorphisms $\tau=\Ad(t)\vt$
are in bijection with proper subdiagrams of the affine diagram of type ${^2E_6}$; 
these subalgebras all contain $\ft^\vt$ as a Cartan subalgebra. 

The automorphisms $\theta$ and $\tau$ commute, 
so  $\theta$ acts on the standard subalgebra 
$\fk:=\fg^\tau$. If $\tau$ is inner and $\vt$ acts nontrivially on the subdiagram for $\fk$ then $\theta\vert_\fk$ is outer, because $\theta$ 
permutes a basis of the root-system of $\ft$ in $\fk$. And if $\tau$ is outer then 
$\theta\vert_\fk$ must be inner, because $\theta$ acts trivially on the Cartan subalgebra $\ft^\vt$ of $\fk$. 

Suppose now that $\rank(\theta\vert_\fk)=\rank(\theta)$, so that there is a Cartan subspace $\fc$ for $\theta$ such that $\fc\subset\fk$. Let $K=\Aut(\fk)^\circ$ and let $\widetilde K$ be the connected subgroup of $H$ corresponding to $\fk$. These groups are normalized by $\theta$ and the natural map $\widetilde K\to K$ restricts to a surjection 
$$\widetilde K_0:=(\widetilde K^\theta)^\circ\lra 
(K^\theta)^\circ=: K_0
$$
which induces an isomorphism 
$$N_{\widetilde K_0}(\fc)/Z_{\widetilde K_0}(\fc)\simeq 
N_{K_0}(\fc)/Z_{K_0}(\fc).
$$
It follows that we have an embedding of little Weyl groups
$$W_K(\fc,\theta\vert_\fk)\hra W_H(\fc,\theta).$$

With the exception of number $2_c$, the next-to-right-most column of Table 23 below gives 
the Kac diagram of an 
$H$-conjugate $\theta'$ of $\theta$ such that the subdiagram of $1's$ determines a standard subalgebra $\fk$ (given in the last column) such that
$$\rank(\theta\vert_\fk)=\rank(\theta)\qquad \text{and}\qquad 
W_K(\fc,\theta\vert_\fk)=W_H(\fc,\theta),
$$
and such that $\theta\vert_\fk$ satisfies the conditions of Lemma \ref{vinberg:stable}. 
From \cite[Prop. 5.2]{levy:thetap} it follows that $\theta$ admits a Kostant section contained in $\fk$. 

In the table below we  indicate $\fk=\fh^\tau$ as the subdiagram of $1's$ in a Kac diagram of type $E_6$ or ${^2E_6}$ according to whether $\tau$ is inner or outer. Recall that $\theta\vert_\fk$ is then outer or inner, respectively. The superscript ${^2X}$ means that $\theta\vert_\fk$ is outer. The notation ${^2(2A_2)}$ indicates that $\fk\simeq \fsl_3\oplus\fsl_3$ and $\theta$ swaps the two factors.

In the exceptional case $2_c$, previous work on involutions \cite[Prop. 23]{kostant-rallis} 
(for $k=\bc$) and \cite[6.3]{levy:involutions} (for $p\neq 2$) shows that there is a $\theta$-stable subalgebra $\fk\simeq \fsl_3$ containing $\fc$ as a Cartan subalgebra, 
and $W_H(\fc,\theta)$ is just the ordinary Weyl group of $\fc$ in $\fk$. In this case $\theta$ is the unique (up to conjugacy) pinned involution of $\fsl_3$, which is known to have a Kostant section.  

\begin{center}
{\small Table 23: The gradings of positive rank in type $E_6$ (outer case)}
$$
{\renewcommand{\arraystretch}{1.3}
\begin{array}{cccccccc}
\hline
\text{ No.} &\theta\vert_\fh& w\!\in\!-W(E_6)& w\in W(E_7)& 
W_H(\fc,\theta\vert_\fh)&\text{degrees}&\theta'\vert_\fh&\fk\\
 \hline
18_a&\outEVI{1}{1}{1}{1}{1}  &-E_6(a_1)&E_7&\bfmu_{9}&9&
\outEVI{1}{1}{1}{1}{1}&{^2E_6}\\
 
12_b
&\outEVI{1}{1}{0}{1}{1}  &-E_6&E_7(a_2)&\bfmu_{12}&12&
\outEVI{-\!2}{1}{1}{1}{1}&{^2E_6}\\
 
10_b&\outEVI{1}{1}{0}{1}{0}  &-(A_4+A_1)&D_6&\bfmu_{5}&5&
\outEVI{-\!3}{1}{1}{1}{1}&{^2E_6}\\

10_a&\outEVI{1}{0}{1}{0}{1}  &-(A_4+A_1)&D_6&\bfmu_{5}&5&
\outEVI{-1}{1}{1}{1}{-\!\!1}&{^2A_5}\\

10_c&\outEVI{0}{1}{0}{1}{1} &-(A_4+A_1)&D_6&\bfmu_{5}&5&
\outEVI{9}{-\!5}{\ 1}{1}{1}&{^2D_5}\\

8_f&\outEVI{1}{0}{0}{1}{1} &-D_5&D_5+A_1&\bfmu_{8}&8&
\outEVI{-\!4}{1}{1}{1}{1}&{^2E_6}\\

8_c&\outEVI{0}{1}{0}{1}{0}  &-D_5&D_5+A_1&\bfmu_{8}&8&
\outEVI{1}{1}{1}{1}{-\!4}&C_4\\

6_a&\outEVI{1}{0}{0}{1}{0}  &-(3A_2)&E_7(a_4)&G_{25}&6,9,12&
\outEVI{-\!5}{1}{1}{1}{1}&{^2E_6}\\
 
6_c&\outEVI{0}{1}{0}{0}{1} &-(2A_2)&D_6(a_2)+A_1&G(3,1,2)&3,6&
\outEVI{2}{1}{1}{1}{-6}&{^2A_5}\\
 
6_g&\outEVI{0}{0}{0}{1}{1}  &-A_2&D_4+3A_1&\bfmu_{6}&6&
\outEVI{-3}{\ 0}{\ 1}{1}{1}&B_3\\
 
6_i&\outEVI{0}{0}{1}{0}{0}  &-(A_5+A_1)&A_5'&\bfmu_{3}&3&
\outEVI{0}{1}{1}{0}{-\!2}&{^2(2A_2)}\\
 
6_k
&\outEVI{1}{1}{0}{0}{0}  &-(A_5+A_1)&A_5'&\bfmu_{3}&3&
\outEVI{0}{1}{1}{2}{-\!6}&{^2(2A_2)}\\
 
4_b&\outEVI{0}{0}{0}{1}{0} &-D_4(a_1)&2A_3+A_1&G_8&8,12&
\outEVI{-\!6}{1}{1}{1}{1}&{^2E_6}\\
 
4_d&\outEVI{0}{1}{0}{0}{0} &-(A_3+2A_1)&(A_3+A_1)''&\bfmu_{4}&4&
\outEVI{1}{1}{1}{-2}{0}&A_3\\

4_e&\outEVI{1}{0}{0}{0}{1}  &-(A_3+2A_1)&(A_3+A_1)''&\bfmu_{4}&4&
\outEVI{-1}{-\!\!1}{1}{1}{0}&B_2\\
 
2_a&\outEVI{0}{0}{0}{0}{1} &-1&7A_1&W(E_6)&2,5,6,8,9,12&
\outEVI{-\!7}{1}{1}{1}{1}&{^2E_6}\\
 
2_c&\outEVI{1}{0}{0}{0}{0}  &-(4A_1)&(3A_1)'&W(A_2)&2,3&
---&{^2A_2}\\
 
\hline

\end{array}}
$$
\end{center}

\def\noopsort#1{}
\providecommand{\bysame}{\leavevmode\hbox to3em{\hrulefill}\thinspace}

\end{document}